\documentclass[12pt]{amsart}

\usepackage{amssymb,amsmath}
\usepackage{graphpap}
\usepackage{supertabular}
\usepackage{epic}

 
\newtheorem{theorem}[subsection]{Theorem}
\newtheorem*{theorem*}{Theorem}
\newtheorem{lemma}[subsection]{Lemma}
\newtheorem{proposition}[subsection]{Proposition}
\newtheorem{corollary}[subsection]{Corollary}

\newtheorem*{theorem1}{Main theorem}
\newtheorem*{proposition1}{Proposition}

\theoremstyle{definition}
\newtheorem{definition}[subsection]{Definition}
\newtheorem{problem}[subsection]{Problem}
\newtheorem{example}[subsection]{Example}
\newtheorem{examples}[subsection]{Examples}

\theoremstyle{remark}
\newtheorem{remark}[subsection]{Remark}

\newcommand{\mt}[1]{\operatorname{#1}}
\newcommand{\EEE}{{\mathbb E}}
\newcommand{\DDD}{{\mathbb D}}
\newcommand{\AAA}{{\mathbb A}}
\newcommand{\QQ}{{\mathbb Q}}
\newcommand{\ZZ}{{\mathbb Z}}
\newcommand{\CC}{{\mathbb C}}
\newcommand{\OO}{{\mathcal O}}
\newcommand{\RR}{{\mathbb R}}
\newcommand{\PP}{{\mathbb P}}
\newcommand{\FF}{{\mathcal F}}
\newcommand{\NN}{{\mathbb N}}
\newcommand{\FFF}{{\mathbb F}}
\newcommand{\Supp}{\mt{Supp}}
\newcommand{\Sing}{\mt{Sing}}
\newcommand{\Diff}{\mt{Diff}}
\newcommand{\ord}{\mt{ord}}
\newcommand{\wt}{\mt{wt}}
\newcommand{\D}{{\Delta}}
\newcommand{\Hom}{\mt{Hom}}
\newcommand{\Exc}{\mt{Exc}}
\newcommand{\codim}{\mt{codim}}
\newcommand{\mult}{\mt{mult}}

\newcommand{\down}[1]{\llcorner #1 \lrcorner}
\newcommand{\up}[1]{\ulcorner #1 \urcorner}
\newcommand{\fr}[1]{\{ #1\}}

\title{Classification of three-dimensional exceptional log canonical
hypersurface singularities \large{I}}
\author{S.~A.~Kudryavtsev}

\date{}
\address{Department of Algebra, Faculty of Mathematics,
Moscow State Lomonosov University, 117234 Moscow,
Russia}
\email{kudryav@mech.math.msu.su}

\begin{document}
\begin{abstract}
In this paper the three-dimensional exceptional strictly log canonical
hypersurface singularities
are described and the detailed classification of
three-dimensional exceptional canonical
hypersurface singularities is given under the condition of well-formedness.
\end{abstract}
\maketitle

\section*{\bf {Introduction}}
One of the main problems in log Minimal Model Program is the study of
appeared extremal contractions (singularities), i.e.
contractions $f\colon(X,D_X)\to (X',D_{X'})$, where
$-(K_X+D_X)$ is  $f$-ample divisor and $\rho(X/X')=1$.
For the three-dimensional varieties this problem is completely open except the
classification of Fano manifolds obtained by
V.~A.~Iskovskikh at the end of 70th. In solving this problem
one of the main difficulties is the absence of geometrical description of
singularities. The first step in Fano manifold classification is to find a
"good"\ divisor in the anticanonical linear system $|-K_X|$. The existence of a
"good"\ divisor for extremal contraction allows to understand the contraction
structure (example \ref{exam1}). Therefore the main problem discussed in this
paper can be formulated as the follows:

\begin{problem}\label{problem}
To find a "good"\ element in the multiple anticanonical linear system for a
variety, extremal contraction or singularity.
\end{problem}

\begin{example}\label{exam1}
\begin{enumerate}
\item Consider a small extremal contraction of three-dimensional terminal
variety. Then the existence of a divisor with Du Val singularities in the linear
system $|-2K_X|$ implies the flip existence \cite{Kaw1}.
\item If an extremal contraction of three-dimensional terminal variety is a conic
bundle then the existence of a divisor with Du Val singularities in the
anticanonical linear system $|-K_X|$ allows at once
to obtain the complete local classification \cite{Pr1}.
\end{enumerate}
\end{example}
The presence of a "good"\ element is closely connected with the
other problems, for
example with the automorphisms of
$K3$ surfaces \cite[1.12]{Sh2} and with the gorenstein indexes of strictly log
canonical singularities (see example \ref{lc}).
\par
Recently V.~V.~Shokurov gave the approach to solve the problem
\ref{problem} in \cite{Sh1},\cite{Sh2}.
It consists of two parts:
\begin{enumerate}
\item The first idea is an inductive transfer from $n$-dimensional varieties
to $(n-1)$-dimensional ones.
\item The second idea is to separate the varieties to the exceptional and
nonexceptional ones.
\end{enumerate}

\par
Shortly, the first step is the following. First of all one have to find some
divisor $S$ on $X$ or on some blow-up of $X$ such that
pair $(S,\Diff_S(D_X))$ is klt and
$-(K_S+\Diff_S(D_X))$ is ample. Then a "good"\ element of this pair can be
extended to a "good"\ divisor of $K_X$ on $X$.

\par
In studying varieties, extremal contractions, singularities the exceptional
phenomenon study importance follows from the following observation:
\begin{enumerate}
\item If a variety or extremal contraction or singularity is nonexceptional
then the linear system $|-nK_X|$ must have a "good"\ member for small $n$.
For example, we can take $n\in \{1,2\}$ for the two-dimensional singularities
\cite[5.2]{Sh1} and $n\in \{1,2,3,4,6\}$ for the three-dimensional
singularities \cite[7.1]{Sh2}.
\item The exceptional singularities are "bounded"\ and can be classified.
For example, in this paper it will be checked that the
three-dimensional exceptional
hypersurface singularities have the finite number of types.
\end{enumerate}

\par
Generally speaking,
the {\it regular} (nonexceptional) extremal contractions (singularities) can not
be completely classified even in the dimension three (for example,
in \cite{Mar} it was shown that it is impossible to enumerate all normal
form equations of the three-dimensional terminal hypersurface singularities,
which are
the very simple three-dimensional singularities).
On the other hand, the existence of small index complement allows to separate the
extremal contractions (singularities) to the families with common properties.

\par
The exceptional extremal contractions (singularities) can be completely
classified but they can have the large minimal complement index.
\par
In \cite{Pr2} it was proved that the exceptionality remains valid after
an inductive transfer and conversely, if
$(S,\Diff_S(D_X))$ is exceptional then $(X/Z, D_X)$ is exceptional too.
Therefore it became possible to speak about {\it the inductive method of the
algebraic variety classification}.

\par
This paper is devoted to the application of the inductive method
to the hypersurface
singularities. In this paper the three-dimensional exceptional hypersurface
singularities are classified. For this purpose, at first the corresponding
quasihomogeneous parts are described, then the purely log terminal blow-ups are
constructed. Next the log Del Pezzo surfaces obtained are investigated on
the exceptionality. In particular, the minimal index of complements is found for
every singularity.
\par
Also the description of three-dimensional exceptional strictly log canonical
hypersurface singularities is given.

\begin{theorem1}
Let $(X,0) \subset(\CC^4,0)$ be a three-dimensional exceptional canonical
$($respectively strictly log canonical$)$ hypersurface singularity
defined by a polynomial $f$.
Then there exists a biholomorphic coordinate change
$\psi \colon (\CC^4,0) \to (\CC^4_{t,z,x,y},0)$ and unique primitive vector
${\bf p}\in N_{\RR}$ such that just one of the following two possibilities holds:
\begin{enumerate}
\item The quasihomogeneous polynomial $\tilde f_{\bf p}=(f\circ \psi)_{\bf p}$
defines an exceptional canonical $($respectively strictly log canonical
and canonical outside
{\bf 0}$)$ singularity $(X_{\bf p},0) \subset(\CC^4_{t,z,x,y},0)$.
In this case ${\bf p}$-blow-up of $\CC^4$ induces purely log terminal blow-ups
$\varphi \colon (Y,E)\to (X,0)$ and
$\varphi_{\bf p} \colon (Y_{\bf p},E_{\bf p})\to (X_{\bf p},0)$, where
$(E,\Diff_E(0))=(E_{\bf p},\Diff_{E_{\bf p}}(0))$.
That is, these singularities have the same type and in particular the same
complement index.
\par
The canonical singularities satisfying the condition of well-formedness ---
$\Diff_{E/\PP(\bf p)}(0)=0$ are classified in the theorems
\ref{f5}, \ref{f4} and in the tables of chapter \S 4. The polynomial
$\tilde f_{\bf p}$; $(E,\Diff_E(0))$; minimal complement index are written
in the tables.
\par
The strictly log canonical and canonical outside
{\bf 0} quasihomogeneous singularities are always exceptional
$($in any dimension$)$ by theorem \ref{crplt}.
In the three-dimensional case
$(E,\Diff_E(0))=(\tilde f_{\bf p}\subset\PP({\bf p}),0)$ is
$K3$ surface with Du Val singularities and $(X,0)$ is
1-complementary.
\item $\tilde f_{\bf p}=t^3+g_2^2(z,x,y)$, where $g_2$ is an irreducible
homogeneous polynomial of degree two.
In this case the purely log terminal blow-ups are constructed in the theorems
\ref{classific2} and \ref{classific3}. Also it was obtained the similar
classification depending on the type of jets
$\tilde f_5$ and $\tilde f_6$.
\end{enumerate}
\end{theorem1}
For the non-degenerate strictly log canonical singularities it was proved that
they are exceptional if and only if they have purely elliptic type
$(0,d-1)$, where $d$ is a dimension of the singularity \cite{Ishii}.
Also for every non-degenerate strictly log canonical singularity the
quasihomogeneous log canonical part can be found \cite[3.5]{PrI}.
Therefore, in our terms it is exceptional if and only if
its quasihomogeneous part defines a canonical outside
{\bf 0} singularity (see corollary \ref{main}).
Using this fact it is not difficult to get the list of all three-dimensional
exceptional strictly log canonical hypersurface singularities.
The mentioned one is not given in this paper. For the isolated singularities
such list is given in \cite{Yon}. This list corresponds to
95 families of weighted $K3$ surfaces obtained by M.Reid and I.Fletcher in
\cite[4.5]{Reid1}, \cite{Fletcher}.
Indeed, the single exceptional divisor
of purely log terminal blow-up is $K3$ surface.
\par
This paper is organized in the following way.
In chapter \S 1 the main definitions and
preliminary results are collected. In chapter \S 2 the quasihomogeneous
singularities are studied. The main theorem is proved in chapter \S 3.
Also the explanation of three-dimensional canonical singularities
study on the exceptionality is given.
In chapter \S 4 the summarizing tables are written.
\par
I am grateful to Professor V.A.Iskovskikh, Professor Yu.G.Prokhorov and
Professor V.V.Shokurov for
useful discussions and valuable remarks.
\par
The research was partially
supported by a grant 99-01-01132 from the Russian Foundation of Basic Research
and a grant INTAS-OPEN 2000\#269.

\section{\bf Preliminary facts and results}
All varieties are algebraic and are assumed to be defined over
$\CC$, the complex number field. We use the terminology, notations
of log Minimal
Model Program and the main properties of complements given in
\cite{KMM}, \cite{Koetal}, \cite{PrLect}.
{\it A strictly lc singularity}
is called lc singularity, but not klt singularity.

\begin{definition}
We say that the boundary $D=\sum d_iD_i$ has {\it standard coefficients} if
$d_i\in \{1\}\cup \big\{1-\frac1m|m\in \NN \big\}$ for all $i$.
\end{definition}

\begin{definition}
Let $(X,D)$ be a log pair, where $D$ is a subboundary. Then a
\textit{$\QQ$-complement} of $K_X+D$ is a
log divisor $K_X+D'$ such that $D'\ge D$, $K_X+D'$ is lc and
$n(K_X+D')\sim 0$ for some $n\in\NN$.
\end{definition}

\begin{definition}\label{defcompl}
Let $X$ be a normal variety and let $D=S+B$ be a subboundary on
$X$ such that $B$ and $S$ have no common components, $S$ is an
effective integral divisor and $\down{B}\le 0$. Then we say that
$K_X+D$ is \textit{$n$-complementary} if there is a $\QQ$-divisor
$D^+$ such that
\begin{enumerate}
\item
$n(K_X+D^+)\sim 0$ (in particular, $nD^+$ is an integral divisor);
\item
$K_X+D^+$ is lc;
\item
$nD^+\ge nS+\down{(n+1)B}$.
\end{enumerate}
In this situation the
\textit{$n$-complement} of
$K_X+D$ is $K_X+D^+$. The divisor $D^+$ is called
\textit{$n$-complement} too.
\end{definition}

The invariance of the complements with respect to log MMP and their inductive
properties are shown in the following three theorems.

\begin{theorem}\label{down} \cite[5.4]{Sh1}
Let $f\colon X\to Y$ be a birational contraction and let $D$ be a
subboundary on $X$. Assume that $K_X+D$ is $n$-complementary for
some $n\in\NN$. Then $K_Y+f(D)$ is also $n$-complementary.
\end{theorem}

Under additional assumptions we have the inverse statement:
\begin{theorem}\cite[2.13]{Sh2}
Let $f\colon Y\to X$ be a birational contraction
and let $D$ be a subboundary on $Y$ such that
$K_Y+D$ is nef over $X$ and $f_*(D)$ has standard coefficients.
Assume that $K_X+f(D)$ is $n$-complementary. Then $K_Y+D$ is also
$n$-complementary.
\end{theorem}

\begin{theorem} \label{Prodolj}\cite[2.1]{Pr4}
Let $(X/Z\ni P,D=S+B)$ be a log variety. Put $S=\down{D}$ and
$B=\fr{D}$. Assume that
\begin{enumerate}
\item
$K_X+D$ is plt;
\item
$-(K_X+D)$ is nef and big over $Z$;
\item
$S\ne 0$ near $f^{-1}(P)$;
\item
$D$ has standard coefficients.
\end{enumerate}
Further, assume that near $f^{-1}(P)\cap S$ there exists an
$n$-complement $K_S+\Diff_S(B)^+$ of $K_S+\Diff_S(B)$. Then near
$f^{-1}(P)$ there exists an $n$-complement $K_X+S+B^+$ of
$K_X+S+B$ such that $\Diff_S(B)^+=\Diff_S(B^+)$.
\end{theorem}

\begin{definition}
Let $(X/Z\ni P,D)$ be a contraction of varieties, where $D$ is a boundary.
\begin{enumerate}
\item Assume that $Z$ is not a point (local case). Then $(X/Z\ni P,D)$
is said to be \textit{exceptional}
over $P$ if for any $\QQ$-complement $K_X+D'$ of $K_X+D$ near the
fiber over $P$ there exists at most one (not necessarily
exceptional) divisor $E$ such that $a(E,D')=-1$.
\item Assume that $Z$ is a point (global case). Then $(X,D)$ is said to
be \textit{exceptional} if every $\QQ$-complement of $K_X+D$ is
klt.
\end{enumerate}
If $Z=X$ then $(X/X\ni P)$ is a singularity $(X\ni P)$.
\end{definition}

\begin{example}\label{lc}
(1) Let $(X \ni P)$ be a singularity. Suppose that there is an effective divisor
$H$ such that $(X,H)$ is lc and $\down{H} \ne 0$. Then the singularity is not
exceptional. Therefore the three dimensional terminal singularity
is not exceptional
because there is a divisor having only Du Val singularities in the anticanonical
linear system $|-K_X|$ \cite[6.4]{YPG}.
\par
It should be noted that there are the
examples of four-dimensional exceptional terminal singularities \cite{Kud1}.
\par
(2) Let $(X\ni P)$ be a strictly log canonical exceptional singularity. Then
the minimal index of complements is equal to the gorenstein index of
$(X\ni P)$ \cite[1.10]{Kud2}.
The minimal index of complements is bounded for the three dimensional
log canonical singularities \cite[7.1]{Sh1}.
A hypothesis is that this index is not
more then 66.
For the strictly log canonical exceptional singularities it was proved in
\cite{I} and \cite{Fuj}.
For the nonexceptional nonisolated strictly log canonical
singularities the gorenstein index
is not bounded \cite[5.1]{Fuj}.
\end{example}

\begin{definition} \label{def1}
Let $(X\ni P)$ be a log canonical singularity and
let $f\colon Y \to X$ be a blow-up such that
the exceptional locus of $f$ contains only one irreducible divisor
$E\ (\Exc(f)=E)$. Then $f:(Y,E) \to X$ is called  {\it a purely log terminal
(plt) blow-up} if $K_Y+E$ is plt and
$-E$ is $f$-ample.
\end{definition}

\begin{definition} \label{deflc}
Let $(X\ni P)$ be a log canonical singularity and
let $f\colon Y \to X$ be a blow-up such that
the exceptional locus of $f$ contains only divisors
$E\ (\Exc(f)=E)$. Then $f:(Y,E) \to X$ is called  {\it a log canonical
(lc) blow-up} if $K_Y+E$ is lc and $-E$ is $f$-ample.
\end{definition}

\begin{remark} In the definitions \ref{def1} and \ref{deflc} it is demanded
that divisor $E$ must be $\QQ$-Cartier. Hence $Y$ is a $\QQ$-gorenstein variety.
\end{remark}

\begin{remark} \label{smain1}
The notations are the same as in the definition \ref{def1}.
We have the following properties
of plt blow-ups:
\begin{enumerate}
\item If $X$ is klt then $-(K_Y+E)$ is $f$-ample and a plt
blow-up always exists \cite[1.5]{Kud2}.
\item If $X$ is strictly lc then $a(E,0)=-1$ and a plt
blow-up always exists for the exceptional singularities \cite[1.9]{Kud2}.
\item If $X$ is $\QQ$-factorial then $Y$ is also
$\QQ$-factorial and $\rho(Y/X)=1$ \cite[2.2]{Pr1}. Hence, in the definition
\ref{def1} for a $\QQ$-factorial singularity it
is not necessary to demand that the divisor $-E$ is $f$-ample
because amplness takes place always.
\item By the inversion of adjunction $K_E+\Diff_E(0)$ is klt.
If $K_E+\Diff_E(0)$ is
$n$-complementary then $K_Y+E$ is $n$-complementary and
$K_X$ is $n$-complementary too (see theorem \ref{Prodolj} and theorem \ref{down}).
\item Let $f_i:(Y_i,E_i)\to (X \ni P)$ be two
plt blow-ups. If $E_1$ and $E_2$ define the same discrete valuation of
function field $\mathcal{K}(X)$ then $f_1$ and $f_2$
are isomorphic \cite[2.2]{Pr1}.
\end{enumerate}
\end{remark}

\begin{definition}\label{defexc}
Let $(X\ni P)$ be a log canonical singularity. It is said to be
{\it weakly exceptional} if there exists only one plt blow-up
(up to isomorphism).
\end{definition}

The following two theorems \ref{CrExc}, \ref{CrWexc} demonstrate the
inductive transfer from an $n$-dimensional exceptional (resp. weakly exceptional)
klt singularity to the unique $(n-1)$-dimensional exceptional
(resp. weakly exceptional) log Fano variety $(E,\Diff_E(0))$
having the same indexes of complements. The inverse statements take place too.

\begin{theorem}\cite[4.9]{Pr2}\label{CrExc}
Let $(X \ni P)$ be a klt singularity and let $f\colon (Y,E) \to X$ be a plt
blow-up of $P$. Then the following conditions are equivalent:
\begin{enumerate}
\item $(X \ni P)$ is an exceptional singularity;
\item $(E,\Diff_E(0))$ is an exceptional log variety.
\end{enumerate}
\end{theorem}

\begin{theorem} \label{CrWexc}\cite[2.1]{Kud2}
Let $(X \ni P)$ be a klt singularity and let $f\colon (Y,E) \to X$ be a plt
blow-up of $P$. Then the following conditions are equivalent:
\begin{enumerate}
\item $(X\ni P)$ is not a weakly exceptional singularity;
\item there exists an effective $\QQ$-divisor $D\ge\Diff_E(0)$ such that
$-(K_E+D)$ is ample and $(E,D)$ is not klt;
\item there exists an effective $\QQ$-divisor $D\ge\Diff_E(0)$ such that
$-(K_E+D)$ is ample and $(E,D)$ is not lc.
\end{enumerate}
\end{theorem}

According to the example of weakly exceptional singularity but not exceptional
singularity in any dimension \cite[2.4]{Kud2} we can construct
the example of the exceptional
singularity.

\begin{example} Let $(X\ni P)$ be an $n$-dimensional canonical hypersurface
singularity given by the equation
$(x^n_1+x^{n+1}_2+x^{n+1}_3+\cdots+x^{n+1}_{n+1}=0) \subset (\CC^{n+1},0)$.
The weighted blow-up of
$\CC^{n+1}$ with weights $(n+1,n,n,\ldots,n)$ induces a plt blow-up of $P$.
Using the basic theorems about the hypersurfaces in the weighted spaces
\cite{Fletcher} (see also \S 4),
it is easy to calculate that the obtained
log Fano variety $(E,\Diff_E(0))$ is
\begin{gather*}
(\tilde x_1+\tilde x^{n+1}_2+\cdots+\tilde x^{n+1}_{n+1}\subset
\PP(n+1,1,\ldots,1),\frac{n-1}{n}\{\tilde x_1=0\})=\\
=(\PP^{n-1},\frac{n-1}{n}Q_{n+1}),
\end{gather*}
where $Q_{n+1}$ is a smooth hypersurface of degree $n+1$ in $\PP^{n-1}$.
Applying the approach of \cite{Kud1} it is not difficult to check that
given log variety is exceptional for $n\ge3$. The
divisor $\frac{n-1}{n}\{\tilde x_1=0\}+\frac{1}{n}\{\tilde x_2=0\}$ is
$n$-complement of minimal index. Our
singularity is exceptional
by exceptionality criterion \ref{CrExc} and $n$-complementary
by remark \ref{smain1}.
\end{example}

\begin{proposition}\label{lcplt}
An exceptional singularity is weakly exceptional. Assume that there exists
lc blow-up but not plt blow-up of $(X\ni P)$. Then $(X\ni P)$ is not exceptional.
\begin{proof} In general case the first statement was proved in
\cite[1.12]{Kud2}. The second statement follows from the proof \cite[2.4]{PrI}.
In fact, let
$f\colon (Y,E) \to (X\ni P)$ be lc blow-up but not plt blow-up.
We may assume that
$(X\ni P)$ is klt otherwise the proof is trivial. Since
$-(K_Y+E)$ is  $f$-ample then the
linear system $|-n(K_Y+E)|$ is very ample over $X$ for $n\gg 0$.
Let $H\in |-n(K_Y+E)|$ be a general element and let
$B=\frac1n f(H)$.
By Bertini's theorem
$K_Y+E+\frac1nH$ is lc but not plt. Since
$K_Y+E+\frac1nH\sim_{\QQ}0$ then the pair $(X,B)$ is lc but not exceptional.
\end{proof}
\end{proposition}

The next proposition follows from the corollary \cite[1.6]{Kud2} and the theorem
\cite[1.9]{Kud2}.
\begin{proposition} \label{videxc}
Let $f\colon (Y,E)\to (X\ni P)$ be a plt blow-up of
log canonical singularity and let
$\dim f(E)\ge 1$. If  $(X\ni P)$ is klt then it is not weakly exceptional.
If $(X\ni P)$ is strictly lc then it is weakly exceptional but not exceptional.
\end{proposition}

\begin{example}\cite[5.2]{Sh1}
A two dimensional klt singularity is exceptional if and only if it has type
$\EEE_6,\EEE_7,\EEE_8$. Its plt blow-up is the extraction of central vertex
of minimal resolution graph. The corresponding one-dimensional log Fano
variety
$(E,\Diff_E(0))$ is
$(\PP^1,\frac{n_1-1}{n_1}P_1+\frac{n_2-1}{n_2}P_2+\frac{n_3-1}{n_3}P_3)$, where
$(n_1,n_2,n_3)=(2,3,3),(2,3,4)$ or $(2,3,5)$ respectively.
\end{example}

\begin{example}\cite[\S 6]{PrLect}
A two dimensional strictly lc  singularity is exceptional if and only if
it has one of the following types: simple elliptic singularity or
singularities of
$\widetilde{\DDD}_4$, $\widetilde{\EEE}_6, \widetilde{\EEE}_7,
\widetilde{\EEE}_8$ types. Its plt blow-up is the extraction of central vertex
of minimal resolution graph. The corresponding one-dimensional log Fano
variety $(E,\Diff_E(0))$ is a smooth
elliptic curve, $(\PP^1,\frac12P_1+\frac12P_2+\frac12P_3+\frac12P_4)$ or
$(\PP^1,\frac{n_1-1}{n_1}P_1+\frac{n_2-1}{n_2}P_2+\frac{n_3-1}{n_3}P_3)$, where
$(n_1,n_2,n_3)=(3,3,3),(2,4,4)$ or $(2,3,6)$ respectively.
\end{example}

\par
Thus, there are only three types of induced log Fano varieties for the
two-dimensional
klt singularities and there are only five types for the two-dimensional
strictly lc singularities. There is a conjecture that the similar result about
the finite number of types is true in any dimension.
For the three-dimensional hypersurface singularities this conjecture is checked in
this paper. The essential difference from the
two-dimensional singularities is that
the infinite number of nonisomorphic
singularities exists in almost every type of them. But they have
the same resolutions.

\subsection{The main facts of toric geometry}\label{Odatoric}
We make use of toric geometry terminology (see \cite{Oda}).
Let $N$ be a free abelian group $\ZZ^n$ and $M$ is its
dual $\Hom_{\ZZ}(N, \ZZ)$. Denote $N\otimes_{\ZZ}\RR$ and
$M\otimes_{\ZZ}\RR$ by $N_{\RR}$ and $M_{\RR}$ respectively.
We have a canonical pairing $\langle\,\ \rangle \colon
N_{\RR}\times M_{\RR}\to \RR$. Let $\sigma$ be the positive
quadrant $\RR_{\ge 0}^n$ of $N_{\RR}= \RR^n$ and $\sigma^{\vee}$
its dual. Then $\CC^n$ is the toric variety corresponding to the
cone $ \sigma$. For a fan $\D$ in $N$ the corresponding toric
variety is denoted by $T_{N}(\D)$.
For a primitive element ${\bf p}\in N$ of a $1$-dimensional cone
$\tau=\RR_{\ge 0}{\bf p}$ in $\D$, the closure
$\overline{orb(\RR_{\ge 0}{\bf p})}$ is denoted by $D_{\bf p}$,
where $D_{\bf p}$ is a
divisor on $T_{N}(\D)$.

\par
Let $(X,0) \subset (\CC^n, 0)$ be a hypersurface singularity and
${\bf p}=(p_{1},\ldots p_n)$ a primitive element in $N$ with $p_{i_1}>0$ and
$p_{i_2}>0$
for some $i_1\ne i_2$ (such ${\bf p}$ is called a \textit{weight}).
Let $\varphi\colon \CC^n({\bf p})\to \CC^n$ be the blow-up with a
weight ${\bf p}$. Denote the proper transform of $X$ on $\CC^n({\bf p})$ by
$X({\bf p})$.
The blow-up $\varphi\colon \CC^n({\bf p})\to \CC^n$ and its
restriction $\varphi|_{X({\bf p})} \colon X({\bf p})\to X$ are called
the ${\bf p}$-{\it blow-ups}.
The exceptional locus of $\varphi|_{X({\bf p})}$ is denoted by
$E_{\bf p}$.
The exceptional locus $D_{\bf p}$ of
$\varphi$ is the weighted projective space
$\PP({\bf p})$.
These ${\bf p}$-blow-ups are obtained by a
subdivision of the cone $\sigma$. The corresponding fan consists
of the faces of cones $\sigma_{i}$ $(i=1,\ldots,n)$, where
$\sigma_{i}$ is generated by
${\bf e}_1,\ldots,{\bf e}_{i-1},{\bf p},{\bf e}_{i+1},\ldots,{\bf e}_n$.
Here ${\bf e}_1,\ldots,{\bf e}_n$ are the unit vectors
$(1,0,\ldots,0),\ldots,(0,\ldots,0,1)$ which generate $\sigma$.
\par
A monomial $x_{1}^{m_{1}}\cdots x_{n}^{m_{n}}\in
\CC[[x_{1},x_{2},\ldots, x_{n}]]$ is denoted by $x^m$, where
$m=(m_{1},\cdots,m_{n})\in \ZZ^n =M$. For a power series
$f=\sum_{m}a_{m}x^m \in \CC[[x_{1},x_{2},\ldots, x_{n}]]$ we
write $x^m\in f$, if $a_{m}\ne 0$.
For ${\bf p}\in N_{\RR}$ and a power series $f$ we define
$$\mbox{\bf p}(f) =
\min_{x^m\in f}\langle\mbox{\bf p},m\rangle.$$

We denote the {\it leading term}
$\sum_{\langle{\bf p},m\rangle={\bf p}(f)}a_{m}x^{m}$ of $f$
with respect to ${\bf p}$ by $f_{\bf p}$.

\begin{definition}
For a power series $f=\sum_{m}a_{m}x^m \in
\CC[[x_{1},x_{2},\ldots, x_{n}]]$ define \textit{Newton
polyhedron} $\Gamma_{+}(f)$ in $M_{\RR} $ as follows:
\[ \Gamma_{+}(f) =
\text{the convex hull of}\ \bigcup_{x^m\in f}(m+\sigma^{\vee}). \]
The set of the interior points of $\Gamma_{+}(f)$ is denoted by
$\Gamma_{+}(f)^0$.
For each face $\gamma$ of $\Gamma_{+}(f)$ we define the
polynomial $f_{\gamma}$ as follows:
\[ f_{\gamma}=\sum_{m\in\gamma}a_{m}x^m. \]
A power series $f$ is said to be \textit{non-degenerate} if for
every face $\gamma$ the equation $f_{\gamma}=0$ defines a
smooth hypersurface in the complement of the hypersurface
$x_{1}\cdots x_{n}=0$.
\end{definition}

\begin{theorem}\label{ness} \cite{Mar}
Let $(X,0)\subset ( \CC^n,0)$ be a normal hypersurface singularity
defined by a power series $f$ and $\Gamma_{+}(f)$
is its Newton polyhedron. Let
$(X,0)$ be canonical (resp. log canonical). Then
${\bf 1}=(1, 1, \ldots, 1)\in \Gamma_{+}(f)^0$ $($resp. ${\bf 1}=(1, 1,
\ldots, 1)\in \Gamma_{+}(f)$$)$.
The inverse statements take place for the non-degenerate singularities.
\end{theorem}

Now as a preliminary we discuss the {\it nonisolated hypersurface singularities
study}.
\begin{definition} \label{notisol}
Let $(X,0)\subset (\CC^n,0)$ be a nonisolated hypersurface singularity
(maybe nonnormal) defined by a polynomial $f$.
By proposition \ref{videxc} we can investigate this singularity on the
exceptionality in the arbitrary small neighborhood
$\mathcal U$ of the origin. For example, in the
three-dimensional case it means that
$(\Sing X\setminus {\bf 0})\cap \mathcal U \cong \bigsqcup_i\CC^1_i$ is not
arcwise connected union. The variety
$X\cap \mathcal U$ is isomorphic to
$\CC^1 \times$(two dimensional log canonical hypersurface singularity) along
the every component.
In an $n$-dimensional case let
$(\Sing X\setminus {\bf 0})\cap \mathcal U \cong
\bigcup^m_{i=1}(\mathcal A^{n_i}_i\setminus {\bf 0}),$
where $\mathcal A^{n_i}_i$ are the irreducible affine varieties of the dimension
$n_i$. Now we define the mutually different affine varieties
$\{V^{n_j}_j\}$, where $n_j=\dim V^{n_j}_j$. The system of indexes
$I_j\subset \{1,\ldots,m \}$ is called {\it admissible} if
$$\emptyset \ne \bigcap_{i\in I_j}(\mathcal A^{n_i}_i\setminus {\bf 0})\nsubseteq
\mathcal A^{n_k}_k,\  \text{where}\  k \notin I_j. $$
For all admissible systems $I_j$ we can put
$V^{n_j}_j \stackrel{\rm def}{=}\bigcap_{i\in I_j} \mathcal A^{n_i}_i $.
Let $\mathfrak a_j$ be the ideal radical of the variety $V^{n_j}_j$ and let
$f_1,\ldots,f_r$ be the generators of
$\mathfrak a_j$. Then we have the decomposition
$$
f=\Phi_1(f_1,\ldots,f_r)u_1(x_1,\ldots,x_n)+\ldots+
\Phi_m(f_1,\ldots,f_r)u_m(x_1,\ldots,x_n),
$$
where
$u_i \notin \mathfrak a_j$ for all $i$.
The variety $\{u_1=\ldots=u_m=0\}$ is empty or nonsingular by the construction
of $V^{n_j}_j$. Let $W_k=V_j\cap \{u_k=0 \}$. Repeat the similar process for the
irreducible components of the nonempty  sets $W_k$.
Thus we obtain the sets
$W_{k,j}^{n_{k,j}}\subset V^{n_j}_j$, where $n_k$ are their dimension.
There is a map $\psi \colon X\cap \mathcal U\to V^{n_j}_j$ for all $j$. By the
construction, $(n-1-n_j)$-dimensional singularities
$\psi^{-1}(P_1)$ and $\psi^{-1}(P_2)$ have the same type
(resolution) with respect to Newton diagram, where
$P_1,P_2 \in (V^{n_j}_j\setminus (\bigcup_k W^{n_{k,j}}_{k,j}
\bigcup {\bf 0})) \setminus (\bigcup_{i\ne j} V^{n_i}_i )$. It should be noted
that they can be not biholomorphic.
This type is denoted by $\FF_j$.
The similar statement is true for the map
$\psi \colon X\cap \mathcal U\to W^{n_{k,j}}_{k,j}$ and the set
$(W^{n_{k,j}}_{k,j}\setminus {\bf 0})
\setminus \bigcup_{l\ne k} W^{n_{l,j}}_{l,j}$. The corresponding type is denoted
by $\FF_{j,k}$.
\end{definition}

\begin{examples}\label{first}
1). Let us consider the singularity $(x_1^2+x_2^4+(x_3^2+x_4^3)^2x_4=0, 0)
\subset (\CC^5,0)$. Then $V^2_1=\{x_1=x_2=x_3^2+x_4^3=0 \}$,
$W^1_{1,1}=V^2_1\cap \{x_4=0 \}$, $\FF_1=\{(y_1^2+y_2^4+y_3^2=0,0)\subset
(\CC^3,0)\}$, $\FF_{1,1}=\{(y_1^2+y_2^4+y_3^2y_4=0,0)\subset
(\CC^4,0)\}$.\\
2). Let us consider the singularity $(x_1^2+x_2^3+(x_3^5+x_4^5)x_5^2+
x_3^3x_4^3x_5=0, 0) \subset
(\CC^5,0)$. Then $V^1_1=\{x_1=x_2=x_3=x_4=0 \}$, $V^1_2=\{x_1=x_2=x_3=x_5=0 \}$,
$V^1_3=\{x_1=x_2=x_4=x_5=0 \}$,
$\FF_1=\{(y_1^2+y_2^3+y_3^5+y_4^5+cy_3^3y_4^3=0,0)\subset (\CC^4,0)\}$ and for
the different values of $c$ the singularities can be not biholomorphic,
$\FF_2=\FF_3=\{(y_1^2+y_2^3+y_4^2+y_3^3y_4=0,0)\subset(\CC^4,0)\}$.
\end{examples}

\begin{proposition}\label{equiv}
Let $(X,0)\subset (\CC^n,0)$ be a normal nonisolated hypersurface singularity.
It is terminal outside {\bf 0} $($resp. canonical outside {\bf 0},
log canonical outside {\bf 0}$)$ if and only if $\FF_j$ and $\FF_{j,k}$ are
terminal $($resp. canonical, log canonical$)$ for all $j,k$.
\begin{proof} For instance, let us prove this theorem in the terminal case.
In definition \ref{notisol} we made the subdivision of
$\Sing X \backslash {\bf 0}$ to the "cells". Every cell has a dimension
not less then
one. Along every cell
$X$ has the same singularity type.
Let some $\FF_j$ or $\FF_{j,k}$ be not terminal singularity. Then there exists
a blow-up with a center along the corresponding cell and we obtain an exceptional
divisor with discrepancy $\le 0$. Contradiction.
\par
Conversely assume that
$\FF_j$ and $\FF_{j,k}$ are the terminal singularities and $(X,0)$ is not
terminal outside {\bf 0}. Consider the log-resolution $\psi$.
Then there exists an exceptional divisor
$E$ with discrepancy $a(E,0)\le 0$ such that $\dim\psi(E)\ge 1$. The generic
point $P\in \psi(E)$ is lying in the fixed cell or
$P\notin \Sing X$. In both cases we have a contradiction with $a(E,0)\le 0$.
\end{proof}
\end{proposition}

\begin{definition}\label{onedim}
Let $(X,0)\subset (\CC^n,0)$ be a nonnormal hypersurface singularity. It is
called {\it a log canonical singularity} if $(\CC^n,X)$
is lc.
It is called {\it log canonical outside {\bf 0} singularity} if
$(\CC^n,X)$ is log canonical outside {\bf 0}. The last definition is equivalent
to the conditions that $\FF_j$ and $\FF_{j,k}$ are lc.
\end{definition}

\section{\bf Quasihomogeneous hypersurface singularities}
\begin{definition}
Let $f\colon (\CC^n,0) \to (\CC,0)$ be a polynomial. It is called a
{\it quasihomogeneous polynomial of
 degree $d$} if there exist the strictly positive rational numbers
$p_1, \ldots,p_n$ (which are called {\it the weights}) that

\begin{equation}\label{eq1}
 f(\lambda^{p_1}x_1,\ldots,\lambda^{p_n}x_n)=\lambda^df(x_1,\ldots,x_n)
\end{equation}
for any $\lambda$ and all $x_1,\ldots,x_n$.
\par
For the quasihomogeneous singularities (polynomials) the nonnegative
rational numbers
$p_1, \ldots,p_n$ satisfying condition
(\ref{eq1}) are also called  {\it the weights} and $d\ne 0$ is called
{\it a quasihomogeneous degree} too.
\end{definition}

{\bf Convention.} Without loss of generality it can be assumed that
$p_1, \ldots,p_n, d$ are the integer numbers in this chapter.

\begin{lemma}\label{Laz}
Let $(X,0)\subset (\CC^n,0)$ be a quasihomogeneous hypersurface singularity
defined by a polynomial $f$. If
$\frac{\partial f}{\partial x_1}({\bf x})=\cdots=
\frac{\partial f}{\partial x_n}({\bf x})=0$
then $f({\bf x})=0$.
\begin{proof}
The polynomial $f$ is quasihomogeneous if and only if
$f\in I_{\nabla f}$, where
$I_{\nabla f}=(\frac{\partial f}{\partial x_1},\ldots,
\frac{\partial f}{\partial x_n})$ is a gradient ideal
\cite{Lazzeri}. This concludes the proof.
\end{proof}
\end{lemma}

\begin{proposition}\label{out}
Let $(X,0)\subset (\CC^n,0)$ be a normal quasihomogeneous hypersurface
singularity with weights
$p_1,\ldots,p_n$ and degree $d$ defined by a polynomial
$f$. Then
\begin{enumerate}
\item If $\sum p_i\ge d+1$ and $(X,0)$ is canonical outside {\bf 0} then
${\bf 1} \in \Gamma_{+}(f)^0$.
\item If $\sum p_i\ge d$ and $(X,0)$ is log canonical outside  {\bf 0} then
${\bf 1} \in \Gamma_{+}(f)$.
\end{enumerate}
\begin{proof} Assume the converse. Then
${\bf 1} \notin \Gamma_{+}(f)^0$ or
${\bf 1} \notin \Gamma_{+}(f)$ respectively. Then there exists a primitive vector
${\bf p'}$ from
$\sigma^{\vee}$ such that the plane
$\Pi=\{l \in N_{\RR}|\ \langle {\bf p'},l \rangle = {\bf p'}(f)\}$ is the plane of
support for Newton polyhedron. This plane contains
{\bf 1} or lies over {\bf 1}.
Let $\varphi \colon \CC^n({\bf p'}) \to \CC^n$ be a
${\bf p'}$-blow-up. Since the plane of support $\Pi$ is parallel to some
coordinate axis then at least one coordinate of
${\bf p'}$ is equal to zero. Therefore
$\dim \varphi(D_{\bf p'}) \ge 1$.
Also the discrepancy
$a(D_{\bf p'},X) =\langle{\bf p'},1\rangle - {\bf p'}(f)-1$ is less then
--1 in the
first case and strictly less then
--1 in the second case. Since the singularity is canonical outside
{\bf 0} or lc outside {\bf 0} we obtain the contradiction respectively.
\end{proof}
\end{proposition}

The next similar proposition for a nonnormal singularity is proved by the
same arguments.
\begin{proposition}\label{nonout}
Let $(X,0)\subset (\CC^n,0)$ be a nonnormal quasihomogeneous hypersurface
singularity with weights $p_1,\ldots,p_n$ and degree $d$ defined by a polynomial
$f$. Let
$\sum p_i\ge d$ and $(X,0)$ is log canonical outside  {\bf 0}
(see definition \ref{onedim}) then
${\bf 1} \in \Gamma_{+}(f)$.
\end{proposition}

\subsection{The quasihomogeneous singularity resolution construction}
\label{toric}
For every hypersurface singularity
$(X,0)\subset (\CC^n,0)$ defined by a polynomial $f$ there exists
a subdivision $\Delta$ of $\sigma$ in $N$ \cite{Var}) depending only
on Newton polyhedron $\Gamma_{+}(f)$ with the following properties.
The variety
$T_N(\Delta)$ is nonsingular, therefore $\Diff_{X({\Delta})}(0)=0$, where
$X(\Delta)$ denotes the proper transform of $X$ on $T_N(\Delta)$.
If
$(X,0)$ is a non-degenerate singularity then the birational morphism
$\psi \colon T_N(\Delta) \to \CC^n$ is the toric log-resolution of
$(\CC^n, X)$ \cite{Var}. The value of non-degenerate singularities is in the
presence of such toric embedded log-resolution. In any case the birational
morphism
$\psi$ is called {\it a partial log-resolution}. Thus
$K_{T_N(\Delta)}+X(\Delta)+E(\Delta)=\psi^*(K_{\CC^n}+X)$.
\par
{\bf Convention.} Let us assume in addition that
$f(x_1,\ldots,x_n)$ is a quasihomogeneous polynomial. Let us choose its weights
${\bf p}=(p_1,\ldots,p_n)$ that $p_i\ne 0$ for
any $i$.

\begin{example}\label{ex2}
If $(X,0)$ is an isolated singularity then $\psi$ is its embedded toric
log-resolution.
\end{example}

(1). Let $(X,0)$ be a normal singularity. It is equivalent to demand that
$\codim_{\CC^n} \Sing X\ge 3$, i.e. $n_i\le n-3$ for all $i$
(see definition \ref{notisol}).
Then
$X(\Delta)$ is normal and $K_{X(\Delta)}=\psi^*K_X + \sum a(E_i,0)E_i$.
As in the example \ref{ex2} variety $X(\Delta)$ can't have the isolated
singularities. Therefore
$\Sing X(\Delta)=\bigcup_{i\in I} \widetilde{\mathcal A}^{n_i}_i$, where
$\widetilde{\mathcal A}^{n_i}_i$ are the irreducible components and
$\psi(\widetilde{\mathcal A}^{n_i}_i)={\mathcal A}^{n_i}_i$
(see definition \ref{notisol}).

\begin{definition} ${\mathcal A}^{n_i}_i$ is called {\it a non-degenerate
component} of singularity if $i\notin I$.
${\mathcal A}^{n_i}_i$ is called {\it a degenerate component} of
singularity if $i\in I$.
\end{definition}
We can recognize these components in
the following way. Let us slightly modify the coefficients of
$f$ in order that the singularity became non-degenerate. If the component
${\mathcal A}^{n_i}_i$ of our singularity remains the component of new
singularity then it is non-degenerate. If it doesn't occur then this
component is degenerate.
\par
In other words the morphism
$\psi$ is the log-resolution only for the non-degenerate
singularity component.
\par
Let us define
$\widetilde V^{n_j}_j,\widetilde W^{n_{k,j}}_{k,j}$ as before
$V^{n_j}_j,W^{n_{k,j}}_{k,j}$ in the definition \ref{notisol}.
Notice that $\psi(\widetilde V^{n_j}_j)=V^{n_j}_j,
\psi(\widetilde W^{n_{k,j}}_{k,j})=W^{n_{k,j}}_{k,j}$ and
$\widetilde V^{n_j}_j,\widetilde W^{n_{k,j}}_{k,j}\nsubseteq \Exc (\psi)$
for all $j,k$.

\begin{lemma}
Along the sets
$\widetilde V^{n_j}_j\setminus((\bigcup_k \widetilde W^{n_{k,j}}_{k,j})\bigcup(
\bigcup_{i\ne j}\widetilde V^{n_i}_i) )$,
$\widetilde W^{n_{k,j}}_{k,j}\setminus \bigcup_{l\ne k}\widetilde
W^{n_{l,j}}_{l,j}$ we have the same singularity types
$\FF_j,\FF_{k,j}$ respectively (cf. definition \ref{notisol}).
\begin{proof} The exceptional divisor corresponding to
${\bf p}$-blow-up is denoted by
$E_{\bf p}$. We have to prove that the type of hypersurface singularities is the
same one in the intersection point of these sets with
$E_{\bf p}$. This means that the singularity given by a polynomial

$$ \hat f=x_i^{-d}f(x_1x_i^{p_1},\ldots,x_i^{p_i},\ldots,x_nx_i^{p_n})=
f(x_1,\ldots,1,\ldots,x_n)
$$
in the neighborhood
of $E_{\bf p}$ coincides with the singularity given by a polynomial
$f$ in the neighborhood of $x_i=1$.
For this purpose it is enough to prove that
$\frac{\partial f}{\partial x_i}({\bf x})=0$, where
${\bf x}\in S=\Sing\{\hat f=0\}$. Since
\begin{gather*}
\frac{\partial (x^d_i\hat f)}{\partial x_i}=\sum_{j\ne i}
p_j\frac{\partial f}{\partial x_j}(x_1x_i^{p_1},\ldots)x_jx_i^{p_j-1}+
p_i\frac{\partial f}{\partial x_i}(x_1x_i^{p_1},\ldots)x_i^{p_i-1}=\\
=\frac{\partial \big(x^d_if(x_1,\ldots,1,\ldots,x_n)\big)}{\partial x_i}=
dx_i^{d-1}f(x_1,\ldots,1,\ldots,x_n),
\end{gather*}
then substituting any ${\bf x}=(x^0_1,\ldots,1,\ldots,x^0_n)\in S$ we get
$p_i\frac{\partial f}{\partial x_i}(x^0_1,\ldots,1,\ldots,x^0_n)=0$. Q.E.D.
\end{proof}
\end{lemma}
By the same argument $\bigcup \Sing E_i\subset \Sing X(\Delta)$.
Therefore there exists a log-resolution
$\varphi \colon\widetilde X \to X(\Delta)$ such that
$\dim\varphi (E)\ge 1$, where $E$ is any component of $\Exc(\varphi)$. Moreover,
the center of every exceptional divisor doesn't lie in  $\bigcup \Sing E_i$.
Such resolution is called {\it a good one}.
\par
(2). Let $(X,0)$ be a nonnormal reduced singularity, i.e.
$\codim_{\CC^n} \Sing X\ge 2$ ($n_i\le n-2$ for all $i$).
Similarly there is the same definition of {\it the non-degenerate and degenerate
singularity components}.
In the same way we construct {\it a good log-resolution}  $\varphi$ of
$\big(T_N(\Delta),X(\Delta)+ E(\Delta)\big)$ with the same properties.
In particular the center of every exceptional divisor doesn't lie in
$ \Sing E(\Delta)$.

\begin{theorem}\label{crqht}
Let $(X,0)\subset (\CC^n,0)$ be a normal quasihomogeneous hypersurface
singularity with weights
$p_1,\ldots,p_n$ and degree $d$ defined by a polynomial
$f$. Then
\begin{enumerate}
\item $\sum p_i \ge d+1$ and $(X,0)$ is canonical outside {\bf 0} if and only
if $(X,0)$ is canonical.
\item $\sum p_i \ge d$ and $(X,0)$ is log canonical outside {\bf 0} if and only
if $(X,0)$ is log canonical.
\end{enumerate}
\begin{proof} The sufficient conditions are obvious. Now we prove their
necessity. Let
$\psi \colon T_N(\Delta) \to \CC^n$ be a partial log-resolution
(see\ \ref{toric}). Let us prove the first assertion.
By proposition \ref{out} we have ${\bf 1} \in \Gamma_{+}(f)^0$. Hence
$a_i(E_i,0)\ge 0$ for all $i$ \cite{Mar}. Since $(X,0)$ is canonical outside
{\bf 0} then $\FF_i,\FF_{k,j}$ are canonical by proposition \ref{equiv}.
This implies that $X(\Delta)$ is canonical. Therefore
$(X,0)$ is a canonical singularity too. The second assertion is proved
similarly. Arguing as above, we see that
$a(E_i,0)\ge -1$, $\FF_i,\FF_{k,j}$ are the lc singularities and
$X(\Delta)$ is lc. By taking a good log-resolution we conclude that
$(X,0)$ is lc since
$\widetilde V^{n_j}_j,\widetilde W^{n_{k,j}}_{k,j}\nsubseteq \Exc (\psi)$
for all $j,k$.
\end{proof}
\end{theorem}

\begin{remark} By  theorem
\ref{crqht} and proposition \ref{equiv} the question about the concrete
quasihomogeneous singularity type is reduced to the same one about
$\FF_i$ and $\FF_{i,k}$ of lesser dimension.
\end{remark}

The next similar criterion for the nonnormal singularity is proved by the
same arguments.
\begin{theorem}\label{noncrqht}
Let $(X,0)\subset (\CC^n,0)$ be a nonnormal quasihomogeneous hypersurface
singularity with weights
$p_1,\ldots,p_n$ and degree $d$. Then
$\sum p_i \ge d$ and $(X,0)$ is log canonical outside {\bf 0} if and only if
$(X,0)$ is log canonical.
\end{theorem}

\begin{example}
Consider the singularities form the example
\ref{first}.
By theorem \ref{crqht} the first singularity is canonical. The second one
is terminal since it is terminal outside {\bf 0} and its hyperplane
section $\{x_3=x_5\}$ has the canonical singularities.
\end{example}

\begin{theorem}\label{crplt}
Let $(X,0)\subset (\CC^n,0)$ be a normal log canonical
quasihomogeneous hypersurface singularity with weights
$p_1,\ldots,p_n$ defined by a polynomial
$f$. Then
\begin{enumerate}
\item If $(X,0)$ is not canonical outside {\bf 0} then
{\bf p}-blow-up is lc, but not plt blow-up, except the case
described in the example \ref{lcex}(1). In particular,
$(X,0)$ is always nonexceptional singularity.
\item $($cf. \cite[3.3]{PrI}$)$ If $(X,0)$ is canonical outside {\bf 0} then
{\bf p}-blow-up is plt one. Assume in addition that
$(X,0)$ is strictly log canonical then
$(X,0)$ is exceptional.
\end{enumerate}
\begin{proof}
Let $\psi \colon \CC^n({\bf p})\to \CC^n$ be a $\bf p$-blow-up.
Its exceptional locus is
$$\Exc(\psi|_{X(\bf p)})=E_{\bf p}=D_{\bf p}\cap X({\bf p})=
\{f(x_1,\ldots,x_n)\subset \PP(p_1,\ldots,p_n)\}.$$
Since $f$ is an irreducible polynomial then $E_{\bf p}$ is an irreducible
variety. It is clear that
$-E_{\bf p}$ is relatively ample and $X(\bf p)$ is a normal variety.
Consider the partial log-resolution of $X$. We have the corresponding morphism
$\varphi \colon T_N(\Delta) \to \CC^n({\bf p})$.
\begin{equation*}
\begin{array}{ccccccc}
&&T_N(\Delta) & \stackrel{\varphi}{\longrightarrow} & \CC^n({\bf p}) &\stackrel{\psi}{\longrightarrow}  & \CC^n \\
& & \bigcup &  & \bigcup &  & \bigcup \\
\widetilde X & \longrightarrow & X(\Delta) & \stackrel{\varphi|_{X(\Delta)}}{\longrightarrow} & X({\bf p})
&\stackrel{\psi|_{X({\bf p})}}{\longrightarrow}  & X
\end{array}
\end{equation*}
One has
\begin{equation}
K_{T_N(\Delta)} + X(\Delta) + D_{\bf p}=
\varphi^*(K_{\CC^n({\bf p})}+X({\bf p})+D_{\bf p})+
\sum_{\RR_{\ge 0}{\bf q}\in \Delta(1), {\bf q}\ne {\bf p}}\alpha_{\bf q}D_{\bf q}.
\end{equation}
By abuse of notation the divisors corresponding to
${\bf p}$ on $T_N(\Delta)$ and
on $\CC^n({\bf p})$ are both denoted by $D_{\bf p}$. Let
${\bf q} \in \sigma_i$ be a primitive vector. By the proof of lemma
\cite[3.2]{PrI} we have
\begin{equation*}
\alpha_{\bf q}=a(D_{\bf q},X({\bf p}))-\frac{q_i}{p_i}=
\langle{\bf q},1\rangle - {\bf q}(f)-\frac{q_i}{p_i}
(\langle{\bf p},1\rangle-{\bf p}(f))-1.
\end{equation*}
Let us prove the first assertion. Obviously
${\bf 1} \in \Gamma_{+}(f)$.
If $q_i=0$ then
$\alpha_{\bf q}=\langle{\bf q},1\rangle - {\bf q}(f)-1\ge -1$.
If $q_i>0$ then we can define ${\bf a}=(a_1,\ldots,a_n)$ from the proportion
$(q_1:\ldots:q_n)=(p_1+a_1:\ldots:p_n+a_n)$. Then $a_i=0$ and
$a_j\ge 0$ for all $j\ne i$. Substitute
${\bf q}=\frac{q_i}{p_i}({\bf p}+{\bf a})$. Then
$\alpha_{\bf q}= \frac{q_i}{p_i}\big(\langle{\bf a},1\rangle
- {\bf a}(f)\big)-1 \ge -1$.

\begin{lemma}
$\Diff_{X({\bf p})}(0)=0$.
\begin{proof}
If $\codim_{\CC^n({\bf p})} \Sing \CC^n({\bf p})>2$ then there is nothing to prove.
Let
$\codim_{\CC^n({\bf p})} \Sing \CC^n({\bf p})=2$. There are only two cases.
\par
(1) $(p_1,\ldots,p_{j-1},p_{j+1},\ldots,p_n)=d>1$ and $p_j\ne 0$.
Hence $\psi(D_{\bf p})=0$.\\
Consider any chart of $\CC^n(\bf p)$

$$
U_k=\CC^n_{y_1,\ldots,y_n}\bigl/\ZZ_{p_k}(-p_1,\ldots,-p_{k-1},1,-p_{k+1},
\ldots,-p_n),\ k\ne j.
$$

The variety $X({\bf p})$ is given by an equation
$f(y_1,\ldots,y_{k-1},1,y_{k+1},\ldots,y_n)=0$ in this chart.
If $\Diff_{X({\bf p})}(0)\ne 0$ then $f({\bf y})=0$ for all
${\bf y}=(y_1,\ldots,y_{k-1},1,y_{k+1},\ldots,y_{j-1},0,y_{j+1},\ldots,y_k)$.
Thus $f=x_j\cdot f'$. It is impossible.
\par
(2) $p_i\ne 0$ if $i\in I$ and $p_i = 0$ if $i\notin I$, where set
$I\subset \{1,\ldots,n\}$ and $2 \le k=|I|\le n-1$.
Hence $\dim\psi(D_{\bf p})=n-k$.\\
In the small neighborhood  $\mathcal U$ of the generic point
$P\in \psi(D_{\bf p})$ the variety
$X$ is isomorphic to $\CC^{n-k}\times X'$. The morphism $\psi$ induces
${\bf p'}$-blow-up $\CC^k({\bf p}')\to \CC^k$ and $X'({\bf p'})\to X'$, where
$p'_i=p_i$ for $i\in I$. By the above statement $\Diff_{X'({\bf p'})}(0)=0$.
Hence $\Diff_{X({\bf p})}(0)=0$.
\end{proof}
\end{lemma}

From this lemma and equality (2) we have
\begin{equation*}
K_{X(\Delta)} + \widetilde E_{\bf p}=
\varphi^*(K_{X({\bf p})}+E_{\bf p})+
\sum\alpha'_{\bf q}D'_{\bf q},\ \text{where}\ \alpha'_{\bf q}\ge -1.
\end{equation*}

Consider a good log-resolution $\widetilde X \to X(\Delta)$.
All exceptional divisors for the pair
$(X(\Delta),\widetilde E_{\bf p}-\sum\alpha'_{\bf q}D'_{\bf q})$
have the discrepancies not less then  $-1$ since
$\widetilde V^{n_j}_j,\widetilde W^{n_{k,j}}_{k,j}
\nsubseteq \Exc (\psi\circ\varphi)$.
By assumption $(X,0)$ is not canonical outside {\bf 0}. Thus
$(X({\bf p})+E_{\bf p})$ is not plt or the case described in
the example \ref{lcex}(1) takes place. By propositions \ref{lcplt} and
\ref{videxc} $(X,0)$ is not exceptional.
\par
Let us prove the second assertion. Arguing as above it is enough to prove that
$\alpha_{\bf q}=\langle{\bf q},1\rangle - {\bf q}(f)-1> -1$
if $q_i=0$ and
$\langle{\bf a},1\rangle- {\bf a}(f)>0$ if $q_i>0$. Actually, if
$\alpha_{\bf q}=-1$ then by considering ${\bf q}$-blow-up or ${\bf a}$-blow-up we
have the contradiction that our singularity is canonical outside
{\bf 0}. It can happen that
${\bf a}$-blow-up is undefined (see \ref{Odatoric}).
Then ${\bf a}(f)=0$ and there is nothing to prove or $f$ is not irreducible, i.e.
we have the contradiction with the normality of $(X,0)$.
\par

Assume in addition that
$(X,0)$ is strictly lc. Since $(X,0)$ is canonical outside {\bf 0} then
$p_i\ne 0$ for all $i$. Thus
$\dim \psi(E_{\bf p})=0$. By theorem \cite[1.9]{Kud2} $(X,0)$ is
exceptional.
\end{proof}
\end{theorem}

\begin{example}\label{lcex}
(1) Let $(X,0)$ be a normal strictly log canonical quasihomogeneous
hypersurface singularity with weights
$p_1,\ldots,p_n$.
Let ${\bf p}$-blow-up $\psi$ is plt one and
$p_i=0$ for some $i$. It means there exists only one index
$k$ such that ${\FF_k}$ is strictly lc.
Note also that ${\FF_k}$ is exceptional, the corresponding component
$\mathcal A^{n_k}_k$ is non-degenerate and
$\psi(E_{\bf p})=\mathcal A^{n_k}_k$.
\par
(2) Let $(X,0)$ be a normal strictly log canonical quasihomogeneous hypersurface
singularity with weights
$p_1,\ldots,p_n$ defined by a polynomial $f$.
Let $(X_{\bf p},0)=(f_{\bf p}=0,0)$ be a singularity of part (1) for some
${\bf p}\in N_{\RR}$.
Then ${\bf p}$-blow-up is also plt one.
\end{example}

\begin{example}
Consider the singularity
$(x^2_1+x^3_2+x^6_3=0,0)\subset (\CC^4,0)$ with weights ${\bf p}=(3,2,1,a)$.
If $a=0$ then the {\bf p}-blow-up is a plt one.
If $a\ne0$ then the {\bf p}-blow-up is a lc but not plt one.
\end{example}

For the nonnormal singularities there is the similar theorem.
\begin{theorem}\label{noncrplt}
Let $(X,0)\subset (\CC^n,0)$ be a nonnormal log canonical quasihomogeneous
hypersurface
singularity with weights ${\bf p}=(p_1,\ldots,p_n)$ defined by a polynomial
$f$. Consider ${\bf p}$-blow-up
$\psi \colon \CC^n({\bf p})\to \CC^n$. Then
$K_{\CC^n({\bf p})}+X({\bf p})+D_{\bf p}$ is lc but not plt.
\end{theorem}

By theorem \ref{crplt} we have the following important corollary.
\begin{corollary}\label{main}
Let $(X,0)\subset (\CC^n,0)$ be a normal hypersurface singularity defined by
a polynomial $f$ and let ${\bf p} \in N_{\RR}$.
Put $(X_{\bf p},0)=(f_{\bf p}=0,0)$. Then
\begin{enumerate}
\item If $(X_{\bf p},0)$ is canonical (resp. log canonical)
then $(X,0)$ is also canonical (resp. log canonical).
\item
Consider {\bf p}-blow-ups $\varphi \colon (X({\bf p}),E_{\bf p})\to (X,0)$ and
$\varphi_{\bf p} \colon (X_{\bf p}({\bf p}),E'_{\bf p})\to (X_{\bf p},0)$.
Then
$(E_{\bf p},\Diff_{E_{\bf p}}(0))=(E'_{\bf p},\Diff_{E'_{\bf p}}(0))$.
\item If $(X_{\bf p},0)$ is log canonical and not  canonical outside {\bf 0} then
{\bf p}-blow-up is lc but not plt one of
$(X,0)$ except the case described in the example \ref{lcex}(2).
In particular $(X,0)$ is always not exceptional.
\item If $(X_{\bf p},0)$ is log canonical and canonical outside {\bf 0} then
{\bf p}-blow-up is plt one of $(X,0)$.
Assume in addition that
$(X,0)$ is strictly log canonical. Then
$(X,0)$ is exceptional.
\end{enumerate}
\begin{proof}
Consider the deformation $F_{t}=t^{-{\bf p}(f)}f(t^{p_1}x_1,\ldots,t^{p_n}x_n)$.
If $t=0$ then $F_0=f_{\bf p}$. For small $t\ne 0$ singularity $F_t$ has the type
$f$. Let $(X_{\bf p},0)$ is lc then $(\CC^n,X_{\bf p})$ is lc by lemma
\cite[7.1.3]{Kollar1}. By Inversion of Adjunction
and \cite[7.8]{Kollar1} the first statement is true.
The statements (2) and (4) are trivial.
To prove the third statement it is enough to show that
$(\CC^n({\bf p}),X({\bf p})+D_{{\bf p}})$ is lc.
Consider the same deformation. If
$t=0$ then $(\CC^n({\bf p}),X_{\bf p}({\bf p})+D_{{\bf p}})$ is lc by the proof of
theorem \ref{crplt}. Therefore corollary
\cite[7.8]{Kollar1} completes the proof.
\end{proof}
\end{corollary}

\begin{corollary}\label{main2}
Let $(X,0)\subset (\CC^n,0)$ be a normal log canonical hypersurface
singularity defined
by a polynomial $f$ and $(f_{\bf p}=0,0)$ be a nonnormal log canonical
singularity for some
${\bf p} \in N_{\RR}$.
Then
{\bf p}-blow-up is lc but not plt one of
$(X,0)$. In particular $(X,0)$ is not exceptional.
\begin{proof} Since $f_{\bf p}$ is a reduced polynomial then it is easy to
check that $\dim \Sing X({\bf p})\le n-2$ in any chart for
${\bf p}$-blow-up \cite[3.7]{PrLect}. Therefore
$X({\bf p})$ is normal. The rest is proved as above corollary
\ref{main}
(applying theorem \ref{noncrplt}).
\end{proof}
\end{corollary}

\begin{remark}
For every non-degenerate log canonical
(resp. canonical) singularity there always exists
${\bf p}\in N_{\RR}$ such that $f_{{\bf p}}$ defines a log canonical
(resp. canonical) singularity \cite[3.5]{PrI}.
\end{remark}

\section{\bf Three-dimensional log canonical hypersurface singularities}

Now we prove the main theorem formulated in the introduction.
\par {\bf Proof of the main theorem.} According to theorem \ref{findpart}
the exceptional strictly log canonical (resp. canonical) hypersurface
singularities are
divided into two types.
\par
(1) There exists a primitive ${\bf p}\in N_{\RR}$ that
$\tilde f_{\bf p}=(f \circ \psi)_{\bf p}$ gives a strictly log canonical
and canonical outside {\bf 0} (resp. canonical)
singularity.
\par
(2) $\tilde f_{\bf p}=t^3+g_2^2(z,x,y)$, where $g_2$ is an irreducible
homogeneous polynomial of degree 2.
\par
The exceptional singularities of second type are classified in the theorems
\ref{classific2} and \ref{classific3}.
\par
Consider the singularities of the first type. The primitive vector
${\bf p}$ is unique. Indeed, let there exist two different
${\bf p_1}$ and ${\bf p_2}$. Then by lemma
\cite[4.1]{PrI} and corollary \ref{main} they give two
different lc blow-ups. Therefore $(X,0)$ is not exceptional by proposition
\ref{lcplt}. Thus ${\bf p}$-blow-up
$\varphi\colon (X({\bf p}),E_{\bf p})\to (X,0)$ is plt one by corollary
\ref{main}.
\par
Let $\tilde f_{\bf p}$ defines a strictly log canonical singularity.
By corollary \ref{main} $(X,0)$ is exceptional. Since
$\dim \varphi(E_{\bf p})=0$ (theorem \cite[1.9]{Kud2}) and
$\langle{\bf p},1\rangle={\bf p}(f)$ then it is easy to prove that
$\Diff_{E_{\bf p}/\PP(\bf p)}(0)=0$ and
$\Diff_{E_{\bf p}}(0)=0$. Hence
$E_{\bf p}$ is a surface with Du Val singularities and
$K_{E_{\bf p}}\sim 0$. By lemma \cite[1.4.1]{Fletcher}
$E_{\bf p}$ is a singular K3 surface.

\par
Let $\tilde f_{\bf p}$ defines a canonical singularity.
Then under condition ---
$\Diff_{E_{\bf p}/\PP(\bf p)}(0)=0$ such singularities are
classified in the theorems
\ref{f5}, \ref{f4} and in the tables of chapter \S 4. The proof of this
classification begins from the subsection \ref{rest}. Q.E.D.

\begin{definition}\label{typem}
A quasihomogeneous hypersurface singularity $f$ belongs to {\it the type}
$\mathcal M_i$ if the maximal compact face of Newton polyhedron
$\Gamma_+(f)$ has the dimension $i$.
\end{definition}

The next theorem  \ref{findpart} divides the three-dimensional exceptional
hypersurface singularities into two types.
\begin{theorem}\label{findpart}
Let $(Z,0) \subset(\CC^4,0)$ be a three-dimensional exceptional
canonical (resp. strictly log canonical) hypersurface singularity defined by
a polynomial $g$. Then there exists
a biholomorphic coordinate change
$\psi \colon (\CC^4,0) \to (\CC^4_{t,z,x,y},0)$ and unique primitive vector
${\bf p}\in N_{\RR}$ such that just one of the following two possibilities holds:
\begin{enumerate}
\item $\tilde g_{\bf p}=(g \circ \psi)_{\bf p}$ defines an exceptional canonical
$($resp. strictly log canonical and canonical outside {\bf 0}$)$ singularity.
\item $\tilde g_{\bf p}=t^3+g_2^2(z,x,y)$, where $g_2$ is an irreducible
homogeneous polynomial of degree 2.
\end{enumerate}
\begin{proof} At first consider the canonical singularity case. If
$g$ is a quasihomogeneous polynomial then there is nothing to prove. Take
${\bf p}\in N_{\RR}$ such that $g_{\bf p}\in \mathcal M_3$,
${\bf 1}\in \Gamma_{+}(g_{\bf p})^0$ and
${\bf p}(g)=1$.
It is clear that 2 or 3-jet of $g_{\bf p}$ is nonzero. Let
$g_{\bf p}$ be not reduced polynomial, i.e. $g_{\bf p}=(g'_{\bf p})^k$, where
$k\ge 2$. Without loss of generality it can be assumed that
$g'_{\bf p}=t+g''$. After the quasihomogeneous coordinate change
$\psi \colon t \longmapsto t-g''$ we get
$\tilde g=g\circ\psi=t^k+\tilde g'$ and $\tilde g_{\bf p}=t^k$. Repeat the
described above process for $\tilde g$. Since ${\bf 1}\in \Gamma_{+}(g)^0$ then
this process will be finished. Hence we may assume that $g_{\bf p}$ is a reduced
polynomial.
\par
Let $(Z_{\bf p},0)=(\{g_{\bf p}=0\},0)$ be a nonnormal singularity.
Then all components $\mathcal A^2_i$ are degenerate (cf. \cite[3.3]{PrI}).
Indeed, let $\mathcal A^2_j$ be a non-degenerate component.
Take the partial resolution
$\psi$. Then there is an exceptional divisor $E$ with discrepancy $a(E,0)\le -1$
such that
$\psi(E)=\mathcal A^2_j$. We obtain the contradiction with
${\bf 1}\in \Gamma_{+}(g_{\bf p})^0$.
By corollary \ref{main2} we have a nonnormal, not log canonical singularity
along some component
$\mathcal A^2_k$. In the neighborhood  of the generic point of
$\mathcal A^2_k$ we have
$Z_{\bf p}\cong \CC^2_{f_1,f_2}\times \{\Phi(f_1(t,z,x,y),f_2(t,z,x,y))=0\}$,
$\deg \Phi\ge 2$.
Since ${\bf 1}\notin \Gamma_{+}(\Phi)$ and 2 or 3-jet of $g_{\bf p}$ is nonzero
it follows that
after the quasihomogeneous biholomorphic coordinate change we get
$f_1=t$. If $\mult_0 f_2=1$ then it is not difficult to show that monomial
$x$, $y$ or $z$ belongs to $f_2$. Thus, after some biholomorphic coordinate
change we can assume that
$\mult_0 f_2\ge 2$. Let
$\Phi_{\bf q}$ be a quasihomogeneous part of $\Phi$ $({\bf q}(\Phi)=1)$, which
defines an isolated singularity. Then
$\wt z+\wt x+\wt y\le 3\cdot\frac12 \wt f_2$. The equality holds if and only if
$f_2=g_2$ is an irreducible homogeneous polynomial of degree 2.
Using this fact in solving the following system
$$
\left\{
\begin{array}{c}
\wt t + \wt z + \wt x+ \wt y >  1 \\
\wt f_1 + \wt f_2  <  1  \\
\wt f_1, \wt f_2 \in \{\frac12,\frac13\}\cup(0,\frac13)  \\
\end{array}
\right.
$$
we get  $\Phi_{\bf q}=f_1^3+f^2_2$. Therefore
$g_{\bf p}=t^3u_1(t,z,x,y)+g_2^2u_2(t,z,x,y)$. Since
${\bf 1}\in \Gamma_{+}(g_{\bf p})^0$ then $u_1$, $u_2$ are the constants and
we have the second case in the theorem condition.
\par
Let $(Z_{\bf p},0)$ be a normal, not canonical singularity. By theorem
\ref{crqht} $(Z_{\bf p},0)$ is not canonical outside {\bf 0}.
By the same argument we have not canonical singularity $\Phi(f_1,f_2,f_3)=0$
along some degenerate component
$\mathcal A^1_k$. If $\mult_0 f_i=1$ for all $i$ then one can take the
corresponding coordinate change and repeat the process.
Let $\mult_0 f_3\ge 2$. Consider a quasihomogeneous part of
$\Phi$. Then
$\wt t + \wt z + \wt x+ \wt y\le \wt f_1 + \wt f_2 +
\frac12\wt f_3+\frac12\wt f_3 \le 1$. Contradiction with the construction of
$g_{\bf p}$. Therefore this theorem is proved for the canonical singularities.
\par
Consider the strictly log canonical singularity case.
Similarly, there is
${\bf p}\in N_{\RR}$ such that $g_{\bf p}\in \mathcal M_3$,
${\bf 1}\in \Gamma_{+}(g_{\bf p})$ and
${\bf p}(g)=1$. Moreover, if $g_{\bf p}$ defines a not log canonical outside
${\bf 0}$ singularity then ${\bf 1}\in \Gamma_{+}(g_{\bf p})^0$.
Therefore, applying the discussed above arguments
we can prove this theorem.
\end{proof}
\end{theorem}

The next theorem classifies all exceptional canonical singularities in the
second case of theorem \ref{findpart}.
\begin{theorem}\label{classific2}
Let $(X,0)=(f=t^3+g_2^2(z,x,y)+f_{\ge 5}'(z,x,y)=0,0)\subset (\CC^4,0)$ be a
canonical singularity. Then there exists a plt blow-up
$\varphi\colon (Y,S)\to (X\ni 0)$ such that
$(S,\Diff_S(0))=(\overline \FFF_{10}, \frac12H_1+\frac23H_2)$, where
$\overline \FFF_{10}=\PP(10,1,1)$ is a cone, $H_1$ and
$H_2$ are its hyperplane sections. Moreover,  $(X,0)$ is exceptional
(resp. weakly exceptional) and 6-complementary if and only if
any irreducible factor of the nonzero 5-jet
$f'_5$ has the multiplicities at most 3 (resp. at most 4).
\begin{proof} Consider ${\bf p}=(4,3,3,3)$-blow-up
$\phi\colon (X({\bf p}),E)\to X$. Then
$(E,\Diff_E(0))=(t+g_2^2\subset \PP(4,1,1,1),\frac23\{t=0\})=(\PP^2,\frac43 C)$,
where $C$ is an irreducible conic (see construction \ref{ED}). Also
$a(E,0)=0$. If 5-jet of $f'$ is zero then polynomial $f$
has a view
${t'}^3+g_2^2(1,x',y')+{z'}^6f''$ in the chart
$\CC^4_{t',z',x',y'}/\ZZ_3(2,1,0,0)$.
This new singularity is not canonical, therefore 5-jet of $f'$ is nonzero.
Now our theorem is proved by exhaustion of all cases for
5-jets.
\par
Let $f=t^3+g^2_2+z^5+f'_{\ge 6}(z,x,y)$. We have the following picture\\

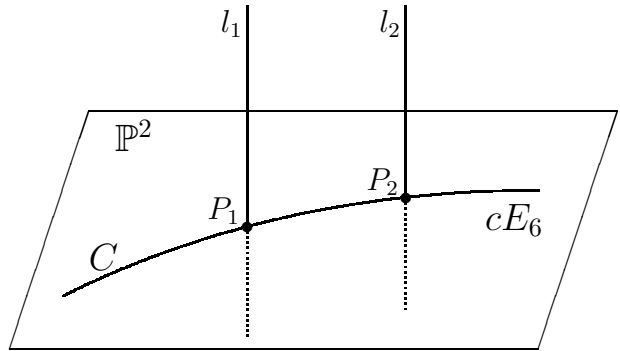
\begin{equation*}
\begin{picture}(260,150)(0,0)
\put(0,20){\line(1,3){30}}
\put(0,20){\line(1,0){200}}
\put(200,20){\line(1,3){30}}
\put(30,110){\line(1,0){200}}
\put(40,95){\large{$\PP^2$}}
\linethickness{0.7pt}
\qbezier[300](20,40)(100,80)(200,80)
\thinlines
\put(30,50){\large{$C$}}
\put(180,65){\large{$cE_6$}}
\linethickness{0.7pt}
\put(90,150){\line(0,-1){84}}
\put(150,150){\line(0,-1){73}}
\put(80,140){{$l_1$}}
\put(140,140){{$l_2$}}
\thinlines
\put(90,66){\circle*{4}}
\put(150,77){\circle*{4}}
\put(75,69){{$P_1$}}
\put(135,80){{$P_2$}}
\linethickness{0.7pt}
\dottedline{2}(90,66)(90,25)
\dottedline{2}(150,77)(150,35)
\end{picture}
\end{equation*}

The variety $X({\bf p})$ in the neighborhood of $P_i$
(the points $P_i$ have the coordinates $(0,1,y_i)$ in $\PP^2$, where
$g_2(0,1,y_i)=0$)
is
$$
\big({t'}^3+\omega^2+{z'}^5{x'}^3+{x'}^6f''=0\subset
(\CC^4_{t',\omega,z',x'},0)\big)/\ZZ_3(2,0,0,1).
$$
Depending on the type of $f'$ this singularity can be nonisolated along the curves
$l_1$ and $l_2$ not lying on $E$. By considering the invariants of group action
we get $X({\bf p})$  to be isomorphic to
$$
(x_1^2x_4^5+x_1x_3^2+x_2^3+f'''=0,0)\subset(\CC^4_{x_1,x_2,x_3,x_4},0)
$$
in the neighborhoods of $P_i$.
\par
The surface $E$ is given by the equations $\{x_1=x_3=0\}$ in these neighborhoods.
Along the curve $C$ the variety
$X({\bf p})$ has $cE_6$ singularity. Let $H$ be a general hyperplane section of
$X({\bf p})$ passing through any point from $C\backslash \{P_1,P_2\}$. Let
$C'=E\cap H$. Consider the minimal resolution of
$H$. We have the following graph
\\
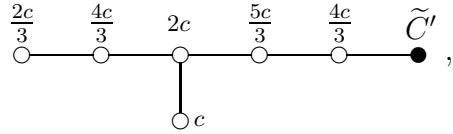
\begin{equation*}
\begin{picture}(200,55)(0,0)
\put(5,30){\circle{6}}
\put(0,39){$\frac{2c}3$}
\put(8,30){\line(1,0){24}}
\put(35,30){\circle{6}}
\put(30,39){$\frac{4c}3$}
\put(38,30){\line(1,0){24}}
\put(65,30){\circle{6}}
\put(60,39){\scriptsize{$2c$}}
\put(68,30){\line(1,0){24}}
\put(65,27){\line(0,-1){19}}
\put(95,30){\circle{6}}
\put(90,39){$\frac{5c}3$}
\put(98,30){\line(1,0){24}}
\put(125,30){\circle{6}}
\put(120,39){$\frac{4c}3$}
\put(128,30){\line(1,0){24}}
\put(155,30){\circle*{6}}
\put(165,28){,}
\put(150,36){{$\widetilde C'$}}
\put(65,5){\circle{6}}
\put(70,3){\scriptsize{$c$}}
\end{picture}
\end{equation*}
where $\widetilde C'$ is a proper transform of $C'$.
The numbers near graph vertexes (exceptional curves) denote their
corresponding discrepancies for the pair $(H,cC')$.
This implies that
log canonical threshold of $(H,cC')$ is equal to 1/2 and the extraction
of only central
graph vertex is an inductive blow-up of
$(H,\frac12 C')$
(see theorem \cite[1.5]{Kud2}).
Let $c=c(X({\bf p}),E)$ be a log canonical threshold. The pair $(X({\bf p}),E)$
has the same singularities in the neighborhoods of the points
$P_1$ and $P_2$. Therefore by connectedness theorem
\cite[17.4]{Koetal} $LCS(X({\bf p}),cE)=C$. Hence $c=1/2$.
\par
Let us take a blow-up of $\CC^4({\bf p})$ along the curve
$C$ with the weights (6,4,3,0).
Then this blow-up induces
$\psi\colon \widetilde X({\bf p})\to X({\bf p})$. Thus
$K_{\widetilde X({\bf p})}+\widetilde S+\frac12E_{\widetilde X({\bf p})}=
\psi^*(K_{X({\bf p})}+\frac12E)$.
It is obvious that
$f'''$ doesn't influence on the exceptional surface
$\Exc\psi=\widetilde S\cong \FFF_{10}$.
In what follows it will be demonstrated that the proper transform
$E_{\widetilde X({\bf p})}$ of $E$ can be contracted to the point.
Therefore $E_{\widetilde X({\bf p})}\cap \widetilde S$ is a minimal section
of $\FFF_{10}$.
One have
$\Diff_{\widetilde S}(\frac12E_{\widetilde X({\bf p})})=
\Diff_{\widetilde S}(0)+\frac12E_{\widetilde X({\bf p})}|_{\widetilde S}=
\frac23l_{\infty}+\frac12l_0+\frac23l'_0+\frac12\cdot\frac13l_{\infty}=
\frac56l_{\infty}+\frac12l_0+\frac23l'_0$, where $l_{\infty}$ is a minimal section,
$l_0$ and $l'_0$ are the different zero sections. The sections
$l_0$ and $l'_0$ have the intersections in the two points with multiplicities 5.
These points are lying over points
$P_1$ and $P_2$. It is clear that log surface
$(\widetilde S,\Diff_{\widetilde S}(\frac12E_{\widetilde X({\bf p})}))$ is klt,
i.e. $\psi$ is an inductive blow-up of
$(X({\bf p}),\frac12E)$. Let
$L_{X({\bf p})}\in|-n(K_{X({\bf p})}+\frac12E)|$ be a general element
from the very ample linear system over $X$ for $(n\gg 0)$. Then
$K_{X({\bf p})}+\frac12E+\frac1nL_{X({\bf p})}$ is lc and numerical trivial
over $X$. One have
$K_{\widetilde X({\bf p})}+\widetilde S+\frac12E_{\widetilde X({\bf p})}+
\frac1nL_{\widetilde X({\bf p})}=
\psi^*\phi^*(K_X+\frac1nL_X)$, where $L_X=\phi(L_{X({\bf p})})$.
Apply $K_{\widetilde X({\bf p})}+\widetilde S+(\frac12+
\varepsilon)E_{\widetilde X({\bf p})}+
\frac1nL_{\widetilde X({\bf p})}$ -- MMP ($0<\varepsilon\ll 1$). Since
$\rho(E_{\widetilde X({\bf p})})=1$ then at the first and unique step of MMP
we have a divisorial contraction
$E_{\widetilde X({\bf p})}$ to the point. Thus we get a plt blow-up
$\varphi\colon (Y,S)\to (X\ni 0)$ and
$(S,\Diff_S(0))=(\overline \FFF_{10}, \frac12H_1+\frac23H_2)$, where
$H_1$ are $H_2$ are the hyperplane sections of cone. Divisor
$K_S+\Diff_S(0)+\frac16l$ is anti-ample and not klt, where
$l$ is a cone generator passing through the intersection point of
$H_1$ and $H_2$ (of multiplicity 5). By theorem \ref{CrWexc}
$(X,0)$ is not weakly exceptional.
\par
In the same way one can consider the remaining cases for
5-jets of $f'$. In these cases the surface
$S$ doesn't change, $H_1$ and $H_2$ are the hyperplane sections too.
Only the multiplicities of $H_1$ and $H_2$ intersection change.
Now our theorem follows from the next trivial lemma.
\end{proof}
\end{theorem}

\begin{lemma} Log Del Pezzo surface
$(\overline \FFF_{10}, \frac12H_1+\frac23H_2)$ is exceptional
(resp. weakly exceptional) and 6-complementary if and only if
its hyperplane sections $H_1$ and $H_2$ have the multiplicities at most 3
(resp. at most 4) in any intersection point.
\end{lemma}

The next theorem classifies all exceptional strictly log canonical
singularities in the second case of theorem \ref{findpart}.
\begin{theorem}\label{classific3}
Let $(X,0)=(f=t^3+g_2^2(z,x,y)+f_{\ge 5}'(z,x,y)=0,0)\subset (\CC^4,0)$ be a
strictly log canonical singularity. Then $f'_5=0$ and $f'_6\ne 0$. Moreover
$(X,0)$ is exceptional and 1-complementary if and only if
$f'_6$ doesn't have an irreducible factor of multiplicity 6.
In this case there exists a plt blow-up
$\varphi\colon (Y,S)\to (X\ni 0)$ such that
$(S,\Diff_S(0))=(S,0)$, where
$S$ is a surface obtained by the section contraction of the
elliptic surface
$S'\to \PP^1$ to Du Val singularity of type $\AAA_1$
$($kodaira dimension $\kappa(S')=0)$.
The degenerate elliptic fibers of $S'$ are always elliptic cuspidal curves.
The singular points of $S$ are Du Val singularities.
They lie only in the cusp vertexes. Depending on the irreducible factors
of $f'_6$ we have the following cases:
$12$(it is a number of degenerate fibers)
$\AAA_0$(it means the singularity type in the cusp vertex;
$\AAA_0$ is a smooth point by the definition);
$2\AAA_2+8\AAA_0$;
$4\AAA_2+4\AAA_0$; $6\AAA_2$; $2\DDD_4+6\AAA_0$;
$2\DDD_4+2\AAA_2+2\AAA_0$; $4\DDD_4$; $2\EEE_6+4\AAA_0$;
$2\EEE_6+2\AAA_2$; $2\EEE_8+2\AAA_0$.
\begin{proof} The proof is similar to the theorem
\ref{classific2}. Consider
 ${\bf p}=(4,3,3,3)$ - blow-up
$\phi\colon (X({\bf p}),E)\to X$. Then
$(E,\Diff_E(0))=(t+g_2^2\subset \PP(4,1,1,1),\frac23\{t=0\})=(\PP^2,\frac43 C)$,
where $C$ is an irreducible conic. Also
$a(E,0)=0$. By theorem \ref{classific2} $f'_5=0$. If $f'_6=0$ then
a polynomial $f$ in the chart
$\CC^4_{t',z',x',y'}/\ZZ_3(2,1,0,0)$ is ${t'}^3+g_2^2(1,x',y')+{z'}^9f''$.
This singularity is not log canonical therefore 6-jet of
$f'$ is nonzero.
Now our theorem is proved by exhaustion of all cases for
6-jets.
Let us prove it in two cases:
$f'_6=z^6+f'_{\ge 7}(z,x,y)$,
$f'_6=z^5x+f'_{\ge 7}(z,x,y)$. Other cases are considered likewise.
\par
Along $C\backslash\bigcup P_i$ the variety
$X({\bf p})$ has nonisolated strictly log canonical singularity
$(x_1^3+x_1x_3^2+x_2^3=0,0)\subset(\CC^4_{x_1,x_2,x_3,x_4},0)$.
\par
Consider the first case. Then
$i=1,2$ and the variety
$X({\bf p})$ in the neighborhood of $P_i$
(the points $P_i$ have coordinates $(0,1,y_i)$ in $\PP^2$, where $g_2(0,1,y_i)=0$)
is
$$
(x_1^3x_4^6+x_1x_3^2+x_2^3+f'''=0,0)\subset(\CC^4_{x_1,x_2,x_3,x_4},0).
$$
That is, $P_i$ are the centers of log canonical singularities of
$X({\bf p})$. Therefore $(X,0)$ is not exceptional by the definition.
\par
Consider the second case. Then
$i=1,2,3,4$.
The variety $X({\bf p})$ in the neighborhoods of $P_1,P_2$
is
$$
(x_1^3x_4^5+x_1x_3^2+x_2^3+f'''=0,0)\subset(\CC^4_{x_1,x_2,x_3,x_4},0),
$$
but in the neighborhoods of $P_3,P_4$ is
$$
(x_1^3x_4+x_1x_3^2+x_2^3+f'''=0,0)\subset(\CC^4_{x_1,x_2,x_3,x_4},0).
$$
Let us take a blow-up of $\CC^4({\bf p})$ along the curve
$C$ with weights (1,1,1,0). Then this blow-up induces
$\psi\colon \widetilde X({\bf p})\to X({\bf p})$. Thus,
$K_{\widetilde X({\bf p})}+ S'=\psi^*K_{X({\bf p})}$.
It is obvious that
$f'''$ doesn't influence on the exceptional elliptic surface
$S'\to \PP^1$. It is easy to check that the degenerate cuspidal fibers lie
over
$P_i$.
The surface $S'$ has the singularities $\EEE_8$ in the degenerate fibers
for $i=1,2$.
Divisor $K_{X({\bf p})}+S'$ is plt.
Apply $K_{\widetilde X({\bf p})}+S'+\varepsilon E_{\widetilde X({\bf p})}$
-- MMP ($0<\varepsilon\ll 1$).
Since
$\rho(E_{\widetilde X({\bf p})})=1$
then on the first and unique step of MMP
we have a divisorial contraction
$E_{\widetilde X({\bf p})}$ to the point.
So we get a required  plt blow-up $\varphi\colon (Y,S)\to (X\ni 0)$.
\end{proof}
\end{theorem}

\subsection{ Classification of three-dimensional exceptional quasihomogeneous
canonical singularities} \label{rest}
According to the proof of the main theorem it remains to classify the
quasihomogeneous exceptional singularities. The classification process
consists of two parts:
\begin{enumerate}
\item At first we have to classify the quasihomogeneous polynomials. For any
such polynomial one can correspond its weight
${\bf p}$. By the main theorem the vector ${\bf p}$ is unique.
It will be shown that there is the finite number of weights and consequently
the finite number of types of exceptional singularities.
\item The second part is the investigation of given singularities on
the exceptionality. Also the minimal complement index will be found.
\end{enumerate}

\begin{theorem}\cite{Reid1}
Let $(X,0) \subset(\CC^4,0)$ be a canonical hypersurface singularity defined by
a polynomial $f$ and
$(X,H)$ is not log canonical for any
hyperplane section $H$. Then $f$ belongs to one of the following types:
\begin{enumerate}
\item[$(\Upsilon_1)$]$f=t^2+z^3+zf_1(x,y)+f_2(x,y)$, where
$\deg f_1 \ge 5, \deg f_2 \ge 7.$
\item[$(\Upsilon_2)$]$f=t^2+f'(z,x,y)$, where $\deg f'\ge4$ and 4-jet is nonzero.
\item[$(\Upsilon_3)$]$f=t^2+f'(z,x,y)$, where $\deg f'\ge5$ and 5-jet is nonzero.
\item[$(\Upsilon_4)$] The rank of quadratic part of $f$ is equal to 0 and
3-jet is nonzero.
\end{enumerate}
\end{theorem}

\par {\bf The further notations.} Further we always work under the
following notations:
$(X,0)=(f(t,z,x,y)=0,0)\subset (\CC^4_{t,z,x,y},0)$ denotes a canonical
quasihomogeneous hypersurface singularity.
Since one have to classify exceptional singularities we can assume that
$(X,H)$ is not log canonical, where $H$ is any hyperplane section
(see example \ref{lc}).
\par
The classification of the quasihomogeneous parts will be made by
"Newton line rotation" \cite[\S 11]{AVG}.
For this purpose we first find the finite number of the parts of type
$\mathcal M_2$ (see definition \ref{typem}). After that we obtain the required
list of quasihomogeneous singularities with the help of the simple
computer program.
Therefore the main problem is to find the parts of type
$\mathcal M_2$.
At first consider the singularities of type $\Upsilon_1$.
\par
{\bf Singularities of type $\Upsilon_1$.}
The following theorem describes such singularities.
\begin{theorem}\label{Newton1}
All exceptional canonical quasihomogeneous singularities of type
$\Upsilon_1$ can be obtained by the rotations of plane in
$\CC^3_{z,x,y}$ passing through the monomial
$z^3$ and one of the following monomials:
$x^7$,$x^8$,$x^7y$,$x^9$,$x^8y$,$x^7y^2$,$x^{10}$,$x^9y$,$x^8y^2$,$x^7y^3$,
$x^{11}$,$x^{10}y$,$x^9y^2$,$x^8y^3$,$x^7y^4$ and $zx^5$,$zx^5y$,$zx^7$,
$zx^6y$,$zx^5y^2$. The rotation condition is that the weight of
$x$ is greater then the weight of $y$
$(\wt x\ge \wt y)$ and the point
${\bf 2}=(2,2,2)$ lies over the rotation plane
(see theorem \ref{crqht}).
\begin{proof} Note that the monomial
$zx^6$ is missing in the list since the monomials
$z^3$, $zx^6$, $x^9$ are lying on the straight line.
One can choose
the quasihomogeneous weights of
$f$ such that the quasihomogeneous degree is equal to 1.
\par
In the theorem statement the monomials $x^m$ of type $zx^ay^b$, $x^ay^b$
satisfy the following conditions:
$a\ge b$; $\wt x\ge \wt y$;
$\wt x+\wt y>\frac16$ (see theorem \ref{crqht}); singularity defined by
$t^2+z^3+x^m$ is not lc.
\par
Let us prove that these monomials give all singularities.
Assume the converse. Let the quasihomogeneous polynomial
$f$ gives an exceptional singularity
$(X,0)$ of type $\Upsilon_1$ and every its monomial doesn't belong to the list
(in this case $\wt x> \wt y$).
Let $zx^{a_1}y^{b_1}$ be a monomial of minimal usual degree
from
$f_1(x,y)$ and $x^ay^b$ be a monomial of minimal usual degree
from $f_2(x,y)$.
Since $\wt x> \wt y$ it follows that the degree of
$x$ in these two monomials is maximal then the one in the remaining monomials
of $f$. Also the degree of
$y$ is minimal.
\par
If $f_1=0$ or $f_2=0$ then $f'=t^2+z^3+x^ay^b$
(or $+zx^{a_1}y^{b_1}$) defines an exceptional log canonical singularity.
In fact there are two different vectors from $N_{\RR}$. Hence, $f'$ gives
an exceptional singularity (see main theorem proof).
By taking
${\bf p}=\big(\frac12,\frac13,\frac1a$ (or $\frac2{3a_1}$)$,1\big)$ we have
$f_{\bf p}=f'$.
Hence $(X,0)$ is not exceptional.
\par
Therefore there are the following 4 cases.
We prove that they are impossible.
\par
{\it Case 1.} $f'=t^2+z^3+x^ay^b$ and $f''=t^2+z^3+zx^{a_1}y^{b_1}$ give the
not log canonical singularities. Then they are not log canonical outside
{\bf 0} along
$\{t=z=y=0\}$. Therefore $(X,0)$ is not log canonical outside {\bf 0}
along $\{t=z=y=0\}$. Contradiction.
\par
{\it Case 2.} $f'=t^2+z^3+x^ay^b$ gives
a not log canonical singularity (i.e. $a<b\ge 7$) and
$f''=t^2+z^3+zx^{a_1}y^{b_1}$ gives
a log canonical singularity (i.e. $a_1\le 4$ and $b_1\le 4$).
Thus if
$b_1=0$ then by the definition of $\Upsilon_1$ we have $a_1\ge 5$. Contradiction.
We have the following linear equation system on the weights:

$$
\left\{
\begin{array}{ccll}
a\wt x + b\wt y & = & 1 &\\
a_1\wt x + b_1\wt y & = & 2/3,& d=ab_1-ba_1\ne 0 \\
\end{array}
\right.
$$
\par
Let $a_1\ge b_1$. To prove theorem in this case it is enough to find
new weights $\wt x, \wt y$ such that
$f_{\bf p}=f''$, where ${\bf p}=(1/2,1/3,\wt x, \wt y)$.
Decrease $\wt x$ by
$\varepsilon$ and increase $\wt y$ by $\sigma=a_1\varepsilon/b_1$
(after it the second system equation remains valid).
These new weights are required. Actually we have to check that
$\langle{\bf p}, x^ay^b\rangle >1$, i.e.
$-d=-ab_1+ba_1>0$. It is true since
$a/b<1$ and $a_1/b_1\ge 1$.
\par
Let $a_1< b_1$. If $-d>0$ then $(X,0)$ is also nonexceptional as above.
Let $d>0$ then $\wt x=(b_1-\frac23b)/d>0$. Hence $b_1>\frac23b\ge 14/3$.
Contradiction.
\par
{\it Case 3.}
$f'=t^2+z^3+x^ay^b$ gives a log canonical singularity
(i.e. $a\le 6$ and $b\le 6$) and
$f''=t^2+z^3+zx^{a_1}y^{b_1}$ gives a
not log canonical singularity (i.e. $a_1<b_1\ge 5$).
For the same reason $b\ne 0$.
If $d>0$ then as before in the case 2) we can find ${\bf p}$ such that
$f_{\bf p}=f'$. If $d<0$ then as before in the case 2) we have
$b_1-\frac23b<0$, i.e. $b>\frac32b_1\ge 15/2$. Contradiction.
\par
{\it Case 4.} $f'=t^2+z^3+x^ay^b$ and $f''=t^2+z^3+zx^{a_1}y^{b_1}$ give the
log canonical singularities. Similarly let us decrease $\wt x$ and increase
$\wt y$ in order that the first or second system equation remains valid.
In one of the cases we get
$f_{\bf p}=f'$ or $f_{\bf p}=f''$.
\end{proof}
\end{theorem}

\par
{\bf The investigation of singularities on the exceptionality.}
At first let us recall the main facts about the hypersurfaces in the
weighted projective spaces
\cite{Fletcher}. The greatest common divisor of the numbers
$a_1,\ldots,a_n$ is denoted by
$(a_1,\ldots,a_n)$.

\begin{proposition} Let $(a_1,\ldots,a_n)=1$ and
$q_i=(a_1,\ldots,\hat{a}_i,\ldots,a_n)$, where $a_i \in \NN$ for all $i$. Then
$\PP=\PP(a_1,\ldots,a_n)\cong \PP({a_1}/q_i,\ldots,a_i,\ldots,a_n/q_i)$.
If $q_i=1$ for all $i$ then $K_{\PP}\sim\OO_{\PP}(-\sum a_i)$ and
$\OO_{\PP}(1)^{n-1}=1/(a_1\cdots a_n)$. In this case the weights of
$\PP$ are called well-defined.
\end{proposition}

\begin{proposition}\label{Fletch} In the same notations let
$(X_d \subset \PP)=(g(x_1,\ldots,x_n)\subset \PP(a_1,\ldots,a_n))$ be a
hypersurface of degree $d$. Then
$(X_d \subset \PP)= (g(x_1,\ldots,x_i^{1/q_i},\ldots,x_n)\subset
\PP({a_1}/q_i,\ldots,a_i,\ldots,a_n/q_i))$. If $\Diff_{X_d/\PP}(0)=0$ and the
weights of
$\PP$ are well-defined then
$K_{X_d}=(K_{\PP}+X_d)|_{X_d}=\OO_{X_d}(d-\sum a_i)$.
\end{proposition}

\begin{definition}\label{w-d}
Assume the weights of the weighted projective space $\PP$ are well-defined. Then
the hypersurface
$X_d\subset \PP$ is called {\it well-formed} if one of the following equivalent
conditions holds:
\begin{enumerate}
\item $\codim_{X_d}(X_d\cap \Sing\PP)\ge 2$,
\item $\Diff_{X_d/\PP}(0)=0$,
\item for all $i\ne j$ we have
$(a_1,\ldots,\hat{a}_i,\ldots,\hat{a}_j,\ldots,a_n)|d$.
\end{enumerate}
\end{definition}

\subsection{\bf The structure of purely log terminal blow-up of
well-formed singularity}\label{ED}
Let ${\bf p}=(a_1,\ldots,a_n)$ be the weights of the quasihomogeneous log canonical
singularity $g(x_1,\ldots,x_n)=0$ and $a_i\in \NN$ for all $i$.
By the main theorem these weights are defined uniquely for the exceptional
singularity (in particular, $a_i>0$ for all $i$). A
${\bf p}$-blow-up of $\CC^n$ induces a plt blow-up
$(Y,E)\to (\{g=0\},0)$. Given log Fano variety is
$\big(E,\Diff_E(0)\big)=\big(E,D\big)=
\big(\tilde g(x_1,\ldots,x_n)\subset \PP(\tilde a_1,\ldots,\tilde a_n),D\big)=
\big( g(x_1^{1/q_1},\ldots,x_n^{1/q_n})\subset
\PP(a_1q_1/(q_1\cdots q_n),\ldots,a_nq_n/(q_1\cdots q_n)),\sum
\frac{q_i-1}{q_i}\{x_i=0\}\big)$.
When calculating
$(E,D)$ it is important that the hypersurface
$E\subset \PP$ be well-formed. If it is not well-formed
(the weights $\tilde a_1,\ldots,\tilde a_n$ are well-defined) then we have
other different $D$ (see example \ref{ex3}).
\par
By abuse of notation the coordinates of $g=0$ and the coordinates of the
weighted projective space
are both denoted equally.
\par
The quasihomogeneous degrees of
$f$, $\tilde f$ are denoted by
$d$ and $\tilde d$ respectively. In our three-dimensional case given
log Del Pezzo surface is always denoted by
$(E,D)$.
\par
The following easy proposition allows to consider only finite number of
cases when studying given log Del Pezzo surface on the exceptionality.

\begin{proposition}\label{DelPezzo1}\cite[2.3]{Sh2}
Let $(S,\Delta)$ be a weakly log Del Pezzo surface,
i.e. $-(K_S+\Delta)$ is nef and big.
Assume that all
1,2,3,4 or 6-complements $\Delta^+$ are exceptional, i.e.
$(S,\Delta^+)$ is klt. Then log surface is exceptional.
\end{proposition}

Using definition \ref{defcompl} we get the next corollary:
\begin{corollary}\label{DelPezzo2}
Notation is as in \ref{DelPezzo1}. Let
$K_S+\frac{\up{(n+1)\Delta}}n$ be ample for all $n=1,2,3,4,6$.
Then $(S,\Delta)$ is exceptional.
\end{corollary}

\begin{corollary}\label{DelPezzo3}
Let $(X,0)$ be a three-dimensional quasihomogeneous log canonical singularity.
Assume its any hyperplane section is not log canonical. Let
$(E,D=\sum d_iD_i)$ be a log Del-Pezzo surface obtained by a plt blow-up.
If $d_k\ge \frac67$ for some $k$ and $\Diff_{E/\PP}(0)= 0$
(it is the condition of well-formedness) then
$(X,0)$ is exceptional.
\begin{proof} Let $(X,0)$ is not exceptional. Then by
proposition
\ref{DelPezzo1} there exists 1,2,3,4 or 6-complement $D^+$.
By the definition of complement
$D^+\ge D_k+\sum_{i\ne k}d_iD_i$, where $D_k=\{x_k=0\}$. Since
$\Diff_{E/\PP}(0)= 0$ then
$\frac1{q_1\cdots q_4}(d-\sum_{i\ne k}a_i)=\tilde d-
\sum_{i\ne k}\frac{\widetilde a_i}{q_i}\le 0$.
Therefore by theorem \ref{crqht} the general hyperplane
section $x_k+h(x_1,\ldots,x_4)=0$ of $(X,0)$ is lc. This contradiction concludes
the proof.
\end{proof}
\end{corollary}

\begin{remark}\label{DelPezzo3-1}
The corollary
\ref{DelPezzo3} is hypothetically true in any dimension $n$.
In the statement we need to assume
$d_k\ge \frac{|RN_{n-1}|}{|RN_{n-1}|+1}$, where $|RN_n|=\max\{m|m\in RN_n\}$.
The set $RN_n$ is defined in the following way:
$m\in RN_n$ if there exists the nonexceptional log Fano variety of dimension
$n$ such that the minimal index of nonexceptional complements is equal to
$m$. A hypothesis is that
$|RN_n|<\infty$ for all $n$.
For $n\le 3$ it was proved in \cite{Sh2}, besides $RN_1=\{1,2\}$ and
$RN_2=\{1,2,3,4,6\}$ \cite{Sh2}.
\end{remark}

\par
{\bf The common principles of investigation on the exceptionality.}
\par(1)
At first one can check that $\Diff_{E/\PP}(0)=0$
(it is convenient to verify the third condition in the definition \ref{w-d}).

\par(2)
If $D$ contains a component with coefficient $\ge \frac67$ then
by corollary \ref{DelPezzo3} the singularity is not exceptional.
\par
Any complement $D^+$ of
$(E,D)$ can be extended to the complement $K_{\PP}+E+\widehat D^+
\sim_{\QQ} 0$, where $\widehat D^+|E=D^+$.
Hence, the structure of Weyl divisor group of weighted projective space $\PP$
implies that
$\Supp D^+$ is cut by polynomials.
\par
For any number
$n\in \NN$ satisfying the following property there exists an $n$-complement.
\par
{\it Property:} Let
$r_n=\sum_{i=1}^4\tilde a_i-\tilde d-\frac1n\sum_{i=1}^4\up
{\frac{q_i-1}{q_i}(n+1)}\tilde a_i$,
i.e. $K_E+\frac{\up{(n+1)D}}n\sim_{\QQ}-r_n\OO_E(1)$.
Then
$r_n\ge 0$ and $r_n=\frac1n(\sum\tilde a_ib_i)$, where $b_i\in \NN$ for all $i$.
\par(3)
By corollary \ref{DelPezzo1} there are two cases: $(E,D)$ is exceptional
and there is no 1,2,3,4 or 6-complement or
there exists 1,2,3,4 or 6-complement. Later on we always assume the second
case to take place.
\par
First we determine the singular points of
$\tilde f\subset \PP$ and the singular points of curves from $\Supp D$.
It is very convenient to apply lemma
\ref{Laz} for these purposes. The answer depends on the coefficients of $f$.
\par (4)
Let $D^+$ be a complement of $(E,D)$. By theorem \ref{CrExc} the question
about the exceptionality is equivalent to the next one: Is pair
$(E,D^+)$ klt?
This question is local. Consider this pair in the chart
$
U_k=\CC^3\bigl/\ZZ_{a_k}(a_1,\ldots,a_{k-1},a_{k+1},\ldots,a_4)
$
of $\PP$.
The klt property of pair remains true by the finite dominant morphisms
\cite[3.16]{Kollar1}. Hence our investigation reduces to the finite number
of questions:
Is pair $(\CC^2_{u,v},B=\sum c_i\{h_i(u,v)=0\})$
klt?  Is pair $(\CC^3_{u,v,w},B=\{h=0\}+\sum c_i\{h_i(u,v)=0\})$ plt,
where $h=0$ defines Du Val singularity
($h$ is a local equation of $E$ in the corresponding chart)? In solving these
problems it is convenient to consider the analog of inductive blow-up
\cite[1.5]{Kud2}
\big(it is a blow-up with the unique exceptional divisor
$C$ such that
$(C,\Diff_C(\widetilde B))$ is klt or plt, where $\widetilde B$ is a
proper transform of $B$ \big)
and apply lemma \ref{inductive}.

\par
{\bf The singularity investigation of type $\Upsilon_1$ on the
exceptionality.}

\par
Further we consider the most interesting cases and the common examples
illustrating
the investigation of given log Del Pezzo surfaces on the exceptionality.

\begin{example}
${\bf f=t^2+z^3+x^{11}+y^{13}}.$ Then $\Big(E,D\Big)=
\Big(\PP^2,\frac12l_1+\frac23l_2+\frac{10}{11}l_3+\frac{12}{13}l_4 \Big)$, where
$l_i$ are the lines in the general position.
By corollary \ref{DelPezzo2} or by corollary \ref{DelPezzo3} given log surface
is exceptional. A minimal complement index is equal to 66. The corresponding
complement is
$D^+=\frac12l_1+\frac23l_2+\frac{10}{11}l_3+\frac{61}{66}l_4$.
A conjecture is that 66 is maximal number for log surfaces with the
standard coefficients and consequently for every three-dimensional local
contraction (singularity) \cite{Sh2}. In other words $|RN_3|=66$
(see remark \ref{DelPezzo3-1}).
\end{example}

\begin{example}
${\bf f=t^2+z^3+g(z,x,y)=t^2+z^3+zyf'_2(x^2,y^3)+}$\\
${\bf+xf'_3(x^2,y^3)}.$
The lower indexes denote the binary form degrees. These forms are different
although they are both denoted by
${\bf f'}$.
Let us remark that we got rid of the monomial $z^2xy^2$ before writing $f$.
Namely we made the coordinate change
$z\longmapsto z+\alpha xy^2$ for some $\alpha$.
The rotation takes place through the unique part
$t^2+z^3+x^7$ of type $\mathcal M_2$. Therefore the coefficient of
$x^7$ is not equal 0. Put $f'_2=ax^{2i}f_{2-i}(x^2,y^3)$,
$f'_3=x^{2j}f_{3-j}(x^2,y^3)$.

\begin{proposition1} $ij=0$ $($if $a=0$ then $i\ne 0$ by definition $)$
if and only if
$(X,0)$ is exceptional and 2-complementary.
\begin{proof} A plt blow-up with weights $(21,14,6,4)$ gives
$\Big(E,D\Big)=\Big(t+z^3+g(z,x,y)\subset \PP(21,7,3,2),\frac12\{t=0\}\Big)=
\Big(\PP(7,3,2),\frac12\{t=0\}\Big)$. Let $i\ge 1$ and $j\ge 1$. Then
$K_E+\frac12\{t=0\}+\frac12\{x=0\}\equiv 0$ is not klt hence singularity
is not exceptional. Check it. Let $i=j=1$. Then this statement is equivalent
to the following one:
$(\CC^2_{u,v},\frac12\{(u^3+a'u^2v+v^3)\cdot v=0\})$ is not klt.
Indeed, take a blow-up with weights
$(1,1)$. Then the exceptional divisor with discrepancy
--1 appears. The cases with other
$i,j$ are considered similarly.
\par
Now, let $i=0$. Consider two more difficult cases:

\begin{enumerate}
\item $g(z,x,y)=zy^7+x^7$
\item $g(z,x,y)=zy(x^2+y^3)^2+x(x^2+y^3)^3$
\end{enumerate}

The first case is exceptional obviously. In the second case it is enough to
check that
$K_E+\frac12\{t=0\}+\frac14\{x^2+y^3=0\}\equiv 0$ is klt. This statement
is equivalent to  $(\CC^2_{u,v},B)$ be
klt, where $B=\frac12\{u^3+a'u^2v+v^3=0\}+
\frac14\{ v=0\}$. To prove this or analogous assertion the next easy lemma
can be applied.

\begin{lemma}\label{inductive}
Let $(Z\ni P,D)$ be a log pair and let
$f\colon Y\to Z$ be a birational morphism.
Assume that the exceptional set of $f$ consists of one irreducible
exceptional divisor $C$.
If $a(C,D)> -1$ (resp. $a(C,D)=-1$) and
$K_Y+C+D_Y$ is plt ($D_Y$ is a proper transform of $D$), then
$(Z\ni P,D)$ is klt (resp. lc) and exceptional.
\end{lemma}

Take the blow-up with the weights
(1,1).  Then the single exceptional divisor
$C\cong\PP^1$ with discrepancy --3/4 appears. The divisor
$K_C+\Diff_C(\widetilde B)=
K_C+\frac12P_1+\frac12P_2+\frac12P_3+\frac14P_4$ is klt, where
$P_i$ are the distinct points and $\widetilde B$ is a proper transform of $B$.
The case $j=0$ is considered similarly. Divisor
$D^+=\frac12\{t=0\}+\frac12\{x=0\}$ is 2-complement.
\end{proof}
\end{proposition1}

Consider quotient singularity $\CC^3/H_{168}$, where $H_{168}$ is a simple
Klein's group of order 168. The singularity pointed out has the same type
as in the previous example. It is exceptional (cf. \cite{PrMar}).
Note also that it is a degenerate singularity.
\end{example}

\begin{example}\label{ex12}
${\bf f=t^2+z^3+\underline{a}zx^if'_{5-i}(x,y^2)+
\underline{b}yx^jf'_{7-j}(x,y^2)}.$
The plane can rotate through two different parts of type
$\mathcal M_2$ for this singularity. The first part is
$t^2+z^3+x^7y$. The second one is $t^2+z^3+zx^5$.
Therefore the coefficient of $zx^5$ or the coefficient of $x^7y$ is not equal
0. On occasion of notation of
${\bf f}$ see the conventions in \S 4.
A plt blow-up with weights $(15,10,4,2)$ gives
$\Big(E,D\Big)=\Big(t+z^3+\ldots \subset \PP(15,5,2,1),\frac12\{t=0\}\Big)=
\Big(\PP(5,2,1),\frac12\{t=0\}\Big)$.
Let $i\ge 3$ and $j\ge 5$. Then
$K_E+\frac12\{t=0\}+\frac14\{x=0\}\equiv 0$ is not klt. Hence
singularity is not exceptional.
Indeed, consider $(E,D)$ in the neighborhood of the point $(0:0:0:1)$.
Then we have to show that
$(\CC^2,\frac12\{u^3+\underline{a}uv^i+\underline{b}v^j=0\}+
\frac14\{v=0\})$ is not klt. Take a blow-up with weights
$(3,2)$. The exceptional divisor with discrepancy
--1 appears. Q.E.D.
\par
If $i\le 2$ or $j\le 4$ then singularity is not exceptional and 2-complementary.
The requirement  $i\le 2 \ ||\ j\le 4$ also means that it is true for the
multiplicities of other common components of type
$x+cy^2$. We write this requirement only for
$x$, because after making the quasihomogeneous coordinate change
$x\longmapsto x'-cy^2$ we have the initial case.
Consider two more difficult cases.
Let $j=4$ and $\underline{a}=0$. It is enough to check that
$(\CC^2,B=\frac12\{u^3+v^4=0\}+
\frac14\{v=0\})$ is klt. Take a blow-up with weights
(4,3). The single exceptional
divisor $C$ with discrepancy --3/4 appears. The pair $(C,\Diff_C(\widetilde B))=
(\PP^1,\frac12P_1+\frac78P_2+\frac23P_3)$ is klt, where
$P_i$ are the distinct points. Lemma \ref{inductive} completes the proof.
The remained difficult case $i=2$ and $\underline{b}=0$ is studied similarly.
\end{example}

\par
{\bf Singularities of type $\Upsilon_2$.}
Let us consider the exceptional singularities of this type.
By corollary \ref{main2} there is a biholomorphic map such that
4-jets can be only the following: $z^4$ (type $\Upsilon^{(1)}_2$),
$z^3x$ (type $\Upsilon^{(2)}_2$), $z^3y$ (type $\Upsilon^{(3)}_2$).

\par
{\it Singularities of type
$\Upsilon^{(1)}_2 (t^2+z^4+z^2f_1(x,y)+zf_2(x,y)+f_3(x,y))$.}
To rotate Newton line (the rotation condition is $\wt x\ge \wt y$) it is
enough to fix the monomials of the following forms: $x^ay^b$ and
$zx^ay^b$.
Using the same arguments just as for singularities of type
$\Upsilon_1$ (theorem \ref{Newton1}) we get the following cases to rotate
Newton line: $x^5$, $x^6$, $x^5y$, $x^7$,
$x^6y$, $x^5y^2$ and $zx^4$, $zx^5$, $zx^4y$.

\begin{example}
${\bf f=t^2+z^4+zx^4y+x^3y^n, n=4,5.}$
The weights of plt blow-up are
$(2(4n-3),(4n-3),3n-4,7)$. Besides
$(E,D)=(f\subset\PP(2(4n-3),(4n-3),3n-4,7) ,0) $ and $K_E=(n-6)\OO_E(1)$.
If $n=5$ then there is no 1,2,3,4 and 6-complement.
Therefore the singularity is exceptional and
$\frac17\{y=0\}$ gives 7-complement of minimal index. Let
$n=4$. There is 4-complement $D^+=\frac14\{x=0\}$.
Let us prove that $(E,D^+)$ is not plt. For this purpose it is enough to prove
that $(\CC^3_{u,v,w},\{u^2+v^4+w^3=0\}+\frac14\{w=0\})$
is not plt. Take the weighted blow-up with weights
(6,3,4). Then exceptional divisor with discrepancy
--1 appears. Thus the singularity is not exceptional.
\end{example}

\par
{\it Singularities of type $\Upsilon^{(2)}_2 (t^2+z^3x+g(z,x,y))$.}
Using the same arguments just as for singularities of type
$\Upsilon_1$ (see theorem \ref{Newton1}) we get the following cases to rotate
Newton line ($\wt x \ge \wt y$):
$x^5$,$x^6$,$x^5y$,$x^7$,$x^6y$,$x^5y^2$,$x^8$,$x^7y$,$x^6y^2$,$x^5y^3$,$x^9$,
$x^8y$,$x^7y^2$,$x^6y^3$,$x^5y^4$
and $zx^4$,$zx^4y$,$zx^6$,$zx^5y$,$zx^4y^2$. Note that
$zx^5$ is a subcase in $x^7$ case.

\begin{example}\label{ex1}
${\bf f=t^2+z^3x+az^2yx^{i_1}(x+y^3)^{i_2}f_{2-i_1-i_2}(x,y^3)
+}$
${\bf +bzy^2x^{j_1}(x+
y^3)^{j_2}f_{3-j_1-j_2}(x,y^3)+
x^{l_1}(x+y^3)^{l_2}f_{5-l_1-l_2}(x,y^3)}$.
The lower indexes denote the binary form degrees.
After the coordinate change
$x\longmapsto x+\alpha y^3$ for some $\alpha$ we got rid of the monomial
$z^3y^3$. The rotation takes place through the unique part $t^2+z^3x+x^5$ of
type
$\mathcal M_2$. Therefore the coefficient of
$x^5$ is not equal 0. Since variable $x$ is present in the monomial
$z^3x$ then we must select other common component $x+y^3$. A plt blow-up is
the blow-up with weights $(15,8,6,2)$.
Thus $(E,D)=(\PP(4,3,1),\frac12\{t=0\})$. There are the following
1,2,3,4 and 6-complements: 2-complement -- $D+\frac12\{y=0\}$ and
6-complement -- $D+\frac16\{x+\beta y^3=0\}$.
Hence, the minimal complement index is equal to 2.
Since the curve $\{y=0\}$ doesn't pass through the singular points
of $\Supp D$ then the first complement is not considered.
Let us consider the remaining cases $\beta=0,\beta=1$.
\par
Let $\beta=0$. If $i_1\ge 2$ ($i_1>2 $ denotes that $a=0$),
$j_1\ge 3$ ($j_1>3 $ denotes that $b=0$), $l_1=3$ (if $l_1\ge 4 $ then our
singularity is not canonical) then  $(\CC^2_{u,v},B)$ is not klt,
where $B=\frac12\{u^3v+v^3+au^2v^2+buv^3=0\}+\frac16\{v=0\}$.
Actually take a blow-up with weights
$(2,3)$. Then the exceptional divisor with discrepancy $-1$ appears.
In other cases its complement is not exceptional.
For example, let us prove it if
$i_1\ge 2, j_1=2, l_1\ge 3$. In this case
$B=\frac12\{u^3v+au^2v^2+uv^2+v^{l_3}=0\}+\frac16\{v=0\}$.
Take a blow-up with weights
$(1,2)$. The discrepancy of the single exceptional divisor $C$ is equal to
$-\frac56$. Besides
$\Diff_C(\widetilde B)=\frac34 P_1+\frac23 P_2+\frac12 P_3$, where
$\widetilde B$ is a proper transform of $B$. Hence
$K_C+\Diff_C(\widetilde B)$ is klt, i.e. by lemma
\ref{inductive} pair $(\CC^2,B)$ is klt.
\par
Let $\beta=1$. Then $B=\frac12 \{h=0\}+\frac16 \{v=0\}$, where
$h(u,v)=u^3+au^2v^{i_2}+buv^{j_2}+v^{l_2}$. It is easy to verify that
$(\CC^2,B)$ is klt if $b\ne 0$. Let $b=0$. Similarly one can verify that
$(\CC^2,B)$ is always klt if $i_2\le 1$.
Let $i_2\ge 2$. If
$l_2\le 4$ then this pair is klt. If $l_2=5$ then it is not klt.
\par
{\it Summary.} The canonical singularity is exceptional and 2-complementary
except the case
$i_1 \ge 2 \ \& \ j_1\ge 3\ \& \  l_1=3$ and the case
$i_2 \ge 2 \ \& \ b=0\ \&\  l_1=3$. These conditions are written taking into
account other common components
(see example \ref{ex12}).
\end{example}

\par
{\it Singularities of type $\Upsilon^{(3)}_2.$} The following theorem takes
place.
\begin{theorem}\label{Newton2-3}
All exceptional quasihomogeneous singularities of this type can be obtained
by the rotations of plane in $\CC^3_{z,x,y}$ passing through the
monomial $z^3y$ and one of the following monomials:
$x^7$,$x^8$,$x^7y$,$x^9$,$x^8y$,$x^7y^2$
and $zx^5$, $zx^6$, $zx^5y$,$z^2x^3$.
The rotation condition is that the weight of $x$ is
greater then the weight of $y$
$(\wt x > \wt y)$ and the point
${\bf 2}=(2,2,2)$ lies over the rotation plane.
\begin{proof} Let $f=t^2+z^3y+z^2f_1(x,y)+zf_2(x,y)+z^2f_3(x,y)$. Denote the
monomial of usual minimal degree from $z^2f_1(x,y)$ by
$z^2x^ay^b$. If $b\ne 0$ then by taking the quasihomogeneous coordinate
change $z\longmapsto z - \alpha f_1(x,y)y^{-1}$ one can assume that $f_1(x,y)=0$.
Applying the same arguments as in the theorem
\ref{Newton1} we complete the proof.
Let $b=0$. Then using $\wt x+\wt y+\wt z > 1/2$ we get
$a\le 3$. Corollary \ref{main2} implies that $a=3$. Q.E.D.
\end{proof}
\end{theorem}

\par
{\bf The singularities of type $\Upsilon_3$.}
For the exceptional singularities of this type there are three cases for
5-jets.
\par THE FIRST CASE. 5-jet is a point (type $\mathcal M_0$),
i.e. $f_5=z^ax^by^c$, where
$a\ge b,a\ge c$ and $a+b+c=5$. By corollary \ref{main2} we have $a\ge 3$. Thus
5-jet is one of the following:
$z^5$ (type $\Upsilon^{(1)}_3$),
$z^4x$ (type $\Upsilon^{(2)}_3$), $z^4y$ (type $\Upsilon^{(3)}_3$),
$z^3x^2$ (type $\Upsilon^{(4)}_3$), $z^3y^2$ (type $\Upsilon^{(5)}_3$),
$z^3xy$ (type $\Upsilon^{(6)}_3$).
As in the theorems \ref{Newton1}, \ref{Newton2-3} there are the following
cases to rotate Newton line
$(\wt x\ge\wt y)$:
\par
{\it Singularities of type $\Upsilon^{(1)}_3$.}
The rotation monomials are $x^6$,
$x^5y$, $x^4y^2$
and $zx^5$, $zx^4y$, $zx^3y^2$.

\par
{\it Singularities of type $\Upsilon^{(2)}_3$.}
The rotation monomials are $x^6$,
$x^5y$, $x^4y^2$
and $zx^5$, $zx^4y$, $zx^3y^2$.

\par
{\it Singularities of type $\Upsilon^{(3)}_3$.}
The rotation monomials $(\wt x>\wt y)$ are  $x^6$,
$x^5y$ and $zx^5$, $zx^4y$. All singularities in the case $zx^4y$ are not
well-formed.

\par
{\it Singularities of type $\Upsilon^{(4)}_3$.}
The rotation monomials are $x^6$, $x^5y$,
$x^4y^2$, $x^3y^3$, $x^7$,
$x^6y$, $x^5y^2$, $x^4y^3$, $x^3y^4$
and $zx^5$, $zx^4y$, $zx^3y^2$.

\par
{\it Singularities of type $\Upsilon^{(5)}_3$.}
The rotation monomials $(\wt x>\wt y)$ are $x^7$ and $zx^5$.

\par
{\it Singularities of type $\Upsilon^{(6)}_3$.}
The rotation monomials are $x^5$, $x^6$,
$x^5y$, $x^7$,
$x^6y$, $x^5y^2$ and $zx^5$, $zx^4y$.
All singularities in the cases $x^5y, x^6y, zx^4y$ are not well-formed.

\par THE SECOND CASE. 5-jet is a line, i.e. $f_5\in\mathcal M_1$
(type $\Upsilon^{(7)}_3$). To rotate
Newton line ($\wt x\ge \wt y$) we have to fix the following parts:
$g_5(z,x)$ is a binary form of degree 5,
$z^5+zx^3y$, $z^5+x^4y$, $z^4x+zx^3y$, $z^4x+x^4y$, $z^4y+z^2x^3$, $z^4y+zx^4$,
$z^3x^2+zx^3y$, $z^3x^2+x^3y^2$, $z^3y^2+z^2x^3$.
This list was obtained in the following way.
The first monomial
$x^{m_1}$ contains a variable of maximal usual degree among
other monomials. This variable is
denoted by
$z$. The second monomial $x^{m_2}=z^ax^by^c$ satisfies the following
conditions:
$b\ge c$ and $t^2+x^{m_1}+x^{m_2}$ gives not log canonical outside
${\bf 0}$ singularity.
\par
To prove that these cases describe all possibilities one may use
the theorem \ref{Newton1} proof.
\par
Consider one-dimensional 5-jet
$f_5=z^5+\alpha z^3xy+\beta zx^2y^2$. Let us prove that this case doesn't
realize.
If $\beta\ne 0$ then $\exists {\bf p}\in N_{\RR}$ such that
$f_{\bf p}=t^2+f_5$ is a log canonical outside {\bf 0} singularity. Therefore by
corollary \ref{main2} our singularity is not exceptional. Actually
let $x^n=z^cx^ay^b$ ($c+a+b\ge 6$) be another monomial and $d=b-a$ be a
determinant of the following linear equation system on the weights.

$$
\left\{
\begin{array}{ccl}
\wt x + \wt y & = & 2/5\\
a\wt x + b\wt y & = & 1-(1/5)c \\
\end{array}
\right.
$$

Let
$d<0$. Since $\wt x\ge \wt y$ then $c+a+b\le 5$. Contradiction.
Let $d>0$. Decrease $\wt x$ on $\varepsilon$ and increase $\wt y$
on $\varepsilon$. Then for new vector
${\bf p}$ we have
$\langle{\bf p}, x^n\rangle >1$. Q.E.D.
\par
Let $\beta=0$ and $x^n=z^cx^ay^b$ ($c+a+b\ge 6$) be another monomial.
As above we have $d>0$. Using conditions $\wt x>0$ and $\wt y>0$
one have $b\ge 4$ for
$c=0,1$ and one have $b\ge 3$ for $c=2,3$. That is,
$t^2+z^3xy+x^n$ gives a not log canonical outside {\bf 0} singularity
along $\{t=z=x=0\}$. Hence the original singularity is not log canonical.
\par
Other cases are considered similarly except the case
$f_5=z^4y+z^3x^2$.
In this case the reader will have no difficulty in proving by
exhaustion that only two new singularities
$t^2+z^4y+z^3x^2+xy^5$ and
$t^2+z^4y+z^3x^2+azy^5+bx^2y^4$ ($b\ne 0$) appear.
It is easy to check that they are not exceptional.

\par
THIRD CASE. 5-jet is a plane, i.e. $f_5\in\mathcal M_2$
(type $\Upsilon^{(8)}_3$).
The classification of such singularities reduces to the description of all
curves of degree 5 in $\PP^2$.

\begin{theorem}\label{f5}
Let $f=t^2+f_5(z,x,y)$ gives a canonical singularity
$(X,0)$, where $f_5$ is a homogeneous polynomial of degree 5.
Then ${\bf p}=(5,2,2,2)$-blow-up induces a plt blow-up
of $(X,0)$. Besides $(E,D)=\big(t+f_5(z,x,y)\subset
\PP(5,1,1,1),\frac12\{t=0\}\big)=(\PP^2,\frac12C)$, where $C=\{f_5=0\}$.
Then $(X,0)$ is exceptional if and only if one of the following
possibilities holds.
\begin{enumerate}
\item The curve $C \subset \PP^2$ is irreducible. Its singular points can be
only ordinary cusps and ordinary double points.
\item The curve $C \subset \PP^2$ is the union of irreducible cubic and
irreducible conic which intersect in 6 distinct points.
The singular point of cubic (if it exists)
can be only ordinary cusp or ordinary double point.
\end{enumerate}
\begin{proof} Since $(X,0)$ is normal singularity then each irreducible
component of the curve
$C=\{f_5=0\}$ has the multiplicity 1.
Hence we have three cases:
$C$ is an irreducible curve,
$C$ is an union of irreducible cubic and irreducible conic,
$C$ is an union of line and quartic.
The corresponding log Del Pezzo surface is
$(E,D)=(\PP^2,\frac12C)$.
\par
Consider the first case. Let the singularity of curve be worse then ordinary cusp.
The reader will easily prove that
$K_{\PP^2}+\frac12C+\frac12l\sim_{\QQ} 0$ is not klt, where
the straight line $l$ passes through the singularity (in the cases of
singularities
$(u^2+v^n=0,0)\subset (\CC^2,0)$, where $n=4,5$ we take $l=\{u=0\}$).
\par
Consider the second case. If the cubic and conic intersect with multiplicity
2 in the point $P$ then the same divisor is not klt, where
the straight line $l$ is a tangent of conic in the point $P$.
\par
Consider the third case. If
$l$ is a straight line of
$\Supp C$ then the same divisor is not klt.
\end{proof}
\end{theorem}

\begin{example}${\bf f=t^2+z^3xy+x^i(x+y^2)^jf_{5-i-j}(x,y^2)=t^2+}$\\
${\bf +z^3xy+h(x,y)}$.
The rotation takes place through the unique part
$t^2+z^3xy+x^5$.
Therefore the coefficient of $x^5$ is nonzero.
A plt blow-up is the blow-up with weights
$(15,7,6,3)$. Then $(E,D)=(t^2+zxy+h(x,y)\subset \PP(5,7,2,1),\frac23\{z=0\})$.
There are the following 1,2,3,4 and 6-complements: 3-complement
$D+\frac13\{y=0\}$ and
6-complement $D+\frac16\{x+\beta y^2=0\}$.
If $i\ge 1$ then 6-complement for $\beta=0$ is 3-complement
$D+\frac13C$, where $C=\{t=x=0\}$.
The minimal complement index is equal to 3 .
As in example \ref{ex1} one can consider only the second complement.
\par
Let $\beta=0$. If $i=3$ (if $i_1\ge 4 $ then the singularity is
not canonical) then the singularity is
not exceptional. Indeed, we have to show that
$(\CC^3_{u,v,w},B)$ is not plt, where
$B=\{u^2+vw+w^3=0\}+\frac23\{v=0\}+\frac16\{w=0\}$.
Take a blow-up with weights
$(3,4,2)$.
Then the exceptional divisor with discrepancy --1 appears.
If $i\le 2$ then $(E,D+\frac16\{x=0\})$ is klt.
For example, prove it for $i=2$. We claim that pair
$(Z,B)$ is klt, where
$Z=\{u^2+vw+w^2=0\} \subset (\CC^3,0)$ and
$B=\frac23\{v=0\}|_Z+\frac16\{w=0\}|_Z$.
Take a blow-up with weights $(1,1,1)$.
Then the single exceptional divisor $C$ with the discrepancy
$a(C,B)=-\frac56$ appears. Also
$\Diff_C(\widetilde B)=\frac23 P_1+\frac23 P_2+\frac13 P_3$.
Hence $K_C+\Diff_C(\widetilde B)$ is klt. Q.E.D.
\par
Let $\beta=1$. For the same reason the singularity is exceptional if
and only if $j\le 4$.
\par
{\it Answer.} Canonical singularity is exceptional and 3-complementary
if and only if
$i \le 2 $ and $j\le 4$.
\end{example}

{\bf Singularities of type $\Upsilon_4$.}
These singularities are considered likewise.
To rotate Newton line ($\wt x\ge \wt y$) there are four cases.
For underlined parts all singularities are not well-formed or they
belong to other cases.

\par (A) $f_3\in \mathcal M_0$. Then $f_3$ is equal to one of the following
monomials: $t^3$, $t^2z$, $t^2x$, $t^2y$.
Let us consider every monomial separately.

\begin{enumerate}
\item After coordinate change we get $f=t^3+tg(z,x,y)+h(z,x,y)$. Since
${\bf 1}\in \Gamma_{+}(f)^0$ then $f_4\ne 0$. We have the following
cases ---
$z^2x^2+x^5$, $z^2x^2+x^4y$, $z^2x^2+x^3y^2$, $\underline{z^2x^2+x^2y^3}$;
$z^2xy+x^5$, $\underline{z^2xy+x^4y}$; $f_4=\varphi_4(z,x)$ is a
binary form of degree 4,
$z^4+zx^2y$, $z^4+x^3y$, $z^4+tx^3$, $z^4+tx^2y$; $z^3x+zx^2y$, $z^3x+x^3y$,
$z^3x+tx^3$, $z^3x+tx^2y$; $z^3y+z^2x^2$, $z^3y+tx^3$; $tz^3+z^2x^2$,
$\underline{tz^3+zx^3}$, $tz^3+zx^2y$, $\underline{tz^3+x^3y}$;
$\underline{z^2x^2+azx^2y+bx^2y^2}$, $z^2x^2+tz^2x$,
$z^2x^2+tx^2y$; $\underline{tz^2x+zx^2y}$; $tz^2x+x^2y^2$;
$f_4\in \mathcal M_2$. The exceptional singularities of last case are
classified in the theorem \ref{f4}.

\item After coordinate change we get  $f=t^2z+tg(x,y)+h(z,x,y)$. Then also
$f_4\ne 0$. We have the following cases --- $f_4=\varphi_4(z,x)$ is a binary
form of degree 4, $z^4+x^3y$, $z^4+x^5$,
$\underline{z^4+x^4y}$, $z^4+x^3y^2$, $z^4+tx^3$, $z^4+tx^2y$; $z^3x+x^3y$,
$z^3x+x^5$, $z^3x+x^4y$, $z^3x+x^3y^2$, $z^3x+tx^3$, $z^3x+tx^2y$,
$z^3y+x^5$; $\underline{z^3y+x^4y}$, $\underline{z^3y+zx^3}$, $z^3y+z^2x^2$,
$z^3y+tx^3$; $z^2x^2+x^5$, $z^2x^2+x^4y$, $z^2x^2+x^3y^2$, $z^2x^2+tx^2y$;
$z^2y^2+x^5$; $z^2xy+x^5$, $\underline{z^2xy+x^4y}$.

\item  After coordinate change we get $f=t^2x+tg(z,y)+h(z,x,y)$. Then also
$f_4\ne 0$. We have the following cases --- $z^4+x^3y$, $z^4+x^5$,
$z^4+x^4y$, $z^4+x^3y^2$, $z^4+zx^2y$, $z^4+zx^4$, $z^4+zx^3y$, $z^4+zx^2y^2$;
$\underline{z^3x+x^2y^2}$, $\underline{z^3x+x^5}$, $\underline{z^3x+x^4y}$,
$\underline{z^3x+x^3y^2}$, $\underline{z^3x+x^2y^3}$, $\underline{z^3x+x^6}$,
$\underline{z^3x+x^5y}$, $\underline{z^3x+x^4y^2}$, $\underline{z^3x+x^3y^3}$,
$\underline{z^3x+x^2y^4}$, $\underline{z^3x+zx^2y}$, $\underline{z^3x+zx^4}$,
$\underline{z^3x+zx^3y}$, $\underline{z^3x+zx^2y^2}$;
$z^3y+x^5$, $z^3y+x^4y$, $z^3y+x^6$, $\underline{z^3y+x^5y}$,
$z^3y+x^4y^2$, $z^3y+zx^4$, $\underline{z^3y+zx^3y}$,
$z^3y+z^2x^2$.

\item After coordinate change we get  $f=t^2y+tg(z,x)+h(z,x,y)$. Then also
$f_4\ne 0$.  We have the following cases --- $z^4+x^5$, $z^4+zx^4$;
$z^3x+x^5$, $z^3x+x^6$, $\underline{z^3x+x^5y}$,
$z^3x+zx^4$, $z^3x+tx^3$.
\end{enumerate}

\par (B) $f_3\in \mathcal M_1$. Then $f_3$ is one of the following polynomials:
\begin{enumerate}
\item $f_3=g_3(t,z)$. The rotation monomials are $x^4$, $x^5$, $x^4y$, $zx^3$.

\item $f_3=t^3+z^2x$. The rotation monomials are
$x^4$, $x^5$, $x^4y$, $x^6$, $x^5y$,
$x^4y^2$, $x^7$, $x^6y$, $x^5y^2$, $x^4y^3$, $x^8$, $x^7y$, $x^6y^2$, $x^5y^3$,
$x^4y^4$, $tx^3$, $tx^3y$, $tx^5$, $tx^4y$, $tx^3y^2$. The monomials
$zx^3$, $zx^4$,
$zx^3y$, $tx^4$ are the subcases of monomials $x^5$, $x^7$, $x^5y^2$, $x^6$
respectively.

\item $f_3=t^3+z^2y$. The rotation monomials are
$x^7$, $x^8$, $\underline{x^7y}$,
$zx^4$, $tx^5$.

\item $f_3=t^2z+z^2x$. The rotation monomials are $x^4$, $x^5$, $x^4y$, $x^6$,
$x^5y$, $x^4y^2$, $tx^3$, $tx^4$, $tx^3y$. The monomial $zx^3$ is a subcase of
monomial $x^5$.

\item $f_3=t^2z+z^2y$. The rotation monomials are $x^5$, $x^6$,
$\underline{x^5y}$, $\underline{zx^3}$, $tx^4$.
\end{enumerate}

\par (C) $f_3\in \mathcal M_2$. Then after coordinate change $f_3$ belongs to
one of the followings types: case
(A), case (B) or a log canonical threshold $c(f_3)=1$.
In the last case by corollaries \ref{main}, \ref{main2} $(X,0)$ is not
exceptional.

\par (D) $f_3\in \mathcal M_3$.
Then after coordinate change $f_3$ belongs to
one of the followings types: case
(A), case (B), case(C) or the general hyperplane section
has log canonical singularities, i.e. $(X,0)$ is not exceptional.

\begin{example}
${\bf f}$ has a type ${\bf f_3(t,z)+y^2f_2(t,z)+(y^4+x^2y)f_1(t,z)+}$
${\bf +f_2(x^2,y^3)}$.
In this case we must assume that the binary form
$ f_3(t,z)$ doesn't have the multiple irreducible factors.
This singularity is obtained by Newton line rotation for the case
$ f_3(t,z)+x^4\in \mathcal M_2$.
After some quasihomogeneous coordinate change (see conventions in \S 4) we get
$f=t^3+ktz^2+lz^3+t(ax^2y+by^4)+z(cx^2y+dy^4)+etzy^2+x^{2i}f_{2-i}(x^2,y^3)$.
Besides $|k|+|l|\ne 0$.
A plt blow-up is the blow-up with weights
$(4,4,3,2)$. Then $(E,D)=(f(t,z,x^{1/2},y)\subset \PP(2,2,3,1),\frac12\{x=0\})$.
There are the following 1,2,3,4,6 complements:
$D+\frac12\{y=0\}$; $\frac23\{x=0\}$;
$D+\frac14\{\alpha t+\beta z+ \gamma y^2=0\}$;
$D+\frac16\{y=0\}+\frac16\{\alpha t+\beta z=0\}$.
Since the coefficient of $x^4$ is nonzero then
it is enough to consider the second and the third complement ($\gamma=0$).
\par
Let $e\ne 0$. Assume that
$l=0 (k\ne 0)\ \&\ z(cx^2y+dy^4)=0\ \&\ i=2$ then
$D^+=\frac14\{t=0\}+\frac12\{x=0\}$ is non-klt 4-complement.
\par
Let $e=0$.
If $|b|+|d|= 0\ \&\ i\ge 1$ then
$D^+=\frac23\{x=0\}$ is non-klt 4-complement.
For
$e\ne 0$ the previous condition on the exceptionality remains only necessary one.
Moreover
$D^+=\frac14\{t+\alpha z=0\}+\frac12\{x=0\}$ is non-klt
4-complement if
$i=2\ \&\ \text{G.C.D}\big(t^3+ktz^2+lz^3,
t(ax^2y+by^4)+z(cx^2y+dy^4)\big)=t+\alpha z$.
\end{example}

\begin{example}\label{ex3}
${\bf f=t^2x+z^3y+x^5+zxy^3}$.
A plt blow-up is the blow-up with weights
$(16,11,8,7)$. Then $(E,D)=(f(t,z,x,y)\subset \PP(16,11,8,7),0)$.
Let us consider the neighborhood of the point
$(0:0:0:1)$ in some chart. Then
$(\CC^3,S+\frac14T)=(\CC^3_{u,v,w},\{u^2w+v^3+w^5+vw=0\}+\frac14\{w=0\})$
is not plt. Indeed, take a blow-up
$\psi$ with the weights $(1,2,4)$. Hence
the singularity is not exceptional.
The singularities of such kind have the complex structure.
For example, in our case $(S,0)$ is Du Val point of type
$\AAA_5$. So to prove the exceptionality we have to consider
the inductive blow-up.
\par
In the capacity of the example let us show that
$(C,B)=(C,\Diff_C(\frac14\widetilde T'))$ is klt in our case,
where $C$ is an exceptional curve of blow-up
$\psi|_{\widetilde S}\colon \widetilde S\to S$ and
$\widetilde T'$ is a proper transform of
$T'=T|_S$. The curve $(C\subset \PP)=(uw+v^3+vw\subset \PP(1,1,2))$ is not
well-formed. Therefore we can't use the results \ref{ED} to calculate
$(C,B)$. The singularities of
$\widetilde S$ are
$\big(u^2+v^3+v\subset \CC^3_{u,v,w}\big)/\ZZ_4(-1,2,1)=\AAA_3$ and
$\big(u^2w+1+w\subset \CC^3_{u,v,w}\big)/\ZZ_2(1,1,0)=\AAA_1$.
Since
$\deg_C\widetilde T'=3$ and $\widetilde T'\bigcap \Sing \widetilde S=\emptyset$
then
$(C,B)=(\PP^1, \frac34P_1+\frac12P_2+\frac34P_3)$. Q.E.D.
\end{example}

\begin{theorem}\label{f4}
Let $f=t^3+f_4(z,x,y)$ gives a canonical singularity
$(X,0)$, where $f_4$ is a homogeneous polynomial of degree 4.
Then ${\bf p}=(4,3,3,3)$-blow-up induces a plt blow-up of
$(X,0)$ and $(E,D)=\big(t+f_4(z,x,y)\subset
\PP(4,1,1,1),\frac23\{t=0\}\big)=(\PP^2,\frac23C)$, where $C=\{f_4=0\}$.
Then $(X,0)$ is exceptional if and only if one of the following
possibilities holds.
\begin{enumerate}
\item The curve $C \subset \PP^2$ is irreducible. Its singular points can be
only ordinary double points.
\item Curve $C \subset \PP^2$ is the union of two irreducible conics
which intersect in 4 distinct points.
\end{enumerate}
\begin{proof} This theorem is proved as theorem \ref{f5}
\end{proof}
\end{theorem}

\section{\bf Classification of three-dimensional exceptional canonical
hypersurface singularities}
\par {\bf The general comment to the tables.}
The classification of exceptional singularities is given in the
tables of the following view.
In the first column the singularity number is written.
In the second column the singularity equation and the
exceptionality condition are written.
Almost in all cases the polynomials, that are grouped according to the similar
log Del Pezzo surface $(E,D)$ (see \ref{ED}), are written.
In the third column corresponding log Del Pezzo surfaces are written.
In the upper row the well-formed hypersurface
$E\subset \PP$ is written. If $E$ is  a linear cone
(i.e. $\tilde d=\tilde a_k$ for some $k$) then
$E=\PP(\tilde a_1,\ldots,\tilde a_{k-1},\tilde a_{k+1},\ldots,\tilde a_n)$ is
written in the second row.
As well the different $D$ is written in the second raw.
The notation is
$\Diff=(\frac{q_1-1}{q_1},\frac{q_2-1}{q_2},\frac{q_3-1}{q_3},\frac{q_4-1}{q_4})$.
In the fourth column the minimal complement index is written.
Let us consider the second column again.
Different cases in it are separated by  a comma or by a semicolon.
If it is a comma then the following singularity has the same minimal
complement index as the previous one.
If it is a semicolon then the following singularity has the next minimal
complement index.
Also in the description of different cases we use usual logical symbols:
$\&$ -- "and", $||$ -- "or". For example, the notation
$|(n=4;5)\ \&\  b\ne 0 \ |\ ...\ |\  2,4| $ means that in the case
$n=4$ and $b\ne 0$ there is
2-complement, but in the case $n=5$ and $b\ne 0$ there is 4-complement.
The symbol $|||$ also means the logical "or", but it is used
to separate the
polynomials (monomials) only. For example, the notation
$h+ (x^8\  ||| \  x^7y )$ means that we have to consider two cases:
$h+x^8$, $h+x^7y$. The exceptionality conditions are also separated by
the symbol $|||$.
\par
Very often the conditions that the singularity is canonical are not
written.
The reader can find them himself. The nonisolated singularity components can be
founded with the help of lemma
\ref{Laz}. By theorem \ref{crqht} the singularity is canonical if and only if
it has Du Val singularities along all components.
As well the requirement that singularity is normal
-- $\codim_{\CC^4}\Sing X\ge 3$ must be taken into account.
For example, the singularity given by
$(t^3-tg(z,x,y)+ag^3(z,x,y)=0,0)$ is nonnormal if
$a=\frac4{3\sqrt{3}}$, $\frac2{3\sqrt{3}}$.
\par
If we know the minimal complement index then the required complement can
be easily found.
\par {\bf Quasihomogeneous transformations.}
Let the quasihomogeneous singularity with weights
${\bf p}$ and degree $d$ is given. It can happen that there exist
a lot of monomials corresponding to this weighting.
On the quasihomogeneous part of the weighting the group of
quasihomogeneous biholomorphic maps naturally acts.
It is denoted by
$H({\bf p})$.
Therefore the classification problem is reduced to  describe the orbits of action
and to find most simple element of every orbit.
The general method that can be applied is given in
\cite{AVG}. The group $H({\bf p})$ is generated by the toric and unipotent parts.
The toric part is the linear transformations of coordinates.
The unipotent part  is the quasihomogeneous coordinate change of form
$\alpha \longmapsto \alpha + h(\beta,\gamma,\delta)$.
In our three-dimensional case the singularities don't have the complex
structures. Therefore it is easy to find which transformations must be
applied. In the tables the most convenient form of singularity is given.
After this we separate to those orbits which have different conditions
on the exceptionality.
\par {\bf Conventions.} According to described above we introduce the
following conventions.

\par (A) If $g(z,x,y)$ doesn't depend on $z$ then
sometimes we use the abridged notation
$g(x,y)$.
In the third column $g(z^{1/q_2},x^{1/q_3},y^{1/q_4})$ is not often given
in details because it takes much place in the table.

\par (B) In the second column the notation
$n=3\ i\le 2,\ 4\ i\le 1$ means the following condition.
If $n=3$ then the condition on the exceptionality is $i\le 2$, but
if $n=4$ then the condition on the exceptionality is $i\le 1$
(let us remember that they have the same complement index).

\par (C) The notation $ i,j\le 1$ is equivalent to $i\le 1 \ \& \ j\le 1$.

\par (D) $ f_n(x,y)$ is a binary form of degree $n$.
In the notation of singularity equation the several binary forms are both
denoted by the same one.
We can consider that they are always different,
i.e. doesn't depend on each other.
\par
Consider notation
$x^{ai}f_{n-i}(x^a,y^b)$. One always assumes that
the coefficient of monomial
$x^{an}$ is nonzero. If there is no coefficient before a binary form then
we always suppose that the coefficient of
$x^{an}$ is nonzero. Let $a=1$ and the remaining part of $f'$ doesn't depend
on $x$ (i.e. $\text{G.C.D}(f',x)=1$).
Then we suppose that the irreducible factors of
$f_{n-i}(x,y^b)$ have the multiplicities at most $i$.

\par (E) Consider notation
$x^i(x+y^b)^jf_{n-i-j}(x,y^b)$. Let
the remaining part of $f'$ depends on
$x$ and doesn't depend on $x+y^b$ (i.e. $\text{G.C.D.}(f',x+y^b)=1$).
Then we suppose that the irreducible factors of
$f_{n-i-j}(x,y^b)$ have the multiplicities at most
$j$.

\par (F) Consider notation
$cx^{ai}(x^a+y^b)^jf_{n-i-j}(x^a,y^b)$. The presence of coefficient $c$
means that the coefficients of binary form are arbitrary.
If $i+j>n$ then $c=0$.
If $i+j\le n$ then $c\ne 0$.
If there is $c\ne 0$ in the exceptionality condition then we always suppose
that the coefficient
$c$ is absent (in other words the coefficient of
$x^{an}$ is nonzero). Let $a=1$. If the remaining part of singularity doesn't
depend on
$x$ then the irreducible factors of $f_{n-i}(x,y^b)$ have the multiplicities
at most $i$. If it depends on $x$, but not on $x+y^b$ then
the irreducible factors have the multiplicities
at most $j$.

\par (I) It can happen that the quasihomogeneous polynomial has
more then one part of type
$\mathcal M_2$, which belongs to Newton line rotation list.
They have the common part of type
$\mathcal M_1$ and differ by some monomials. In the nontrivial cases
their coefficients will be underlined.
Assume that $k$ coefficients are underlined. Then we must consider
$k$ cases. In each one the corresponding coefficient is absent.
It is of a great importance for the binary forms
(see (D)--(F)).

\par (J) If there is a footnote in the first column and in the second raw then
the exceptionality condition is given at the end of the corresponding table.

\par (K) The condition on the exceptionality is written for all
common components.

\begin{example} (1)
$ f=t^2+z^3+\underline{a}zx^if_{5-i}(x,y^2)+
\underline{b}yx^jf_{7-j}(x,y^2).$ The exceptionality condition
$i\le 2 \ ||\ j\le 4$ (see example \ref{ex12}) has the following meaning.
Assume that $x+cy^2$ has the multiplicity $i_1$ in $f_{5-i}(x,y^2)$ and
the multiplicity $j_1$ in $f_{7-j}(x,y^2)$.
Then the exceptionality condition means that
$i_1\le 2 \ ||\ j_1\le 4$.
Recall that if $\underline{a}=0$ or
$\underline{b}=0$ then according to the point (F) one assume that $i>5$
or $j>8$ respectively.
\par
(2) $ f=t^2+z^3x+x^i(x+y^2)^jf_{5-i}(x,y^2).$
The exceptionality condition is
$i\le 1 \ \& \ j\le 3$. Here the remaining part
$t^2+z^3x$ depends on $x$.
Therefore the condition $i\le 1$ is fixed and the requirement
$j\le 3$ as in the previous point means that the multiplicity of any
irreducible factor of
$f_{5-i}(x,y^2)$ is at most 3.

\end{example}
\newpage
\begin{center}
\small{1. Singularity -- $t^2+z^3+g(z,x,y)$}
\end{center}

\begin{center}
\hoffset=-1cm

\tabletail{\hline}
\tablehead{\hline}

\scriptsize{
\begin{supertabular}{|c|c|c|c|}
\hline
No. & $g(z,x,y)$ & log Del Pezzo surface & Index\\
\hline
1& $x^7 + azy^{4n} +by^{6n}$ &
$t^2 + z^3 + x + azy^4+by^6 \subset \PP(3,2,6,1)$  & 7\\
& $n=2,3,4,5,6$  & $\PP(3,2,1)$,\
$\Diff$=\scriptsize{$(0,0,\frac{6}{7},\frac{n-1}{n})$} & \\
\hline
2& $x^7+azy^{2n}+by^{3n}$ &
$ t+ z^3+x+azy^2+by^3 \subset \PP(3,1,3,1) $  &7,8,14 \\
& $n=3;5;9,11,13$ & $\PP(3,1,1)$,\
$\Diff$=\scriptsize{$(\frac{1}{2},0,\frac{6}{7},\frac{n-1}{n})$}&\\
\hline
3& $x^7 + y^{2n},n=4,5;8;$ &
$ t^2 + z + x + y^2 \subset \PP(1,2,2,1) $  &7,9,15,\\
 &
$10,11,13;16,17,19,20$ & $\PP(1,2,1)$,\
$\Diff$=\scriptsize{$(0,\frac{2}{3},\frac{6}{7},\frac{n-1}{n})$} &21\\
\hline
4 & $x^7 + y^{n}$ &$t+z+x+y\subset \PP(1,1,1,1)$
 &12,14,18,
\\
\end{supertabular}
}
\end{center}


\vspace{1cm}

\begin{center}
\small{2. Singularity -- $t^2+z^4+g(z,x,y)$}
\end{center}
\hoffset=-1cm
\begin{center}

\tabletail{\hline}
\tablehead{\hline}
\scriptsize{
\begin{supertabular}{|c|c|c|c|}
\hline
1 & $az^{2}y^n+x^{5}+by^{2n}$ &
$t^2 +z^2 +azy+x+by^2 \subset \PP(1,1,2,1)$  &5,10
\\
& $n=3;7,9$ & $\PP^2$,\
$\Diff$=\scriptsize{$(0,\frac{1}{2},\frac{4}{5},\frac{n-1}{n})$}
&\\
\hline
2& $x^5+y^n$ &
$t^2 +z^2 +x+y \subset \PP(1,1,2,2)$  & 7,10,11,
\\
&$n=7;9;11;13;17,19$ &  $ \PP(1,1,2)$,\
$\Diff$=\scriptsize{$(0,\frac{1}{2},\frac{4}{5},\frac{n-1}{n})$} &
15,20
\\
\hline
3 & $az^{2}y^{2n}+bzy^{3n}+$ &
$t^2 +z^4 +g(z,x^{1/5},y^{1/n})\subset \PP(2,1,4,1)$  &5
\\
& $+cy^{4n}+x^5,n=2,3,4$ & $\PP(2,1,1)$,\
$\Diff$=\scriptsize{$(0,0,\frac{4}{5},\frac{n-1}{n})$} &\\

\hline
4 & $az^2yx^{i}f_{2-i}(x,y^2)+$ &
$ t^2+z^2+g(z^{1/2},x,y) \subset\PP(5,5,2,1)$  &2
\\
& $+x^jf_{5-j}(x,y^2),i\le 1 || j\le 2 $ &
$\Diff$=\scriptsize{$(0,\frac{1}{2},0,0)$}  &  \\
\end{supertabular}
}
\end{center}
\footnotesize{
\par (1) There are two cases.
A). $t^2+z^4+z^2(ax^3+by^4)+zy^2(cx^3+dy^4)+ex^6+kx^3y^4+ly^8$.
Then the singularity is exceptional if the following condition
is satisfied:
$|b|+|d|+|l|\ne 0$\ \& \ $\big(e\ne 0\ || \ (a\ne 0\ \& \
|c|+|k|+|l|\ne 0)\big)$.\\
B). $t^2+z^3y^2+z^2x^3+zy^2(ax^3+by^4)+cx^6+dx^3y^4+ey^8$.
The exceptionality condition is
$|ab|+|c|+|d|+|e|\ne 0\ \& \ |a|+|b|+|d|+|e|\ne 0$ .
}

\begin{center}

\vspace{1cm}

\begin{center}
\small{3. Singularity -- $t^2+z^3x+g(z,x,y)$}
\end{center}

\tabletail{\hline}
\tablehead{\hline}
\scriptsize{
\begin{supertabular}{|c|c|c|c|}
\hline
1 & $azx^2y^{2n}+x^{5}+by^{6n}$ &
$t^2 +z^3x +azx^2y^2+x^5+by^6 \subset \PP(15,8,6,5)$  &4,5
\\
& $n=2\ b\ne 0;3,4$ &
$\Diff$=\scriptsize{$(0,0,0,\frac{n-1}{n})$}
&\\
\hline
2 & $azx^2y^n+x^{5}+by^{3n}$ &
$t+z^3x+azx^2y+x^5+by^3 \subset \PP(15,4,3,5)$  &4,10
\\
& $n=3\ b\ne 0;7,9$ &$\PP(4,3,5)$,\
$\Diff$=\scriptsize{$(\frac{1}{2},0,0,\frac{n-1}{n})$}
&\\
\hline
3& $x^5+y^n$ &
$t+zx+x^5+y \subset \PP(5,4,1,5)$  & 7,12,14,
\\
&$n=7;11;13;17;19;$ &  $ \PP(4,1,5)$,\
$\Diff$=\scriptsize{$(\frac{1}{2},\frac{2}{3},0,\frac{n-1}{n})$} &
18,22,24,
\\
&23;29&&30\\
\hline
4& $x^5+y^{2n}$ &
$t^2+zx+x^5+y^2 \subset \PP(5,8,2,5)$  & 7,9,15
\\
\end{supertabular}
}
\end{center}

\footnotesize{
\par (1) See example 3.22.
\par (2) The singularity is
$t^2+z^3x+\underline{a}zx^i(x+y^2)^jf_{4-i-j}(x,y^2)+bz^2y^5+
\underline{c}x^ky(x+y^2)^lf_{5-k-l}(x,y^2)$. The exceptionality condition is
$b\ne 0\ ||\ \big( b=0 \& (i\le 1 || k\le 2) \& (j\le 2 || l\le 4)\big)$.
\par (3) The singularity is
$t^2+z^3x+\underline{a}zx^iy^j(x+y)^kf_{5-i-j-k}(x,y)+
\underline{b}x^ly^m(x+y)^nf_{7-l-m-n}(x,y)+cz^2y^3$. The common exceptionality
condition for these two cases is
$(c\ne 0) \ ||\ \Big((i=0 || l= 0) \& (k\le 2 || n\le 2)\Big)$.
The first case is $j\le 1$ and the
coefficient $\underline{a}$ is absent (see the comments to the tables).
The second case is $m\le 2$ and the
coefficient  $\underline{b}$ is absent.
}

\begin{center}

\vspace{1cm}

\begin{center}
\small{4. Singularity -- $t^2+z^3y+g(z,x,y)$}
\end{center}

\tabletail{\hline}
\tablehead{\hline}
\scriptsize{
\begin{supertabular}{|c|c|c|c|}
\hline
1& $x^7+y^{2n}$ &
$t^2+zy+x+y^{2n} \subset \PP(n,2n-1,2n,1)$  & 7,9,15,
\\
&$n=4;6;9,10;12,13$ &  $ \PP(n,2n-1,1)$,\
$\Diff$=\scriptsize{$(0,\frac{2}{3},\frac{6}{7},0)$} &
21\\
\hline
2 & $x^7+y^n$ &
$t+zy+x+y^n \subset \PP(n,n-1,n,1)$  &14,16,22
\\
& $n=9;11;15;17;23;27$ &$\PP(n-1,n,1)$,\
$\Diff$=\scriptsize{$(\frac{1}{2},\frac{2}{3},\frac{6}{7},0)$}
&28,36,42\\
\hline
3& $ay^{2n}z^if_{2-i}(z,y^{2n-1})+x^7$ &
$t^2+z^3y+g(z,x^{1/7},y) \subset \PP(3n-1,2n-1,6n-2,1)$  & 7
\\
&$n=2,3,4$ &$\PP(3n-1,2n-1,1)$,\
$\Diff$=\scriptsize{$(0,0,\frac{6}{7},0)$} &
\\
\hline
4& $ay^nz^if_{2-i}(z,y^{n-1})+x^7$ &
$t+z^3y+g(z,x^{1/7},y) \subset \PP(3n-2,n-1,3n-2,1)$  &8,14
\\
&$n=5;7,9$ &$\PP(n-1,3n-2,1)$,\
$\Diff$=\scriptsize{$(\frac{1}{2},0,\frac{6}{7},0)$} &
\\
\end{supertabular}
}


\vspace{1cm}

\begin{center}
\small{5. Singularity -- $t^2+z^5+g(z,x,y)$}
\end{center}

\tabletail{\hline}
\tablehead{\hline}
\scriptsize{
\begin{supertabular}{|c|c|c|c|}
\hline
1 & $x^iy^jf_{6-i-j}(x,y)$ &
$t^2 +z +g(z,x,y) \subset \PP(3,6,1,1)$  &5
\\
&$i\le 2 \& j\le 2$ & $\PP(3,1,1)$,\
$\Diff$=\scriptsize{$(0,\frac{4}{5},0,0)$}
&\\
\hline
2& $x^6+y^7$ &
$t^2 +z +x^2+y \subset \PP(1,2,1,2)$  & 15
\\
& &  $ \PP(1,1,2)$,\
$\Diff$=\scriptsize{$(0,\frac{4}{5},\frac{2}{3},\frac{6}{7})$} &
\\
\hline
3 & $x^6+xy^6$ &
$t^2 +z +x^6+xy\subset \PP(3,6,1,5)$  &25
\\
& & $\PP(3,1,5)$,\
$\Diff$=\scriptsize{$(0,\frac{4}{5},0,\frac{5}{6})$} &\\
\hline
4 & $x^6+zy^5$ &
$ t^2 +z^5 +x^2+zy \subset\PP(5,2,5,8)$  &6
\\
& &$\Diff$=\scriptsize{$(0,0,\frac{2}{3},\frac{4}{5})$}
&  \\
\hline
5 & $x^6+zxy^4$ &
$ t^2 +z^5 +x^6+zxy \subset\PP(15,6,5,19)$  &8
\\
& &$\Diff$=\scriptsize{$(0,0,0,\frac{3}{4})$}
& \\
\hline
6 & $x^6+zx^2y^3$ &
$t^2+z^5+x^6+zx^2y \subset\PP(15,6,5,14)$  &7
\\
& &$\Diff$=\scriptsize{$(0,0,0,\frac{2}{3})$}
&  \\
\hline
7 & $ax^6+bzx^3y^2+z^2y^4$ &
$t^2+z^5+ax^2+bzxy+z^2y^2 \subset\PP(5,2,5,3)$  &9
\\
& &$\Diff$=\scriptsize{$(0,0,\frac{2}{3},\frac{1}{2})$}
&  \\
\hline
8 & $x^6+z^2xy^3$ &
$t^2+z^5+x^6+z^2xy \subset\PP(15,6,5,13)$  &13
\\
& &$\Diff$=\scriptsize{$(0,0,0,\frac{2}{3})$}
&\\
\hline
9 & $az^2xy^3+\underline{b}zx^3y^2+\underline{c}x^5y+$ &
$ t+z^5+g(z,x,y) \subset\PP(35,7,6,5)$  &10
\\
& $dy^7, |a|+|d|\ne 0 $ &$\PP(7,6,5)$,\
$\Diff$=\scriptsize{$(\frac{1}{2},0,0,0)$}  &  \\
\hline
10 & $x^5y+xy^6$ &
$t+z+x^5y+xy^6\subset \PP(29,29,5,4)$  &30
\\
& & $\PP(29,5,4)$,\
$\Diff$=\scriptsize{$(\frac{1}{2},\frac{4}{5},0,0)$}
&\\
\hline
11 & $x^5y+x^2y^5$ &
$t+z+x^5y+x^2y^5\subset \PP(23,23,4,3)$  &30
\\
& & $\PP(23,4,3)$,\
$\Diff$=\scriptsize{$(\frac{1}{2},\frac{4}{5},0,0)$}
&\\
\hline
12 & $x^5y+zy^5$ &
$t+z^5+xy+zy^5\subset \PP(25,5,21,4)$  &30
\\
& & $\PP(5,21,4)$,\
$\Diff$=\scriptsize{$(\frac{1}{2},0,\frac{4}{5},0)$}
&\\
\hline
13 & $x^5y+azy^6+bz^3y^3$ &
$t+z^5+xy+azy^6+bz^3y^3\subset \PP(15,3,13,2)$  &20
\\
& & $\PP(3,13,2)$,\
$\Diff$=\scriptsize{$(\frac{1}{2},0,\frac{4}{5},0)$}
&\\
\hline
14 & $x^5y+zxy^4$ &
$t+z^5+x^5y+zxy^4\subset \PP(95,19,16,15)$  &6
\\
& & $\PP(19,16,15)$,\
$\Diff$=\scriptsize{$(\frac{1}{2},0,0,0)$}
&\\
\hline
15 & $x^5y+zx^2y^3$ &
$t+z^5+x^5y+zx^2y^3\subset \PP(65,13,11,10)$  &14
\\
& & $\PP(13,11,10)$,\
$\Diff$=\scriptsize{$(\frac{1}{2},0,0,0)$}
&\\
\hline
16 & $x^5y+z^2y^4$ &
$t^2+z^5+xy+z^2y^4 \subset\PP(10,4,17,3)$  &10
\\
& &$\Diff$=\scriptsize{$(0,0,\frac{4}{5},0)$}
&\\
\hline
17& $x^4y^2+y^n$ &
$t^2 +z +x^2y^2+y^n \subset \PP(n,2n,n-2,2)$  & 15,20
\\
&$n=7;9$ &  $ \PP(n,n-2,2)$,\
$\Diff$=\scriptsize{$(0,\frac{4}{5},\frac{1}{2},0)$} &
\\
\hline
18& $x^4y^2+ax^2y^5+by^8$ &
$t^2 +z +x^2y^2+axy^5+by^8 \subset \PP(4,8,3,1)$  & 10
\\
& &  $ \PP(4,3,1)$,\
$\Diff$=\scriptsize{$(0,\frac{4}{5},\frac{1}{2},0)$} &
\\
\hline
19 & $x^4y^2+xy^n$ &
$t^2 +z +x^4y^2+xy^n\subset \PP(2n-1,4n-2,n-2,3)$  &10,15
\\
&$n=6;7$ & $\PP(2n-1,n-2,3)$,\
$\Diff$=\scriptsize{$(0,\frac{4}{5},0,0)$} &\\
\hline
20 & $x^4y^2+zy^n$ &
$ t^2 +z^5 +x^2y^2+zy^n \subset\PP(5n,2n,5n-8,8)$  &11,16
\\
&$n=5;7$ &$\Diff$=\scriptsize{$(0,0,\frac{1}{2},0)$}
&  \\
\hline
21 & $x^4y^2+zy^6$ &
$ t^2 +z^5 +x^2y^2+zy^6 \subset\PP(15,6,11,4)$  &8
\\
& &$\Diff$=\scriptsize{$(0,0,\frac{1}{2},0)$}
&  \\
\hline
22& $x^4y^2+zxy^n$ &
$ t^2 +z^5 +g(z,x,y) \subset\PP(10n-5,4n-2,5n-8,11)$  &6,11
\\
&$n=4;5$ &$\Diff=\emptyset$
& \\
\hline
23 & $x^4y^2+az^2y^4+bzx^2y^3$ &
$ t^2 +z^5 +g(z,x^{1/2},y) \subset\PP(10,4,7,3)$  &6
\\
& &$\Diff$=\scriptsize{$(0,0,\frac{1}{2},0)$}
&  \\
\hline
24 & $x^4y^2+z^2y^5$ &
$ t^2 +z^5 +x^2y^2+z^2y^5 \subset\PP(25,10,19,6)$  &12
\\
& &$\Diff$=\scriptsize{$(0,0,\frac{1}{2},0)$}
&  \\
\hline
25& $x^4y^2+z^2xy^3$ &
$ t^2 +z^5 +g(z,x,y) \subset\PP(25,10,9,7)$  &7
\\
& &$\Diff=\emptyset$
& \\
\hline
26& $zx^5+y^7$ &
$t +z^5 +zx+y \subset \PP(5,1,4,5)$  & 50
\\
& &  $ \PP(1,4,5)$,\
$\Diff$=\scriptsize{$(\frac{1}{2},0,\frac{4}{5},\frac{6}{7})$} &
\\
\hline
27 & $zx^5+z^2y^4$ &
$t^2 +z^5 +zx+z^2y^2\subset \PP(5,2,8,3)$  &15
\\
& &
$\Diff$=\scriptsize{$(0,0,\frac{4}{5},\frac{1}{2})$} &\\
\hline
28& $zx^5+z^2xy^3$ &
$ t +z^5 +zx^5+z^2xy \subset\PP(25,5,4,11)$  &22
\\
& &$\PP(5,4,11)$,\
$\Diff$=\scriptsize{$(\frac{1}{2},0,0,\frac{2}{3})$}
&  \\
\hline
29& $zx^4y+y^7$ &
$t^2 +z^5 +zx^2y+y^7 \subset \PP(35,14,23,10)$  & 20
\\
& &
$\Diff$=\scriptsize{$(0,0,\frac{1}{2},0)$} &
\\
\hline
30& $zx^4y+xy^6$ &
$ t +z^5 +zx^4y+xy^6 \subset\PP(115,23,19,16)$  &32
\\
& &$\PP(23,19,16)$,\
$\Diff$=\scriptsize{$(\frac{1}{2},0,0,0)$}
&  \\
\hline
31 & $zx^4y+z^2y^4$ &
$t^2 +z^5 +zxy+z^2y^4\subset \PP(10,4,13,3)$  &12
\\
& &
$\Diff$=\scriptsize{$(0,0,\frac{3}{4},0)$} &\\
\hline
32& $zx^4y+z^2xy^3$ &
$ t +z^5 +zx^4y+z^2xy^3 \subset\PP(55,11,9,8)$  &16
\\
& &$\PP(11,9,8)$,\
$\Diff$=\scriptsize{$(\frac{1}{2},0,0,0)$}
&  \\
\hline
33& $zx^3y^2+y^8$ &
$t^2 +z^5 +zxy+y^4 \subset \PP(10,4,11,5)$  & 15
\\
& &
$\Diff$=\scriptsize{$(0,0,\frac{2}{3},\frac{1}{2})$} &
\\
\hline
34& $zx^3y^2+y^9$ &
$ t +z^5 +zxy^2+y^9 \subset\PP(45,9,26,5)$  &30
\\
& &$\PP(9,26,5)$,\
$\Diff$=\scriptsize{$(\frac{1}{2},0,\frac{2}{3},0)$}
&  \\
\hline
35 & $zx^3y^2+xy^6$ &
$t^2 +z^5 +zx^3y+xy^3\subset \PP(20,8,7,11)$  &11
\\
& &
$\Diff$=\scriptsize{$(0,0,0,\frac{1}{2})$} &\\
\hline
36& $zx^3y^2+xy^7$ &
$ t +z^5 +zx^3y^2+xy^7 \subset\PP(95,19,18,11)$  &22
\\
& &$\PP(19,18,11)$,\
$\Diff$=\scriptsize{$(\frac{1}{2},0,0,0)$}
&  \\
\hline
37& $zx^3y^2+x^2y^5$ &
$ t +z^5 +zx^3y^2+x^2y^5 \subset\PP(55,11,10,7)$  &14
\\
& &$\PP(11,10,7)$,\
$\Diff$=\scriptsize{$(\frac{1}{2},0,0,0)$}
&  \\
\hline
38& $zx^3y^2+z^2y^5$ &
$ t +z^5 +zxy^2+z^2y^5 \subset\PP(25,5,14,3)$  &18
\\
& &$\PP(5,14,3)$,\
$\Diff$=\scriptsize{$(\frac{1}{2},0,\frac{2}{3},0)$}
&  \\
\hline

\end{supertabular}
}


\vspace{1cm}

\begin{center}
\small{6. Singularity -- $t^2+z^4x+g(z,x,y)$}
\end{center}

\tabletail{\hline}
\tablehead{\hline}
\scriptsize{
\begin{supertabular}{|c|c|c|c|}
\hline
1 & $x^iy^j f_{6-i-j}(x,y)$ &
$t^2 +zx +g(z,x,y) \subset \PP(3,5,1,1)$  &4
\\
&$i\le 1\ j\le 2$ &
$\Diff$=\scriptsize{$(0,\frac{3}{4},0,0)$}
&\\
\hline
2& $x^6+y^7$ &
$t^2 +zx +x^6+y \subset \PP(3,5,1,6)$  & 8
\\
& &  $ \PP(3,5,1)$,\
$\Diff$=\scriptsize{$(0,\frac{3}{4},0,\frac{6}{7})$} &
\\
\hline
3& $x^6+x^2y^5$ &
$t^2 +zx +x^6+x^2y \subset \PP(3,5,1,4)$  & 8
\\
& &$\Diff$=\scriptsize{$(0,\frac{3}{4},0,\frac{4}{5})$} &
\\
\hline
4 & $x^6+zy^n$ &
$t^2 +z^4x +x^6+zy\subset \PP(12,5,4,19)$  &5,19
\\
&$n=5;6$ &$\Diff$=\scriptsize{$(0,0,0,\frac{n-1}{n})$} &\\
\hline
5 & $ax^6+bzx^3y^2+z^2y^4$ &
$t^2 +z^4x +ax^6+bzx^3y+z^2y^2\subset \PP(12,5,4,7)$  &7
\\
& &$\Diff$=\scriptsize{$(0,0,0,\frac{1}{2})$} &\\
\hline
6 & $x^5y+y^7$ &
$ t^2+z^2x+x^5y+y^7 \subset\PP(35,29,12,10)$  &8
\\
& &$\Diff$=\scriptsize{$(0,\frac{1}{2},0,0)$}  &  \\
\hline
7 & $x^5y+y^8$ &
$ t^2+zx+x^5y+y^8 \subset\PP(20,33,7,5)$  &20
\\
& &$\Diff$=\scriptsize{$(0,\frac{3}{4},0,0)$}  &  \\
\hline
8 & $x^5y+x^2y^5$ &
$ t^2+z^2x+x^5y+x^2y^5 \subset\PP(23,19,8,6)$  &12
\\
& &$\Diff$=\scriptsize{$(0,\frac{1}{2},0,0)$}  &  \\
\hline
9 & $x^5y+zy^n$ &
$t+z^4x+x^5y+zy^n\subset \PP(20n+1,4n+1,4n-3,16)$  &6,14
\\
&$n=5;6$ & $\PP(4n+1,4n-3,16)$,\
$\Diff$=\scriptsize{$(\frac{1}{2},0,0,0)$}
&\\
\hline
10& $x^5y+z^2y^4$ &
$ t^2 +z^4x +x^5y+z^2y^4 \subset\PP(41,17,14,12)$  &6
\\
& &$\Diff=\emptyset$& \\
\hline
11 & $x^5y+z^3y^3$ &
$t+z^4x+x^5y+z^3y^3\subset \PP(63,13,11,8)$  &16
\\
& & $\PP(13,11,8)$,\
$\Diff$=\scriptsize{$(\frac{1}{2},0,0,0)$}
&\\
\hline
12& $x^4y^2+y^n$ &
$t^2 +zx +x^4y^2+y^n \subset \PP(2n,3n+2,n-2,4)$  & 12,16
\\
&$n=7;9$ &
$\Diff$=\scriptsize{$(0,\frac{3}{4},0,0)$} &
\\
\hline
13& $x^4y^2+y^8$ &
$t^2 +zx +x^4y^2+y^8 \subset \PP(8,13,3,2)$  & 8
\\
& &$\Diff$=\scriptsize{$(0,\frac{3}{4},0,0)$} &
\\
\hline
14& $x^4y^2+zy^n$ &
$ t^2 +z^4x +g(z,x,y) \subset\PP(8n+1,3n+2,4n-6,13)$  &9,13
\\
&$n=5,6;7$ &$\Diff=\emptyset$
& \\
\hline
15& $x^4y^2+azx^2y^3+z^2y^4$ &
$ t^2 +z^4x +g(z,x,y) \subset\PP(17,7,6,5)$  &5
\\
& &$\Diff=\emptyset$
& \\
\hline
16 & $x^4y^2+z^2y^5$ &
$ t^2 +z^2x +x^4y^2+zy^5 \subset\PP(21,17,8,5)$  &10
\\
& &$\Diff$=\scriptsize{$(0,\frac{1}{2},0,0)$}
&  \\
\hline
17& $x^4y^2+z^3y^3$ &
$ t^2 +z^4x +g(z,x,y) \subset\PP(27,11,10,7)$  &7
\\
& &$\Diff=\emptyset$
& \\
\hline
18& $zx^5+y^7$ &
$t +z^4x +zx^5+y \subset \PP(19,4,3,19)$  & 14
\\
& &  $ \PP(4,3,19)$,\
$\Diff$=\scriptsize{$(\frac{1}{2},0,0,\frac{6}{7})$} &
\\
\hline
19& $zx^5+x^2y^5$ &
$t +z^4x +zx^5+x^2y \subset \PP(19,4,3,13)$  & 26
\\
& &  $ \PP(4,3,13)$,\
$\Diff$=\scriptsize{$(\frac{1}{2},0,0,\frac{4}{5})$} &
\\
\hline
20 & $zx^5+z^2y^4$ &
$t^2 +z^4x +zx^5+z^2y^2\subset \PP(19,8,6,11)$  &11
\\
& &
$\Diff$=\scriptsize{$(0,0,0,\frac{1}{2})$} &\\
\hline
21& $zx^4y+y^7$ &
$t +z^4x +zx^4y+y^7 \subset \PP(105,22,17,15)$  & 10
\\
& &$\PP(22,17,15)$,\
$\Diff$=\scriptsize{$(\frac{1}{2},0,0,0)$} &
\\
\hline
22& $zx^4y+x^2y^5$ &
$t +z^4x +zx^4y+x^2y^5 \subset \PP(67,14,11,9)$  & 18
\\
& &$\PP(14,11,9)$,\
$\Diff$=\scriptsize{$(\frac{1}{2},0,0,0)$} &
\\
\hline
23& $zx^4y+z^2y^4$ &
$ t^2 +z^4x +g(z,x,y) \subset\PP(31,13,10,9)$  &9
\\
& &$\Diff=\emptyset$
& \\
\hline
24& $zx^3y^2+y^n$ &
$ t +z^4x +zx^3y^2+y^n \subset\PP(11n,2n+2,3n-8,11)$  &16,22
\\
&$n=7;9$ &$\PP(2n+2,3n-8,11)$,\
$\Diff$=\scriptsize{$(\frac{1}{2},0,0,0)$}
&  \\
\hline
25 & $zx^3y^2+y^8$ &
$t^2 +z^4x +zx^3y+y^4\subset \PP(22,9,8,11)$  &11
\\
& &
$\Diff$=\scriptsize{$(0,0,0,\frac{1}{2})$} &\\
\hline
26& $zx^3y^2+x^2y^5$ &
$ t +z^4x +zx^3y^2+x^2y^5 \subset\PP(39,8,7,5)$  &10
\\
& &$\PP(8,7,5)$,\
$\Diff$=\scriptsize{$(\frac{1}{2},0,0,0)$}
&  \\
\hline
27& $zx^3y^2+z^2y^5$ &
$ t +z^4x +zx^3y^2+z^2y^5 \subset\PP(59,12,11,7)$  &14
\\
& &$\PP(12,11,7)$,\
$\Diff$=\scriptsize{$(\frac{1}{2},0,0,0)$}
&  \\
\hline
\end{supertabular}
}


\vspace{1cm}

\begin{center}
\small{7. Singularity -- $t^2+z^4y+g(z,x,y)$}
\end{center}

\tabletail{\hline}
\tablehead{\hline}
\scriptsize{
\begin{supertabular}{|c|c|c|c|}
\hline
1& $x^6+azx^3y^2+bz^2y^4+$ &
$t^2 +z^4y +g(z,x^{1/3},y) \subset \PP(7,3,7,2)$  & 6
\\
&$+cy^7,\ |b|+|c|\ne 0 $&
$\Diff$=\scriptsize{$(0,0,\frac{2}{3},0)$} &
\\
\hline
2& $x^6+ax^3y^4+by^8$ &
$t^2 +zy +x^2+axy^4+by^8 \subset \PP(4,7,4,1)$  & 12
\\
& &$\Diff$=\scriptsize{$(0,\frac{3}{4},\frac{2}{3},0)$} &
\\
\hline
3& $x^6+xy^n$ &
$t^2 +zy +x^6+xy^n \subset \PP(3n,6n-5,n,5)$  & 8,20
\\
&$n=6;7$ &$\Diff$=\scriptsize{$(0,\frac{3}{4},0,0)$} &
\\
\hline
4& $x^6+azxy^3+bx^2y^5$ &
$t^2 +z^2y +g(z,x,y) \subset \PP(15,13,5,4)$  & 8
\\
& &$\Diff$=\scriptsize{$(0,\frac{1}{2},0,0)$} &
\\
\hline
5 & $x^6+zy^n$ &
$t^2 +z^4y +x^2+zy^n\subset \PP(4n-1,2n-2,4n-1,6)$  &6,9
\\
&$n=5;6$ &$\Diff$=\scriptsize{$(0,0,\frac{2}{3},0)$} &\\
\hline
6& $x^6+zxy^n$ &
$ t^2 +z^4y +g(z,x,y)\subset\PP(12n-3,6n-5,4n-1,14)$  &5,14
\\
&$n=4;5$ &$\Diff=\emptyset$
& \\
\hline
7& $x^6+zx^2y^3$ &
$ t^2 +z^4y +g(z,x,y)\subset\PP(33,14,11,10)$  &5
\\
&&$\Diff=\emptyset$& \\
\hline
8 & $az^3x^2+bz^2xy^3+czx^3y^2+$ &
$t+z^4y+g(z,x,y)\subset\PP(23,5,4,3)$  &6
\\
&$dzy^6+ex^5y+kx^2y^5+lz^4y$ & $\PP(5,4,3)$,\
$\Diff$=\scriptsize{$(\frac{1}{2},0,0,0)$}
&\\
&\scriptsize{can.} \& $|e|+|c|\ne 0$&&\\
\hline
9& $zx^5+az^2y^4+by^7$ &
$t^2 +z^4y +zx+az^2y^4+by^7 \subset \PP(7,3,11,2)$  & 10
\\
& &$\Diff$=\scriptsize{$(0,0,\frac{4}{5},0)$} &
\\
\hline
10& $zx^5+xy^6$ &
$t +z^4y +zx^5+xy^6 \subset \PP(121,26,19,17)$& 24
\\
& &  $ \PP(26,19,17)$,\
$\Diff$=\scriptsize{$(\frac{1}{2},0,0,0)$} &
\\
\hline
11& $zx^5+x^2y^5$ &
$t +z^4y +zx^5+x^2y^5 \subset \PP(51,11,8,7)$& 16
\\
& &  $ \PP(11,8,7)$,\
$\Diff$=\scriptsize{$(\frac{1}{2},0,0,0)$} &
\\
\hline
12& $zx^5+x^3y^4$ &
$t +z^4y +zx^5+x^3y^4 \subset \PP(83,18,13,11)$& 22
\\
& &  $ \PP(18,13,11)$,\
$\Diff$=\scriptsize{$(\frac{1}{2},0,0,0)$} &
\\
\hline
\end{supertabular}
}


\vspace{1cm}

\begin{center}
\small{8. Singularity -- $t^2+z^3x^2+g(z,x,y)$}
\end{center}

\tabletail{\hline}
\tablehead{\hline}
\scriptsize{
\begin{supertabular}{|c|c|c|c|}
\hline
1 & $x^iy^j f_{6-i-j}(x,y)$ &
$t^2 +zx^2 +g(z,x,y) \subset \PP(3,4,1,1)$  &3
\\
&$i\le 1\ j\le 3$ &
$\Diff$=\scriptsize{$(0,\frac{2}{3},0,0)$}
&\\
\hline
2& $x^6+y^7$ &
$t^2 +zx^2 +x^6+y \subset \PP(3,4,1,6)$  & 7
\\
& &  $ \PP(3,4,1)$,\
$\Diff$=\scriptsize{$(0,\frac{2}{3},0,\frac{6}{7})$} &
\\
\hline
3& $ax^6+bx^3y^4+y^8$ &
$t^2 +zx^2 +ax^6+bx^3y+y^2 \subset \PP(3,4,1,3)$  & 9
\\
& &$\Diff$=\scriptsize{$(0,\frac{2}{3},0,\frac{3}{4})$} &
\\
\hline
4& $z^2y^4+bzx^3y^2+cx^6+dxy^6$ &
$t^2 +z^3x^2 +g(z,x,y^{1/2})\subset \PP(9,4,3,5)$  & 5
\\
&$|b|+|c|\ne 0 $&
$\Diff$=\scriptsize{$(0,0,0,\frac{1}{2})$} &\\
\hline
5& $x^6+xy^7$ &
$t^2 +zx^2 +x^6+xy \subset \PP(3,4,1,5)$  & 15
\\
& &$\Diff$=\scriptsize{$(0,\frac{2}{3},0,\frac{6}{7})$} &
\\
\hline
6& $x^6+zy^5$ &
$t^2 +z^3x^2 +x^6+zy \subset \PP(9,4,3,14)$  & 5
\\
& &$\Diff$=\scriptsize{$(0,0,0,\frac{4}{5})$} &
\\
\hline
7& $x^6+zy^6$ &
$t^2 +z^3x^2 +x^6+zy^2 \subset \PP(9,4,3,7)$  & 7
\\
& &$\Diff$=\scriptsize{$(0,0,0,\frac{2}{3})$} &
\\
\hline
8 & $x^6+zxy^n$
& $t^2 +z^3x^2 +x^6+zxy\subset \PP(9,4,3,11)$ &4,11
\\
&$n=4;5$ &$\Diff$=\scriptsize{$(0,0,0,\frac{n-1}{n})$} &\\
\hline
9& $x^5y+y^n$ &
$ t +zx^2 +x^5y+y^n\subset\PP(5n,3n+2,n-1,5)$  &15,30
\\
&$n=7;9$ &$\PP(3n+2,n-1,5)$,\
$\Diff$=\scriptsize{$(\frac{1}{2},\frac{2}{3},0,0)$}
&  \\
\hline
10 & $x^5y+y^8$ &
$ t^2+zx^2+x^5y+y^8 \subset\PP(20,26,7,5)$  &13
\\
& &$\Diff$=\scriptsize{$(0,\frac{2}{3},0,0)$}  &  \\
\hline
11& $x^5y+xy^6$ &
$ t +zx^2 +x^5y+xy^6\subset\PP(29,19,5,4)$  &6
\\
& &$\PP(19,5,4)$,\
$\Diff$=\scriptsize{$(\frac{1}{2},\frac{2}{3},0,0)$}
&  \\
\hline
12& $az^3y^3+bx^5y+cx^3y^4+dxy^7$ &
$t +zx^2 +g(z,x,y)\subset \PP(17,11,3,2)$  & 12
\\
&$|b|+|c|\ne 0,|a|+|d|\ne 0 $&$\PP(11,3,2)$,\
$\Diff$=\scriptsize{$(\frac{1}{2},\frac{2}{3},0,0)$} &\\
\hline
13 & $x^5y+zy^n$ &
$t+z^3x^2+x^5y+zy^n\subset\PP(15n+2,3n+2,3n-2,12)$  &11,24
\\
&$n=5;7$ & $\PP(3n+2,3n-2,12)$,\
$\Diff$=\scriptsize{$(\frac{1}{2},0,0,0)$}
&\\
\hline
14 & $x^5y+zxy^4$ &
$t+z^3x^2+x^5y+zxy^4\subset\PP(59,13,10,9)$  &4
\\
&& $\PP(13,10,9)$,\
$\Diff$=\scriptsize{$(\frac{1}{2},0,0,0)$}
&\\
\hline
15 & $x^5y+zxy^5$ &
$ t^2+z^3x^2+x^5y+zxy^5 \subset\PP(37,16,13,9)$  &9
\\
& &$\Diff=\emptyset$& \\
\hline
16 & $x^5y+z^2y^4$ &
$ t^2+z^3x^2+x^5y+z^2y^4 \subset\PP(32,14,11,9)$  &7
\\
& &$\Diff=\emptyset$& \\
\hline
17 & $x^5y+z^2y^5$ &
$t+z^3x^2+x^5y+z^2y^5\subset\PP(79,17,14,9)$  &18
\\
&& $\PP(17,14,9)$,\
$\Diff$=\scriptsize{$(\frac{1}{2},0,0,0)$}
&\\
\hline
18& $az^4y+bz^2y^4+czx^2y^3+$ &
$t^2 +z^3x^2 +g(z,x^{1/2},y)\subset \PP(7,3,5,2)$  & 4
\\
&$+x^4y^2+dy^7 $&
$\Diff$=\scriptsize{$(0,0,\frac{1}{2},0)$} &\\
\hline
19& $az^3y^3+x^4y^2+bx^2y^5+$ &
$t^2 +zx +g(z^{1/3},x^{1/2},y)\subset \PP(4,5,3,1)$  & 6
\\
&$+cy^8,|a|+|c|\ne 0 $&
$\Diff$=\scriptsize{$(0,\frac{2}{3},\frac{1}{2},0)$} &\\
\hline
20& $x^4y^2+y^9$ &
$t^2 +zx +x^2y^2+y^9 \subset \PP(9,11,7,2)$  & 12
\\
& &$\Diff$=\scriptsize{$(0,\frac{2}{3},\frac{1}{2},0)$} &
\\
\hline
21& $x^4y^2+xy^n$ &
$ t^2 +zx^2 +g(z,x,y)\subset\PP(2n-1,2n+2,n-2,3)$  &6,9
\\
&$n=6;7$ &$\Diff$=\scriptsize{$(0,\frac{2}{3},0,0)$}
& \\
\hline
22& $x^4y^2+zy^n$ &
$ t^2 +z^3x +x^2y^2+zy^n\subset\PP(3n+1,n+2,3n-4,5)$  &7,10
\\
&$n=5;7$ &$\Diff$=\scriptsize{$(0,0,\frac{1}{2},0)$}
& \\
\hline
23& $az^2xy^3+bzy^6+x^4y^2$ &
$ t^2 +z^3x^2 +g(z,x,y)\subset\PP(19,8,7,5)$  &5
\\
& &$\Diff=\emptyset$& \\
\hline
24& $x^4y^2+zxy^n$ &
$ t^2 +z^3x^2 +g(z,x,y)\subset\PP(6n-1,2n+2,3n-4,7)$  &4,7
\\
&$n=4;5$ &$\Diff=\emptyset$& \\
\hline
25 & $x^4y^2+z^2y^5$ &
$ t^2 +z^3x +x^2y^2+z^2y^5 \subset\PP(17,7,13,4)$  &8
\\
& &$\Diff$=\scriptsize{$(0,0,\frac{1}{2},0)$}
&  \\
\hline
26& $x^3y^3+y^n$ &
$ t +zx^2 +x^3y^3+y^n\subset\PP(3n,n+6,n-3,3)$  &10,18
\\
&$n=7;11$ &$\PP(n+6,n-3,3)$,\
$\Diff$=\scriptsize{$(\frac{1}{2},\frac{2}{3},0,0)$}
&  \\
\hline
27 & $x^3y^3+y^{2n}$ &
$ t^2+zx^2+x^3y^3+y^{2n} \subset\PP(3n,2n+6,2n-3,3)$  &7,9
\\
&$n=4;5$ &$\Diff$=\scriptsize{$(0,\frac{2}{3},0,0)$}  &  \\
\hline
28& $az^3y^5+x^3y^3+bxy^8$ &
$ t +zx^2 +g(z,x,y)\subset\PP(21,11,5,2)$  &12
\\
& &$\PP(11,5,2)$,\
$\Diff$=\scriptsize{$(\frac{1}{2},\frac{2}{3},0,0)$}
&  \\
\hline
29& $x^3y^3+zy^n$ &
$ t +z^3x^2 +x^3y^3+zy^n\subset\PP(9n+6,n+6,3n-6,8)$  &8,10
\\
&$n=5;7$ &$\PP(n+6,3n-6,8)$,\
$\Diff$=\scriptsize{$(\frac{1}{2},0,0,0)$}
&  \\
\hline
30& $az^5y+bz^2xy^4+czy^8+x^3y^3$ &
$ t +z^3x^2 +g(z,x,y)\subset\PP(39,7,9,4)$  &8
\\
& &$\PP(7,9,4)$,\
$\Diff$=\scriptsize{$(\frac{1}{2},0,0,0)$}
&  \\
\hline
31& $x^3y^3+zy^9$ &
$ t +z^3x^2 +x^3y+zy^3\subset\PP(29,5,7,8)$  &16
\\
& &$\PP(5,7,8)$,\
$\Diff$=\scriptsize{$(\frac{1}{2},0,0,\frac{2}{3})$}
&  \\
\hline
32& $az^4y+bzxy^4+x^3y^3$ &
$t +z^3x^2 +g(z,x,y)\subset \PP(33,7,6,5)$  & 4
\\
& &$\PP(7,6,5)$,\
$\Diff$=\scriptsize{$(\frac{1}{2},0,0,0)$} &\\
\hline
33& $az^4y^2+bzxy^5+x^3y^3$ &
$ t^2 +z^3x^2 +g(z,x,y)\subset\PP(21,8,9,5)$  &5
\\
& &$\Diff=\emptyset$
& \\
\hline
34& $az^4y^3+bzxy^6+x^3y^3$ &
$ t +z^3x^2 +g(z,x,y^{1/3})\subset\PP(17,3,4,5)$  &10
\\
& &$\PP(3,4,5)$,\
$\Diff$=\scriptsize{$(\frac{1}{2},0,0,\frac{2}{3})$}
&  \\
\hline
35& $z^2y^4+x^3y^3$ &
$ t^2 +z^3x^2 +g(z,x,y)\subset\PP(24,10,9,7)$  &5
\\
& &$\Diff=\emptyset$
& \\
\hline
36 & $x^3y^3+z^2y^n$ &
$t+z^3x^2+g(z,x,y)\subset\PP(9n+12,n+6,3n-3,7)$  &8,14
\\
&$n=5;7$ & $\PP(n+6,3n-3,7)$,\
$\Diff$=\scriptsize{$(\frac{1}{2},0,0,0)$}
&\\
\hline
37 & $x^3y^3+z^2y^6$ &
$t^2+z^3x^2+x^3y+z^2y^2\subset\PP(11,4,5,7)$  & 7
\\
& & $\Diff$=\scriptsize{$(0,0,0,\frac{2}{3})$}
&\\
\hline
38 & $x^iy^j f_{7-i-j}(x,y)$ &
$t +zx^2 +g(z,x,y) \subset \PP(7,5,1,1)$  &6
\\
&$i\le 1\ j\le 4$ &
$\Diff$=\scriptsize{$(\frac{1}{2},\frac{2}{3},0,0)$}
&\\
\hline
39 & $x^7+y^8$ &
$t^2 +zx^2 +x^7+y^2 \subset \PP(7,10,2,7)$  & 12
\\
& &$\Diff$=\scriptsize{$(0,\frac{2}{3},0,\frac{3}{4})$} &
\\
\hline
40& $x^7+xy^7$ &
$t +zx^2 +x^7+xy \subset \PP(7,5,1,6)$  & 36
\\
& &  $ \PP(5,1,6)$,\
$\Diff$=\scriptsize{$(\frac{1}{2},\frac{2}{3},0,\frac{6}{7})$} &
\\
\hline
41& $az^2xy^3+bzy^6+x^7$ &
$t +z^3x^2 +g(z,x,y^{1/3}) \subset \PP(21,5,3,8)$  & 16
\\
& &  $ \PP(5,3,8)$,\
$\Diff$=\scriptsize{$(\frac{1}{2},0,0,\frac{2}{3})$} &
\\
\hline
42& $x^7+zxy^5$ &
$t +z^3x^2 +x^7+zxy \subset \PP(21,5,3,13)$  & 26
\\
& &  $ \PP(5,3,13)$,\
$\Diff$=\scriptsize{$(\frac{1}{2},0,0,\frac{4}{5})$} &
\\
\hline
43& $x^7+z^2y^4$ &
$t^2 +z^3x^2 +x^7+z^2y^2\subset \PP(21,10,6,11)$  & 11
\\
& & $\Diff$=\scriptsize{$(0,0,0,\frac{1}{2})$} &\\
\hline
44 & $x^6y+y^8$ &
$t^2 +zx +x^3y+y^8 \subset \PP(12,17,7,3)$  &18
\\
& &$\Diff$=\scriptsize{$(0,\frac{2}{3},\frac{1}{2},0)$}
&\\
\hline
45& $x^6y+xy^7$ &
$t +zx^2 +x^6y+xy^7 \subset \PP(41,29,6,5)$  & 30
\\
& &  $ \PP(29,6,5)$,\
$\Diff$=\scriptsize{$(\frac{1}{2},\frac{2}{3},0,0)$} &
\\
\hline
46& $az^2xy^3+bzy^6+x^6y$ &
$t +z^3x^2 +g(z,x,y) \subset \PP(55,13,8,7)$  & 14
\\
& &  $ \PP(13,8,7)$,\
$\Diff$=\scriptsize{$(\frac{1}{2},0,0,0)$} &
\\
\hline
47& $x^6y+zxy^5$ &
$t +z^3x^2 +x^6y+zxy^5 \subset \PP(89,21,13,11)$  & 22
\\
& &  $ \PP(21,13,11)$,\
$\Diff$=\scriptsize{$(\frac{1}{2},0,0,0)$} &
\\
\hline
48& $x^6y+z^2y^4$ &
$t^2 +z^3x +x^3y+z^2y^4\subset \PP(19,9,11,5)$  & 10
\\
& & $\Diff$=\scriptsize{$(0,0,\frac{1}{2},0)$} &\\
\hline
49 & $x^5y^2+y^8$ &
$t^2 +zx^2 +x^5y+y^4 \subset \PP(10,14,3,5)$  &15
\\
& &$\Diff$=\scriptsize{$(0,\frac{2}{3},0,\frac{1}{2})$}
&\\
\hline
50& $x^5y^2+xy^7$ &
$t +zx^2 +x^5y^2+xy^7 \subset \PP(33,23,5,4)$  & 24
\\
& &  $ \PP(23,5,4)$,\
$\Diff$=\scriptsize{$(\frac{1}{2},\frac{2}{3},0,0)$} &
\\
\hline
51& $x^5y^2+az^2xy^3+bzy^6$ &
$t +z^3x^2 +g(z,x,y) \subset \PP(47,11,7,6)$  & 12
\\
& &  $ \PP(11,7,6)$,\
$\Diff$=\scriptsize{$(\frac{1}{2},0,0,0)$} &
\\
\hline
52& $x^5y^2+zxy^5$ &
$t +z^3x^2 +x^5y^2+zxy^5 \subset \PP(73,17,11,9)$  & 18
\\
& &  $ \PP(17,11,9)$,\
$\Diff$=\scriptsize{$(\frac{1}{2},0,0,0)$} &
\\
\hline
53& $x^5y^2+z^2y^4$ &
$t^2 +z^3x^2 +x^5y+z^2y^2\subset \PP(17,8,5,9)$  & 9
\\
& & $\Diff$=\scriptsize{$(0,0,0,\frac{1}{2})$} &\\
\hline
54 & $x^4y^3+y^8$ &
$t^2 +zx +x^2y^3+y^8 \subset \PP(8,11,5,2)$  &12
\\
& &$\Diff$=\scriptsize{$(0,\frac{2}{3},\frac{1}{2},0)$}
&\\
\hline
55& $x^4y^3+xy^7$ &
$t +zx^2 +x^4y^3+xy^7 \subset \PP(25,17,4,3)$  & 18
\\
& &  $ \PP(17,4,3)$,\
$\Diff$=\scriptsize{$(\frac{1}{2},\frac{2}{3},0,0)$} &
\\
\hline
56& $az^2xy^3+bzx^5+czy^6+dx^4y^3$ &
$t +z^3x^2 +g(z,x,y^{1/3})\subset \PP(13,3,2,5)$  & 10
\\
&$|b|+|d|\ne 0,|a|+|c|\ne 0 $&$\PP(3,2,5)$,\
$\Diff$=\scriptsize{$(\frac{1}{2},0,0,\frac{2}{3})$} &\\
\hline
57& $x^4y^3+zxy^5$ &
$t +z^3x^2 +g(z,x,y) \subset \PP(57,13,9,7)$  & 14
\\
& &  $ \PP(13,9,7)$,\
$\Diff$=\scriptsize{$(\frac{1}{2},0,0,0)$} &
\\
\hline
58& $x^4y^3+z^2y^4$ &
$t^2 +z^3x^2 +x^2y^3+z^2y^4\subset \PP(15,7,9,4)$  & 8
\\
& & $\Diff$=\scriptsize{$(0,0,\frac{1}{2},0)$} &\\
\hline
59& $x^3y^4+y^9$ &
$t +zx^2 +x^3y^4+y^9 \subset \PP(27,17,5,3)$  & 18
\\
& &  $ \PP(17,5,3)$,\
$\Diff$=\scriptsize{$(\frac{1}{2},\frac{2}{3},0,0)$} &
\\
\hline
60& $x^3y^4+zy^7$ &
$t +z^3x^2 +x^3y^4+zy^7 \subset \PP(71,15,13,8)$  & 16
\\
& &  $ \PP(15,13,8)$,\
$\Diff$=\scriptsize{$(\frac{1}{2},0,0,0)$} &
\\
\hline
61& $az^4y+bzxy^5+x^3y^4$ &
$t +z^3x^2 +g(z,x,y) \subset \PP(41,9,7,5)$  & 10
\\
& &  $ \PP(9,7,5)$,\
$\Diff$=\scriptsize{$(\frac{1}{2},0,0,0)$} &
\\
\hline
62& $x^3y^4+z^2y^5$ &
$t +z^3x^2 +x^3y^4+z^2y^5 \subset \PP(61,13,11,7)$  & 14
\\
& &  $ \PP(13,11,7)$,\
$\Diff$=\scriptsize{$(\frac{1}{2},0,0,0)$} &
\\
\hline
63& $zx^5+y^7$ &
$t +z^3x^2 +zx^5+y \subset \PP(13,3,2,13)$  & 14
\\
& &  $ \PP(3,2,13)$,\
$\Diff$=\scriptsize{$(\frac{1}{2},0,0,\frac{6}{7})$} &
\\
\hline
64& $zx^5+y^8$ &
$t^2 +z^3x^2 +zx^5+y^2\subset \PP(13,6,4,13)$  & 13
\\
& & $\Diff$=\scriptsize{$(0,0,0,\frac{3}{4})$} &\\
\hline
65& $zx^5+xy^6$ &
$t^2 +z^3x^2 +zx^5+xy^2\subset \PP(13,6,4,11)$  & 6
\\
& & $\Diff$=\scriptsize{$(0,0,0,\frac{2}{3})$} &\\
\hline
66& $zx^5+xy^7$ &
$t +z^3x^2 +zx^5+xy \subset \PP(13,3,2,11)$  & 22
\\
& &  $ \PP(3,2,11)$,\
$\Diff$=\scriptsize{$(\frac{1}{2},0,0,\frac{6}{7})$} &
\\
\hline
67& $z^2y^4+azx^5+bx^3y^4$ &
$t^2 +z^3x^2 +g(z,x,y^{1/2})\subset \PP(13,6,4,7)$  & 7
\\
& & $\Diff$=\scriptsize{$(0,0,0,\frac{1}{2})$} &\\
\hline
68& $zx^4y+y^7$ &
$ t^2 +z^3x^2 +g(z,x,y)\subset\PP(35,16,11,10)$  &5
\\
& &$\Diff=\emptyset$
& \\
\hline
69& $zx^4y+y^8$ &
$t^2 +z^3x +zx^2y+y^8\subset \PP(20,9,13,5)$  & 10
\\
& & $\Diff$=\scriptsize{$(0,0,\frac{1}{2},0)$} &\\
\hline
70& $zx^4y+xy^n$ &
$t +z^3x^2 +zx^4y+xy^n \subset\PP(10n-3,2n+1,2n-3,8)$  & 6,16
\\
&$n=6;7$ &$\PP(2n+1,2n-3,8)$,\
$\Diff$=\scriptsize{$(\frac{1}{2},0,0,0)$} &
\\
\hline
71& $z^2y^4+zx^4y+ax^2y^5$ &
$t^2 +z^3x +g(z,x^{1/2},y) \subset\PP(11,5,7,3)$  & 6
\\
& &$\Diff$=\scriptsize{$(0,0,\frac{1}{2},0)$} &
\\
\hline
72& $az^2xy^3+bzx^4y+czy^6+dx^3y^4$ &
$t +z^3x^2 +g(z,x,y)\subset \PP(31,7,5,4)$  & 8
\\
&$|b|+|d|\ne 0 \& |a|+|c|\ne 0 $ & $\PP(7,5,4)$,\
$\Diff$=\scriptsize{$(\frac{1}{2},0,0,0)$} &\\
\hline
73& $zx^3y^2+y^n$ &
$ t +z^3x^2 +zx^3y^2+y^n\subset\PP(7n,n+4,2n-6,7)$  &10,14
\\
&$n=7;9$ &$\PP(n+4,2n-6,7)$,\
$\Diff$=\scriptsize{$(\frac{1}{2},0,0,0)$}
&  \\
\hline
74 & $zx^3y^2+y^8$ &
$t^2 +z^3x^2 +zx^3y+y^4\subset\PP(14,6,5,7)$  &7
\\
& &
$\Diff$=\scriptsize{$(0,0,0,\frac{1}{2})$} &\\
\hline
75& $az^2y^5+zx^3y^2+bxy^7$ &
$ t +z^3x^2 +g(z,x,y)\subset\PP(43,9,8,5)$  &10
\\
& &$\PP(9,8,5)$,\
$\Diff$=\scriptsize{$(\frac{1}{2},0,0,0)$}
&  \\
\hline

\end{supertabular}
}


\vspace{1cm}

\begin{center}
\small{9. Singularity -- $t^2+z^3y^2+g(z,x,y)$}
\end{center}

\tabletail{\hline}
\tablehead{\hline}
\scriptsize{
\begin{supertabular}{|c|c|c|c|}
\hline
1& $x^7+y^n$ &
$t +zy^2 +x+y^n \subset \PP(n,n-2,n,1)$  & 28,42
\\
&$n=9;13$ &  $ \PP(n-2,n,1)$,\
$\Diff$=\scriptsize{$(\frac{1}{2},\frac{2}{3},\frac{6}{7},0)$} &
\\
\hline
2& $x^7+y^{10}$ &
$t^2 +zy^2 +x+y^{10} \subset \PP(5,8,10,1)$  & 15
\\
& &  $ \PP(5,8,1)$,\
$\Diff$=\scriptsize{$(0,\frac{2}{3},\frac{6}{7},0)$} &
\\
\hline
3& $x^7+ay^5z^if_{2-i}(z,y^3)$ &
$t +z^3y^2 +g(z,x^{1/7},y) \subset \PP(11,3,11,1)$  & 14
\\
&$i\le 1$ &  $ \PP(3,11,1)$,\
$\Diff$=\scriptsize{$(\frac{1}{2},0,\frac{6}{7},0)$} &
\\
\hline
4& $x^7+y^{12}$ &
$t^2 +zy +x+y^6 \subset \PP(3,5,6,1)$  & 21
\\
& &  $ \PP(3,5,1)$,\
$\Diff$=\scriptsize{$(0,\frac{2}{3},\frac{6}{7},\frac{1}{2})$} &
\\
\hline
5& $x^7+xy^n$ &
$ t +zy^2 +x^7+xy^n\subset\PP(7n,7n-12,n,6)$  &22,36
\\
&$n=7;11$ &$\PP(7n-12,n,6)$,\
$\Diff$=\scriptsize{$(\frac{1}{2},\frac{2}{3},0,0)$}
&  \\
\hline
6& $x^{3i+1}f_{2-i}(x^3,y^4)$ &
$t^2 +zy +g(x,y^{1/2}) \subset \PP(7,11,2,3)$  & 9
\\
&$i\le 1$ & $\Diff$=\scriptsize{$(0,\frac{2}{3},0,\frac{1}{2})$} &
\\
\hline
7& $x^{2i+1}f_{3-i}(x^2,y^3)$ & $t +zy^2 +g(x,y) \subset
\PP(21,17,3,2)$  & 12
\\
&$i\le 2$ &$\PP(17,3,2)$,\
$\Diff$=\scriptsize{$(\frac{1}{2},\frac{2}{3},0,0)$} &
\\
\hline
8& $x^{3i+1}f_{2-i}(x^3,y^5)$ &
$t +zy^2 +g(x,y) \subset \PP(35,29,5,3)$  & 18
\\
&$i\le 1$ &$\PP(29,5,3)$,\
$\Diff$=\scriptsize{$(\frac{1}{2},\frac{2}{3},0,0)$} &
\\
\hline
9& $x^7+x^2y^6$ &
$t^2 +zy^2 +x^7+x^2y^6 \subset \PP(21,32,6,5)$  & 9
\\
& &$\Diff$=\scriptsize{$(0,\frac{2}{3},0,0)$} &
\\
\hline
10& $azx^3y^3+x^7+bx^2y^7$ &
$ t +z^3y^2 +g(z,x,y)\subset\PP(49,13,7,5)$  & 10
\\
& &$\PP(13,7,5)$,\
$\Diff$=\scriptsize{$(\frac{1}{2},0,0,0)$}
&  \\
\hline
11& $x^7+x^2y^8$ &
$t^2 +zy +x^7+x^2y^4 \subset \PP(14,23,4,5)$  & 12
\\
& &$\Diff$=\scriptsize{$(0,\frac{2}{3},0,\frac{1}{2})$} &
\\
\hline
12& $x^7+x^2y^9$ &
$t +zy^2 +x^7+x^2y^9 \subset \PP(63,53,9,5)$  & 30
\\
& &$\PP(53,9,5)$,\
$\Diff$=\scriptsize{$(\frac{1}{2},\frac{2}{3},0,0)$} &
\\
\hline
13& $bzx^2y^4+x^7+cx^3y^5$ &
$t +z^3y^2+g(z,x,y)\subset \PP(35,9,5,4)$  & 8
\\
&&$\PP(9,5,4)$,\
$\Diff$=\scriptsize{$(\frac{1}{2},0,0,0)$} &\\
\hline
14 & $x^7+x^3y^7$ &
$ t+zy^2+x^7+x^3y^7\subset\PP(49,41,7,4)$  &24
\\
&&$\PP(41,7,4)$,\
$\Diff$=\scriptsize{$(\frac{1}{2},\frac{2}{3},0,0)$} &\\
\hline
15 & $x^7+zy^n$ &
$t+z^3y^2+x+zy^n\subset\PP(3n-2,n-2,3n-2,2)$  &22,28
\\
&$n=7;9$& $\PP(n-2,3n-2,2)$,\
$\Diff$=\scriptsize{$(\frac{1}{2},0,\frac{6}{7},0)$}
&\\
\hline
16& $x^7+zxy^n$ &
$ t +z^3y^2 +x^7+zxy^n\subset\PP(21n-14,7n-12,3n-2,11)$  &16,22
\\
&$n=5;7$ &$\PP(7n-12,3n-2,11)$,\
$\Diff$=\scriptsize{$(\frac{1}{2},0,0,0)$}
&  \\
\hline
17 & $x^7+zxy^6$ &
$ t^2+z^3y+x^7+zxy^3 \subset\PP(28,15,8,11)$  &11
\\
& &$\Diff$=\scriptsize{$(0,0,0,\frac{1}{2})$}  &  \\
\hline
18& $x^7+zx^2y^5$ &
$ t +z^3y^2 +x^7+zx^2y^5\subset\PP(91,25,13,8)$  &16
\\
& &$\PP(25,13,8)$,\
$\Diff$=\scriptsize{$(\frac{1}{2},0,0,0)$}
&  \\
\hline
19& $zx^5+y^n$ &
$t +z^3y^2 +zx+y^n \subset \PP(3n,n-2,2n+2,3)$  & 15,30
\\
&$n=7;13$ &  $ \PP(n-2,2n+2,3)$,\
$\Diff$=\scriptsize{$(\frac{1}{2},0,\frac{4}{5},0)$} &
\\
\hline
20& $zx^5+czx^2y^4+$ &
$ t +z^3y^2 +g(z,x,y)\subset\PP(27,7,4,3)$  &6
\\
&$+ex^6y+dx^3y^5+ly^9$ &$\PP(7,4,3)$,\
$\Diff$=\scriptsize{$(\frac{1}{2},0,0,0)$}
&  \\
&$|c|+|d|+|l|\ne 0$&&\\
\hline
21& $zx^5+y^{10}$ &
$t^2 +z^3y^2 +zx+y^{10} \subset \PP(15,8,22,3)$  & 11
\\
& &$\Diff$=\scriptsize{$(0,0,\frac{4}{5},0)$} &
\\
\hline
22& $by^5z^if_{3-i}(z,y^3)+zx^5+$ &
$t +z^3y^2 +g(z,x^{1/5},y) \subset \PP(11,3,8,1)$  & 10
\\
&$+ax^5y^3,i\le 2\ || a\ne 0$ &  $ \PP(3,8,1)$,\
$\Diff$=\scriptsize{$(\frac{1}{2},0,\frac{4}{5},0)$} &
\\
\hline
23 & $zx^5+y^{12}$ &
$t^2 +z^3y^2 +zx+y^6 \subset \PP(9,5,13,3)$  &15
\\
& &$\Diff$=\scriptsize{$(0,0,\frac{4}{5},\frac{1}{2})$}
&\\
\hline
24 & $zx^5+xy^{4n+2}$ &
$t^2 +z^3y +g(z,x,y^{1/2}) \subset \PP(15n+8,10n+1,4n+3,13)$  &9,13
\\
&$n=1;2$ &$\Diff$=\scriptsize{$(0,0,0,\frac{1}{2})$}
&\\
\hline
25& $zx^5+xy^n$ &
$t +z^3y^2 +zx^5+xy^n \subset \PP(15n+2,5n-8,2n+2,13)$  & 16,22,26
\\
&$n=7;9;11$ &  $ \PP(5n-8,2n+2,13)$,\
$\Diff$=\scriptsize{$(\frac{1}{2},0,0,0)$} &
\\
\hline
26& $zx^5+xy^8$ &
$ t^2 +z^3y^2+g(z,x,y)\subset\PP(61,32,18,13)$  &9
\\
& &$\Diff=\emptyset$
& \\
\hline
27& $zx^5+x^2y^n$ &
$t +z^3y^2 +zx^5+x^2y^n \subset \PP(15n+4,5n-6,2n+2,11)$  & 12,18,22
\\
&$n=5;7;9$ &  $ \PP(5n-6,2n+2,11)$,\
$\Diff$=\scriptsize{$(\frac{1}{2},0,0,0)$} &
\\
\hline
28& $zx^5+x^2y^6$ &
$ t^2 +z^3y^2+g(z,x,y)\subset\PP(47,24,14,11)$  &7
\\
& &$\Diff=\emptyset$
& \\
\hline
29 & $zx^5+x^2y^8$ &
$t^2 +z^3y +zx^5+x^2y^4 \subset \PP(31,17,9,11)$  &11
\\
& &$\Diff$=\scriptsize{$(0,0,0,\frac{1}{2})$}
&\\
\hline
30& $zx^5+x^3y^4$ &
$ t^2 +z^3y^2+g(z,x,y)\subset\PP(33,16,10,9)$  &5
\\
& &$\Diff=\emptyset$
& \\
\hline
31 & $zx^5+x^3y^6$ &
$t^2 +z^3y +zx^5+x^3y^3 \subset \PP(24,13,7,9)$  &9
\\
& &$\Diff$=\scriptsize{$(0,0,0,\frac{1}{2})$}
&\\
\hline
32& $zx^5+x^3y^7$ &
$t +z^3y^2 +zx^5+x^3y^7 \subset \PP(111,31,16,9)$  & 18
\\
& &  $ \PP(31,16,9)$,\
$\Diff$=\scriptsize{$(\frac{1}{2},0,0,0)$} &
\\
\hline
33& $zx^5+x^4y^n$ &
$t +z^3y^2 +zx^5+x^4y^n \subset \PP(15n+8,5n-2,2n+2,7)$  & 10,14
\\
&$n=3;5$ &  $ \PP(5n-2,2n+2,7)$,\
$\Diff$=\scriptsize{$(\frac{1}{2},0,0,0)$} &
\\
\hline
34 & $zx^5+x^4y^4$ &
$t^2 +z^3y +zx^5+x^4y^2 \subset \PP(17,9,5,7)$  &7
\\
& &$\Diff$=\scriptsize{$(0,0,0,\frac{1}{2})$}
&\\
\hline

\end{supertabular}
}


\vspace{1cm}

\begin{center}
\small{10. Singularity -- $t^2+z^3xy+g(z,x,y)$}
\end{center}

\tabletail{\hline}
\tablehead{\hline}
\scriptsize{
\begin{supertabular}{|c|c|c|c|}
\hline
1 & $x^5+y^{2n}$ &
$t^2 +zxy +g(x,y) \subset \PP(5n,8n-5,2n,5)$  &7,9,15
\\
&$n=3;6;8,9$ &
$\Diff$=\scriptsize{$(0,\frac{2}{3},0,0)$}
&\\
\hline
2& $x^5+y^n$ &
$ t +zxy +x^5+y^n\subset\PP(5n,4n-5,n,5)$  &10,14,22
\\
&$n=7;9;13;19$ &$\PP(4n-5,n,5)$,\
$\Diff$=\scriptsize{$(\frac{1}{2},\frac{2}{3},0,0)$}
&30  \\
\hline
3 & $azx^2y^{2n+1}+x^5+by^{6n+2}$ &
$ t^2+z^3xy+g(z,x,y) \subset\PP(15n+5,8n+1,6n+2,5)$  &4,5
\\
&$n=1\ b\ne 0; 2$ &$\Diff=\emptyset$& \\
\hline
4 & $x^i(x+y^2)^jf_{5-i-j}(x,y^2)$ &
$t^2 +zxy +g(x,y) \subset \PP(5,7,2,1)$  &3
\\
&$i\le 2\& j\le 4$ &
$\Diff$=\scriptsize{$(0,\frac{2}{3},0,0)$}
&\\
\hline
5& $azx^2y^{2n}+x^5+by^{6n-1}$ &
$ t +z^3xy +g(z,x,y)\subset\PP(30n-5,8n-3,6n-1,5)$  &10
\\
&$n=2;3$ &$\PP(8n-3,6n-1,5)$,\
$\Diff$=\scriptsize{$(\frac{1}{2},0,0,0)$}
&  \\
\hline
6 & $x^if_{5-i}(x,y^3)$ &
$t +zxy +g(x,y) \subset \PP(15,11,3,1)$  &6
\\
&$i\le 3$ &$\PP(11,3,1)$,\
$\Diff$=\scriptsize{$(\frac{1}{2},\frac{2}{3},0,0)$}
&\\
\hline
7 & $x^5+x^2y^{2n}$ &
$t^2 +zxy +g(x,y) \subset \PP(5n,8n-3,2n,3)$  &9
\\
&$n=4,5$ &
$\Diff$=\scriptsize{$(0,\frac{2}{3},0,0)$}
&\\
\hline
8& $x^5+x^2y^n$ &
$ t +zxy +g(x,y)\subset\PP(5n,4n-3,n,3)$  &10,18
\\
&$n=5,7;11$ &$\PP(4n-3,n,3)$,\
$\Diff$=\scriptsize{$(\frac{1}{2},\frac{2}{3},0,0)$}
& \\
\hline
9 & $x^5+zy^{2n+1}$ &
$ t^2+z^3xy+g(z,x,y) \subset\PP(15n+5,8n-1,6n+2,11)$  &3,5,11
\\
&$n=2;3;4,5,6$ &$\Diff=\emptyset$& \\
\hline
10& $x^5+zy^n$ &
$ t +z^3xy +g(z,x,y)\subset\PP(15n-5,4n-5,3n-1,11)$  &8,14,22
\\
&$n=6;8,10;12,14$ &$\PP(4n-5,3n-1,11)$,\
$\Diff$=\scriptsize{$(\frac{1}{2},0,0,0)$}
&  \\
\hline
11 & $x^5+z^2y^{2n}$ &
$ t^2+z^3xy+g(z,x,y) \subset\PP(15n-5,8n-5,6n-2,7)$  &6,7
\\
&$n=2;3,4$ &$\Diff=\emptyset$& \\
\hline
12& $x^5+z^2y^n$ &
$ t +z^3xy +g(z,x,y)\subset\PP(15n-10,4n-5,3n-2,7)$  &6,14
\\
&$n=5;7,9$ &$\PP(4n-5,3n-2,7)$,\
$\Diff$=\scriptsize{$(\frac{1}{2},0,0,0)$}
&  \\
\hline
13& $x^iy^j(x+y)^lf_{6-i-j-l}(x,y)$ &
$t^2 +zxy +g(x,y) \subset \PP(3,4,1,1)$  &3
\\
&$i\le 1\& j\le 1\& l\le 3$ &
$\Diff$=\scriptsize{$(0,\frac{2}{3},0,0)$}
&\\
\hline
14 & $x^6+y^n$ &
$t^2 +zxy +g(x,y) \subset \PP(3n,5n-6,n,6)$  &10,18
\\
&$n=7;11$ &
$\Diff$=\scriptsize{$(0,\frac{2}{3},0,0)$}
&\\
\hline
15& $x^{3i}f_{2-i}(x^3,y^n)$ &
$t^2 +zxy +g(x,y) \subset \PP(3n,5n-3,n,3)$  &6,9
\\
&$n=4\ i=0;5$ &
$\Diff$=\scriptsize{$(0,\frac{2}{3},0,0)$}
&\\
\hline
16& $x^{2i}f_{3-i}(x^2,y^3)$ &
$t^2 +zxy +g(x,y) \subset \PP(9,13,3,2)$  &6
\\
&$ i\le 1$ &
$\Diff$=\scriptsize{$(0,\frac{2}{3},0,0)$}
&\\
\hline
17 & $az^2y^4+bzx^3y^2+x^6+$ &
$ t^2+z^3xy+g(z,x,y) \subset\PP(15,7,5,4)$  &4
\\
&$+cx^2y^5,|a|+|c|\ne 0$ &$\Diff=\emptyset$& \\
\hline
18 & $x^6+x^2y^7$ &
$t^2 +zxy +g(x,y) \subset \PP(21,31,7,4)$  &12
\\
& &$\Diff$=\scriptsize{$(0,\frac{2}{3},0,0)$}
&\\
\hline
19 & $x^6+zy^n$ &
$ t^2+z^3xy+g(z,x,y) \subset\PP(9n-3,5n-6,3n-1,13)$  &7,8,10
\\
&$n=5;6;7;8$ & $\Diff=\emptyset$ &13 \\
\hline
20 & $x^6+zx^2y^4$ &
$ t^2+z^3xy+g(z,x,y) \subset\PP(33,16,11,7)$  &7
\\
& & $\Diff=\emptyset$ & \\
\hline
21 & $x^6+z^2y^5$ &
$ t^2+z^3xy+g(z,x,y) \subset\PP(39,19,13,8)$  &8
\\
& & $\Diff=\emptyset$ & \\
\hline
22& $x^iy^j(x+y)^lf_{7-i-j-l}(x,y)$ &
$t +zxy +g(x,y) \subset \PP(7,5,1,1)$  &6
\\
&$i\le 2\& j\le 2\& l\le 4$ &$\PP(5,1,1)$,\
$\Diff$=\scriptsize{$(\frac{1}{2},\frac{2}{3},0,0)$}
&\\
\hline
23 & $x^7+y^8$ &
$t^2 +zxy +x^7+y^8 \subset \PP(28,41,8,7)$  &12
\\
& &$\Diff$=\scriptsize{$(0,\frac{2}{3},0,0)$}
&\\
\hline
24& $x^7+y^9$ &
$t +zxy +x^7+y^9 \subset \PP(63,47,9,7)$  &42
\\
& &$\PP(47,9,7)$,\
$\Diff$=\scriptsize{$(\frac{1}{2},\frac{2}{3},0,0)$}
&\\
\hline
25& $x^7+x^2y^6$ &
$t^2 +zxy +x^7+x^2y^6 \subset \PP(21,31,6,5)$  &15
\\
& &$\Diff$=\scriptsize{$(0,\frac{2}{3},0,0)$}
&\\
\hline
26& $x^7+zy^6$ &
$ t +z^3xy +x^7+zy^6\subset\PP(119,29,17,15)$  &10
\\
& &$\PP(29,17,15)$,\
$\Diff$=\scriptsize{$(\frac{1}{2},0,0,0)$}
&  \\
\hline
27& $x^7+zx^2y^4$ &
$ t +z^3xy +x^7+zx^2y^4\subset\PP(77,19,11,9)$  &18
\\
& &$\PP(19,11,9)$,\
$\Diff$=\scriptsize{$(\frac{1}{2},0,0,0)$}
&  \\
\hline
28 & $x^7+z^2y^4$ &
$ t^2+z^3xy+x^7+z^2y^4 \subset\PP(35,17,10,9)$  &9
\\
& &$\Diff=\emptyset$& \\
\hline
29 & $x^5y^2+y^{2n}$ &
$t^2 +zxy +x^5y^2+y^{2n} \subset \PP(5n,8n-3,2n-2,5)$  &9,15
\\
&$n=4;5$ &$\Diff$=\scriptsize{$(0,\frac{2}{3},0,0)$}
&\\
\hline
30& $azx^2y^4+x^5y^2+by^9$ &
$ t +z^3xy +g(z,x,y)\subset\PP(45,11,7,5)$  &10
\\
& &$\PP(11,7,5)$,\
$\Diff$=\scriptsize{$(\frac{1}{2},0,0,0)$}
&  \\
\hline
31& $x^5y^2+y^{11}$ &
$t +zxy +x^5y^2+y^{11} \subset \PP(55,41,9,5)$  &30
\\
& &$\Diff$=\scriptsize{$(\frac{1}{2},\frac{2}{3},0,0)$}
&\\
\hline
32 & $x^5y^2+x^2y^6$ &
$t^2 +zxy +x^5y^2+x^2y^6 \subset \PP(13,19,4,3)$  &9
\\
& &$\Diff$=\scriptsize{$(0,\frac{2}{3},0,0)$}
&\\
\hline
33& $x^5y^2+x^2y^7$ &
$t +zxy +x^5y^2+x^2y^7 \subset \PP(31,23,5,3)$  &18
\\
& &$\PP(23,5,3)$,\
$\Diff$=\scriptsize{$(\frac{1}{2},\frac{2}{3},0,0)$}
&\\
\hline
34& $x^5y^2+zy^n$ &
$ t +z^3xy +g(z,x,y)\subset\PP(15n-3,4n-3,3n-5,11)$  &16,22
\\
&$n=6;8$ &$\PP(4n-3,3n-5,11)$,\
$\Diff$=\scriptsize{$(\frac{1}{2},0,0,0)$}
&  \\
\hline
35 & $x^5y^2+zy^7$ &
$ t^2+z^3xy+x^5y^2+zy^7 \subset\PP(51,25,16,11)$  &11
\\
& &$\Diff=\emptyset$& \\
\hline
36 & $x^5y^2+z^2y^4$ &
$ t^2+z^3xy+x^5y^2+z^2y^4 \subset\PP(27,13,8,7)$  &7
\\
& &$\Diff=\emptyset$& \\
\hline
37& $x^5y^2+z^2y^5$ &
$ t +z^3xy +x^5y^2+z^2y^5\subset\PP(69,17,11,7)$  &14
\\
& &$\PP(17,11,7)$,\
$\Diff$=\scriptsize{$(\frac{1}{2},0,0,0)$}
&  \\
\hline
38 & $zx^5+y^n$ &
$ t^2+z^3xy+zx^5+y^n \subset\PP(7n,4n-5,2n+1,14)$  &5,7,14
\\
&$n=7;8;9$ &$\Diff=\emptyset$& \\
\hline
39 & $zx^5+x^2y^n$ &
$ t^2+z^3xy+zx^5+x^2y^n \subset\PP(7n+1,4n-3,2n+1,10)$  &5,6
\\
&$n=5;6$ &$\Diff=\emptyset$& \\
\hline
40 & $z^2y^4+zx^5+bx^4y^3$ &
$ t^2+z^3xy+g(z,x,y) \subset\PP(23,11,7,6)$  &6
\\
& &$\Diff=\emptyset$& \\
\hline

\end{supertabular}
}


\vspace{1cm}

\begin{center}
\small{11. Singularity --
$t^2+f(z,x)+g(z,x,y)=t^2+z^ix^jf_{5-i-j}(z,x)+g(z,x,y),(i\le j)$}
\end{center}

\tabletail{\hline}
\tablehead{\hline}
\scriptsize{
\begin{supertabular}{|c|c|c|c|}
\hline
1 & $y^{2n}$ &
$t^2 +f(z,x)+y^2 \subset \PP(5,2,2,5)$  &3,5
\\
& $n=3\ j=1;4\  j\le 2$ &
$\Diff$=\scriptsize{$(0,0,0,\frac{n-1}{n})$}
&\\
\hline
2 & $y^n$ &
$t+f(z,x)+y \subset\PP(5,1,1,5)$  &8,10
\\
& $(j\le 2)\ \& \ (n=7;9)$ &$\PP(1,1,5)$,\
$\Diff$=\scriptsize{$(\frac{1}{2},0,0,\frac{n-1}{n})$}
&\\
\hline
3& $zy^n,n=5 \& j=1;$ &
$t+f(z,x)+zy \subset \PP(5,1,1,4)$  & 6,8
\\
&$n=7\ \& \ i\le 2\ \& \ j\le 3$ &  $ \PP(1,1,4)$,\
$\Diff$=\scriptsize{$(\frac{1}{2},0,0,\frac{n-1}{n})$} &
\\
\hline
4& $y^3f'_3(z,x)+zy^6$ &
$t+f(z,x)+g(z,x,y^{1/3}) \subset \PP(5,1,1,2)$  &6
\\
(1)& &$\PP(1,1,2)$,\
$\Diff$=\scriptsize{$(\frac{1}{2},0,0,\frac{2}{3})$} &
\\
\hline
5& $zxy^4\ ||| \ z^2y^4$ &
$t^2+f(z,x)+g(z,x,y^{1/2}) \subset \PP(5,2,2,3)$  & 3
\\
&$j\le 1\ |||\ (i\le 1 \& j\le 2)$ &
$\Diff$=\scriptsize{$(0,0,0,\frac{1}{2})$} &
\\
\hline
6& $zxy^5\ ||| \ z^2y^5$ &
$t+f(z,x)+g(z,x,y^{1/5}) \subset \PP(5,1,1,3)$  &6
\\
&$j\le 2\ |||\ (i\le 1\& j\le 2)$ &$\PP(1,1,3)$,\
$\Diff$=\scriptsize{$(\frac{1}{2},0,0,\frac{4}{5})$} &
\\
\hline
\end{supertabular}
}

\end{center}

\footnotesize{
\par (1) The singularity is
$t^2+z^ix^jf_{5-i-j}(z,x)+ay^3z^kx^lf_{3-k-l}(z,x)+bzy^6$.
Assume that $f_5\in \mathcal M_1$ after each
quasihomogeneous coordinate change. The canonical singularity of this kind is
exceptional, except the cases:
$(i=3 \&  k=1)$; $(i=2\&  k\ge 1)$; $(j=3 \&  l=1)$.
}

\begin{center}


\vspace{1cm}

\begin{center}
\small{12. Singularity -- $t^2+g(z,x,y), \text{where}\ g_5\in \mathcal M_1$}
\end{center}

\tabletail{\hline}
\tablehead{\hline}
\scriptsize{
\begin{supertabular}{|c|c|c|c|}
\hline
1& $z^5+zx^3y+y^{2n}$ & $t^2 +z^5+zxy+y^{2n} \subset \PP(5n,2n,8n-5,5)$  & 7,9,15
\\
&$n=3;6;8,9$ &
$\Diff$=\scriptsize{$(0,0,\frac{2}{3},0)$} &
\\
\hline
2 & $z^5+zx^3y+y^n$ & $t +z^5 +zxy+ y^n \subset
\PP(5n,n,4n-5,5)$  &10,14,22
\\
&$n=7;9;13;19$ & $\PP(n,4n-5,5)$,\
$\Diff$=\scriptsize{$(\frac12,0,\frac23,0)$}
&30\\
\hline
3& $z^5+zx^3y+az^2xy^{2n+1}+$ & $ t^2 +g(z,x,y)
\subset\PP(15n+5,6n+2,8n+1,5)$  &4,5
\\
&$+by^{6n+2},\ n=1\ b\ne 0;2$ &$\Diff=\emptyset$
& \\
\hline
4& $z^5+zx^3y+y^4z^if_{3-i}(z,y^2)+$ &
$t^2 +g(z,x^{1/3},y) \subset \PP(5,2,7,1)$  & 3
\\
&$ax^3y^3, i\le 2 || a\ne 0$ &
$\Diff$=\scriptsize{$(0,0,\frac{2}{3},0)$} &
\\
\hline
5& $z^5+zx^3y+az^2xy^{n+1}+$ & $ t +g(z,x,y)
\subset\PP(15n+10,3n+2,4n+1,5)$  &10
\\
&$by^{3n+2},\ n=3,5$ & $\PP(3n+2,4n+1,5)$,\
$\Diff$=\scriptsize{$(\frac12,0,0,0)$}
& \\
\hline
6& $z^5+zx^3y+by^6z^if_{3-i}(z,y^3)+$ &
$t +g(z,x^{1/3},y) \subset \PP(15,3,11,1)$  & 6
\\
&$ax^3y^4$ &$\PP(3,11,1)$,\
$\Diff$=\scriptsize{$(\frac12,0,\frac{2}{3},0)$} &
\\
\hline
7& $z^5+zx^3y+xy^{2n+1}$ & $ t^2 +g(z,x,y)
\subset\PP(15n+5,6n+2,8n-1,11)$  &3,5,11
\\
&$n=2;3;4,5,6$ &$\Diff=\emptyset$
& \\
\hline
8& $z^5+zx^3y+xy^n$ & $ t +g(z,x,y)
\subset\PP(15n-5,3n-1,4n-5,11)$  &8,14,22
\\
&$n=6;8,10;12,14$ & $\PP(3n-1,4n-5,11)$,\
$\Diff$=\scriptsize{$(\frac12,0,0,0)$}
& \\
\hline
9& $z^5+zx^3y+x^2y^{2n}$ & $ t^2 +g(z,x,y)
\subset\PP(15n-5,6n-2,8n-5,7)$  &6,7
\\
&$n=2;3,4$ &$\Diff=\emptyset$
& \\
\hline
10& $z^5+zx^3y+x^2y^n$ & $ t +g(z,x,y)
\subset\PP(15n-10,3n-2,4n-5,7)$  &6,14
\\
&$n=5;7,9$ & $\PP(3n-2,4n-5,7)$,\
$\Diff$=\scriptsize{$(\frac12,0,0,0)$}
& \\
\hline
11 & $z^5+zx^3y+z^2y^{2n}$ & $t^2 +z^5 +zxy+ z^2y^{2n} \subset
\PP(5n,2n,8n-3,3)$  &9
\\
&$n=4,5$ &
$\Diff$=\scriptsize{$(0,0,\frac23,0)$}
&\\
\hline
12 & $z^5+zx^3y+z^2y^n$ & $t +z^5 +zxy+ z^2y^n \subset
\PP(5n,n,4n-3,3)$  &10,18
\\
&$n=5,7;11$ & $\PP(n,4n-3,3)$,\
$\Diff$=\scriptsize{$(\frac12,0,\frac23,0)$}
&\\
\hline
13& $z^5+x^4y+y^{2n}$ & $t^2 +z+xy+y^{2n} \subset \PP(n,2n,2n-1,1)$  & 5,12,20
\\
&$n=3,4;6;7$ & $\PP(n,2n-1,1)$,\
$\Diff$=\scriptsize{$(0,\frac45,\frac34,0)$} &
\\
\hline
14& $z^5+x^4y+ax^2y^{n+1}+$ &
$t^2 +g(z^{1/5},x,y) \subset \PP(2n+1,4n+2,n,2)$  & 5
\\
&$by^{2n+1},\ n=3,5$ &$\PP(2n+1,n,2)$,\
$\Diff$=\scriptsize{$(0,\frac45,0,0)$} &
\\
\hline
15& $z^5+x^4y+ay^nx^if_{3-i}(x,y^{n-1})$ &
$t +g(z^{1/5},x,y) \subset \PP(4n-3,4n-3,n-1,1)$  & 6,10
\\
&$n=3\ i\le 2;4$ &$\PP(4n-3,n-1,1)$,\
$\Diff$=\scriptsize{$(\frac12,\frac45,0,0)$} &
\\
\hline
16& $z^5+x^4y+ay^4z^if_{3-i}(z,y^2)$ &
$t^2 +g(z,x^{1/4},y) \subset \PP(5,2,9,1)$  & 4
\\
& &
$\Diff$=\scriptsize{$(0,0,\frac34,0)$} &
\\
\hline
17& $z^5+x^4y+xy^n$ &
$t +z+x^4y+xy^n \subset \PP(4n-1,4n-1,n-1,3)$  & 6,10,30
\\
&$n=5;6,8;11$ &$\PP(4n-1,n-1,3)$,\
$\Diff$=\scriptsize{$(\frac12,\frac45,0,0)$} &
\\
\hline
18& $z^5+x^4y+az^3xy^2+bz^2y^7$ &
$t +g(z,x,y) \subset \PP(35,7,8,3)$  & 6
\\
&$czx^2y^4+dxy^9,|b|+|c|+|d|\ne 0$ &$\PP(7,8,3)$,\
$\Diff$=\scriptsize{$(\frac12,0,0,0)$} &
\\
\hline
19& $z^5+x^4y+zy^n$ &
$t^2 +z^5+x^2y+zy^n \subset \PP(5n,2n,5n-4,8)$  & 3,12,16
\\
&$n=5;7,9;11$ &
$\Diff$=\scriptsize{$(0,0,\frac12,0)$} &
\\
\hline
20& $z^5+x^4y+az^3y^n+$ &
$t^2 +g(z,x^{1/2},y) \subset \PP(5n,2n,5n-2,4)$  &4,8
\\
&$bzy^{2n},n=3 \ b\ne 0;5$ &
$\Diff$=\scriptsize{$(0,0,\frac12,0)$} &
\\
\hline
21& $z^5+x^4y+zxy^n$ &
$t +g(z,x,y) \subset \PP(20n-5,4n-1,5n-4,11)$  & 6,18,22
\\
&$n=4,5;6,7;8$ &$\PP(4n-1,5n-4,11)$,\
$\Diff$=\scriptsize{$(\frac12,0,0,0)$} &
\\
\hline
22& $z^5+x^4y+az^2y^5+$ & $ t^2 +g(z,x,y)
\subset\PP(25,10,11,6)$  &3
\\
&$bzx^2y^3$ &$\Diff=\emptyset$
& \\
\hline
23& $z^5+x^4y+z^2y^{2n}$ &
$t^2 +z^5+xy+z^2y^{2n} \subset \PP(5n,2n,10n-3,3)$  &10,12
\\
&$n=2;3$ &
$\Diff$=\scriptsize{$(0,0,\frac34,0)$} &
\\
\hline
24& $z^5+x^4y+z^2xy^n$ &
$t +g(z,x,y) \subset \PP(20n-5,4n-1,5n-3,7)$  & 10,14
\\
&$n=3,4;5$ &$\PP(4n-1,5n-3,7)$,\
$\Diff$=\scriptsize{$(\frac12,0,0,0)$} &
\\
\hline
25& $z^4x+zx^3y+y^{2n}$ & $ t^2 +g(z,x,y)
\subset\PP(11n,4n+1,6n-4,11)$  &5,7,11
\\
&$n=3,4;5;6,7$ &$\Diff=\emptyset$
& \\
\hline
26& $z^4x+zx^3y+y^n$ &
$t +g(z,x,y) \subset \PP(11n,2n+1,3n-4,11)$  & 10,12,16
\\
&$n=7;9;11;13;15$ &$\PP(2n+1,3n-4,11)$,\
$\Diff$=\scriptsize{$(\frac12,0,0,0)$} &18,22
\\
\hline
27& $z^4x+zx^3y+az^3y^{2n-1}+$ & $ t^2 +g(z,x,y)
\subset\PP(11n-4,4n-1,6n-4,5)$  &4,5
\\
&$bx^2y^{2n}, n=2\ ab\ne 0;3$ &$\Diff=\emptyset$
& \\
\hline
28& $z^4x+zx^3y+az^3y^{n-1}+$ &
$t +g(z,x,y) \subset \PP(11n-8,2n-1,3n-4,5)$  & 6,10
\\
&$bx^2y^n, n=5\ b\ne 0;7$ &$\PP(2n-1,3n-4,5)$,\
$\Diff$=\scriptsize{$(\frac12,0,0,0)$} &
\\
\hline
29& $z^4x+zx^3y+z^2y^{2n}$ & $ t^2 +g(z,x,y)
\subset\PP(11n+1,4n+1,6n-2,7)$  &7
\\
&$n=3,4$ &$\Diff=\emptyset$
& \\
\hline
30& $z^4x+zx^3y+z^2y^n$ &
$t +g(z,x,y) \subset \PP(11n+2,2n+1,3n-2,7)$  & 8,10,14
\\
&$n=5;7;9$ &$\PP(2n+1,3n-2,7)$,\
$\Diff$=\scriptsize{$(\frac12,0,0,0)$} &
\\
\hline
31& $z^4x+x^4y+y^n$ &
$t^2 +zx+x^4y+y^n \subset \PP(2n,3n+1,n-1,4)$  &5,12,16
\\
&$n=6,8;10;12$ &
$\Diff$=\scriptsize{$(0,\frac34,0,0)$} &
\\
\hline
32& $z^4x+x^4y+ax^2y^{n+1}+$ &
$t^2 +g(z^{1/4},x,y) \subset \PP(2n+1,3n+2,n,2)$  &4,8
\\
&$by^{2n+1},\ n=3 \ b\ne 0;5$ &
$\Diff$=\scriptsize{$(0,\frac34,0,0)$} &
\\
\hline
33& $z^4x+x^4y+by^3x^if_{3-i}(x,y^2)$ &
$t^2 +g(z^{1/2},x,y) \subset \PP(9,7,4,2)$  &4
\\
&$+az^4y^2,\ i\le 2 ||  a\ne 0$ &
$\Diff$=\scriptsize{$(0,\frac12,0,0)$} &
\\
\hline
34& $z^4x+x^4y+zy^n$ &
$t +g(z,x,y) \subset \PP(16n+1,3n+1,4n-3,13)$  & 6,10,22
\\
&$n=5;6,7,8;9;10$ &$\PP(3n+1,4n-3,13)$,\
$\Diff$=\scriptsize{$(\frac12,0,0,0)$} & 26
\\
\hline
35& $z^4x+x^4y+az^2y^{2n-1}+$ &
$t +g(z,x,y) \subset \PP(16n-7,3n-1,4n-3,5)$  & 6,10
\\
&$bzx^2y^n, n=3\ a\ne 0;4$ &$\PP(3n-1,4n-3,5)$,\
$\Diff$=\scriptsize{$(\frac12,0,0,0)$} &
\\
\hline
36& $z^4x+x^4y+z^2y^n$ & $ t^2 +g(z,x,y)
\subset\PP(8n+1,3n+1,4n-2,10)$  &5
\\
&$n=4,6$ &$\Diff=\emptyset$
& \\
\hline
37& $z^4x+x^4y+z^3y^n$ &
$t +g(z,x,y) \subset \PP(16n+3,3n+1,4n-1,7)$  & 10,14
\\
&$n=4;5$ &$\PP(3n+1,4n-1,7)$,\
$\Diff$=\scriptsize{$(\frac12,0,0,0)$} &
\\
\hline
38& $z^4y+z^2x^3+y^{10}$ & $t^2 +z^2y+zx+y^{10} \subset \PP(10,9,11,2)$
& 12
\\
&&
$\Diff$=\scriptsize{$(0,\frac12,\frac23,0)$} &
\\
\hline
39& $z^4y+z^2x^3+az^2y^4+$ & $t^2 +g(z,x^{1/3},y) \subset \PP(7,3,8,2)$
& 3
\\
&$bx^3y^3+cy^7, |ab|+|c|\ne 0$&
$\Diff$=\scriptsize{$(0,0,\frac23,0)$} &
\\
\hline
40& $z^4y+z^2x^3+az^2xy^3+bx^4y^2+$ & $t^2 +g(z^{1/2},x,y) \subset
\PP(8,7,3,2)$ & 4
\\
&$cx^2y^5+dy^8, |ab|+|c|+|d|\ne 0$&
$\Diff$=\scriptsize{$(0,\frac12,0,0)$} &
\\
\hline
41& $z^4y+z^2x^3+cy^5z^if_{2-i}(z,y^2)$ &
$t +g(z,x^{1/3},y) \subset \PP(9,2,5,1)$  & 6
\\
&$+azx^3y^2+bx^3y^4$ &$\PP(2,5,1)$,\
$\Diff$=\scriptsize{$(\frac12,0,\frac{2}{3},0)$} &\\
&$|a|+|b|+|c|\ne 0 \& i\le 1$&&\\
\hline
42& $z^4y+z^2x^3+xy^n$ & $ t^2 +g(z,x,y)
\subset\PP(6n+1,3n-2,2n+2,10)$  &3,5
\\
&$n=5;7$ &$\Diff=\emptyset$
& \\
\hline
43& $z^4y+z^2x^3+azx^2y^n+$ &
$t +g(z,x,y) \subset \PP(12n+1,3n-1,2n+1,5)$  & 8,10
\\
&$bxy^{2n}, n=3;4$ &$\PP(3n-1,2n+1,5)$,\
$\Diff$=\scriptsize{$(\frac12,0,0,0)$} &
\\
\hline
44& $z^4y+z^2x^3+x^2y^n$ &
$t^2 +g(z^{1/2},x,y) \subset \PP(3n+1,3n-1,n+1,4)$  &6,8
\\
&$n=4;6$ &
$\Diff$=\scriptsize{$(0,\frac12,0,0)$} &
\\
\hline
45& $z^4y+z^2x^3+zy^8$ &
$t +z^4y+z^2x+zy^8 \subset \PP(31,7,17,3)$  & 18
\\
& &$\PP(7,17,3)$,\
$\Diff$=\scriptsize{$(\frac12,0,\frac{2}{3},0)$} &
\\
\hline
46& $z^4y+z^2x^3+zxy^n$ &
$t +g(z,x,y) \subset \PP(12n-1,3n-2,2n+1,7)$  & 12,14
\\
&$n=5;6$ &$\PP(3n-2,2n+1,7)$,\
$\Diff$=\scriptsize{$(\frac12,0,0,0)$} &
\\
\hline
47& $z^4y+zx^4+y^n$ &
$t^2 +z^4y+zx+y^n \subset \PP(2n,n-1,3n+1,4)$  &5,12,16
\\
&$n=6,8;12;16$ &
$\Diff$=\scriptsize{$(0,0,\frac34,0)$} &
\\
\hline
48& $z^4y+zx^4+az^2y^{n+1}+$ &
$t^2 +g(z,x^{1/4},y) \subset \PP(2n+1,n,3n+2,2)$  & 4,8
\\
&$by^{2n+1}, n=3\ b\ne 0;5$ &
$\Diff$=\scriptsize{$(0,0,\frac34,0)$} &
\\
\hline
49& $z^4y+zx^4+by^5z^if_{2-i}(z,y^2)$ &
$t^2 +g(z,x^{1/2},y) \subset \PP(9,4,7,2)$  & 4
\\
&$+ax^4y^2,\ i\le 2$ &
$\Diff$=\scriptsize{$(0,0,\frac12,0)$} &
\\
\hline
50& $z^4y+zx^4+xy^n$ &
$t +g(z,x,y) \subset \PP(16n+1,4n-3,3n+1,13)$  & 6,10,22
\\
&$n=5;6,7,8;9;10$ &$\PP(4n-3,3n+1,13)$,\
$\Diff$=\scriptsize{$(\frac12,0,0,0)$} &
26\\
\hline
51& $z^4y+zx^4+x^2y^n$ & $ t^2 +g(z,x,y)
\subset\PP(8n+1,4n-2,3n+1,10)$  &5
\\
&$n=4,6$ &$\Diff=\emptyset$
& \\
\hline
52& $z^4y+zx^4+x^2y^{2n+1}$ &
$t +g(z,x,y) \subset \PP(16n+9,4n+1,3n+2,5)$  & 6,10
\\
&$n=2;3$ &$\PP(4n+1,3n+2,5)$,\
$\Diff$=\scriptsize{$(\frac12,0,0,0)$} &
\\
\hline
53& $z^4y+zx^4+x^3y^n$ &
$t +g(z,x,y) \subset \PP(16n+3,4n-1,3n+1,7)$  & 10,14
\\
&$n=4;5$ &$\PP(4n-1,3n+1,7)$,\
$\Diff$=\scriptsize{$(\frac12,0,0,0)$} &
\\
\hline
54& $z^3x^2+zx^3y+y^{2n}$ & $ t^2 +g(z,x,y)
\subset\PP(7n,2n+2,4n-3,7)$  &5,7
\\
&$n=4;5$ &$\Diff=\emptyset$
& \\
\hline
55& $z^3x^2+zx^3y+y^n$ &
$t +g(z,x,y) \subset \PP(7n,n+2,2n-3,7)$  & 8,10,14
\\
&$n=7;9;11$ &$\PP(n+2,2n-3,7)$,\
$\Diff$=\scriptsize{$(\frac12,0,0,0)$} &
\\
\hline
56& $z^3x^2+zx^3y+az^2y^{2n}+$ & $ t^2 +g(z,x,y)
\subset\PP(7n+2,2n+2,4n-1,5)$  &3,5
\\
&$bxy^{2n+1}, n=2\ ab\ne 0;3$ &$\Diff=\emptyset$
& \\
\hline
57& $z^3x^2+zx^3y+az^2y^{n-1}+$ & $ t +g(z,x,y)
\subset\PP(7n-3,n+1,2n-3,5)$  &6,10
\\
&$bxy^n, n=6\ a\ne 0;8$ &$\PP(n+1,2n-3,5)$,\
$\Diff$=\scriptsize{$(\frac12,0,0,0)$}
& \\
\hline
58& $z^3x^2+zx^3y+az^4y^2+bz^2xy^3+$ & $ t^2 +g(z,x,y)
\subset\PP(11,4,5,3)$  &3
\\
&$czy^6+dx^2y^4,\ |c|+|d|\ne 0\ \&$ &$\Diff=\emptyset$
& \\
&$|ad|+|b|+|c|\ne 0$&&\\
\hline
59& $z^3x^2+zx^3y+az^4y^3+bz^2xy^4+$ & $ t+g(z,x,y)
\subset\PP(29,5,7,3)$  &6
\\
&$czy^8+dx^2y^5,|ad|+|b|+|c|\ne 0$ &$\PP(5,7,3)$,\
$\Diff$=\scriptsize{$(\frac12,0,0,0)$}
& \\
\hline
60& $z^3x^2+x^3y^2+y^{10}$ &
$t^2 +zx^2+x^3y^2+y^{10} \subset \PP(15,14,8,3)$  &7
\\
& &
$\Diff$=\scriptsize{$(0,\frac23,0,0)$} &
\\
\hline
61& $z^3x^2+x^3y^2+y^n$ &
$t +zx^2+x^3y^2+y^n \subset \PP(3n,n+4,n-2,3)$  &9,12,18
\\
&$n=7;9;13$ &$\PP(n+4,n-2,3)$,\
$\Diff$=\scriptsize{$(\frac12,\frac23,0,0)$} &
\\
\hline
62& $z^3x^2+x^3y^2+az^6y+bz^3y^6+$ & $ t+g(z^{1/3},x,y)
\subset\PP(11,5,3,1)$  &6
\\
&$cz^3xy^3+dx^2y^5+exy^8+ly^{11}$ &$\PP(5,3,1)$,\
$\Diff$=\scriptsize{$(\frac12,\frac23,0,0)$}
& \\
& $|ad|+|b|+|c|+|e|+|l|\ne 0$&&\\
\hline
63& $z^3x^2+x^3y^2+y^{12}$ &
$t^2 +zx^2+x^3y+y^6 \subset \PP(9,8,5,3)$  &9
\\
& &
$\Diff$=\scriptsize{$(0,\frac23,0,\frac12)$} &
\\
\hline
64& $z^3x^2+x^3y^2+az^3y^{n-2}+$ & $ t +g(z^{1/3},x,y)
\subset\PP(3n-2,n+2,n-2,2)$  &6,12
\\
&$bxy^n, n=5\ a\ne 0;9$ &$\PP(n+2,n-2,2)$,\
$\Diff$=\scriptsize{$(\frac12,\frac23,0,0)$}
& \\
\hline
65& $z^3x^2+x^3y^2+az^5y^2+bz^4xy+$ & $ t+g(z,x,y)
\subset\PP(19,3,5,2)$  &4
\\
(1)&$cz^3y^5+dz^2xy^4+ezx^2y^3+$ &$\PP(3,5,2)$,\
$\Diff$=\scriptsize{$(\frac12,0,0,0)$}
& \\
& $lzy^8+nxy^7$&&\\
\hline
66& $z^3x^2+x^3y^2+zy^n$ &
$t +g(z,x,y) \subset \PP(9n+4,n+4,3n-4,8)$  &7,10,14
\\
&$n=5;7;9;11$ &$\PP(n+4,3n-4,8)$,\
$\Diff$=\scriptsize{$(\frac12,0,0,0)$} &16
\\
\hline
67& $z^3x^2+x^3y^2+az^2xy^n+$ & $ t +g(z,x,y)
\subset\PP(9n+2,n+2,3n-2,4)$  &6,8
\\
&$bzy^{2n}, n=3;5$ &$\PP(n+2,3n-2,4)$,\
$\Diff$=\scriptsize{$(\frac12,0,0,0)$}
& \\
\hline
68& $z^3x^2+x^3y^2+az^4y^2+bzxy^4$ & $ t^2 +g(z,x,y)
\subset\PP(17,6,8,5)$  &3
\\
&$ab\ne 0$ &$\Diff=\emptyset$
& \\
\hline
69& $z^3x^2+x^3y^2+az^4y^{n-2}+$ & $ t +g(z,x,y)
\subset\PP(9n-2,n+2,3n-4,5)$  &8,10
\\
&$bzxy^n,\ n=5;7$ &$\PP(n+2,3n-4,5)$,\
$\Diff$=\scriptsize{$(\frac12,0,0,0)$}
& \\
\hline
70& $z^3x^2+x^3y^2+az^4y^4+$ & $ t^2 +g(z,x,y^{1/2})
\subset\PP(13,4,7,5)$  &5
\\
&$bzxy^6$ &
$\Diff$=\scriptsize{$(0,0,0,\frac12)$}
& \\
\hline
71& $z^3x^2+x^3y^2+z^2y^{4n}$ & $ t^2 +g(z,x,y^{1/2})
\subset\PP(9n+2,2n+2,6n-1,7)$  &6,7
\\
&$n=1;2$ &
$\Diff$=\scriptsize{$(0,0,0,\frac12)$}
& \\
\hline
72& $z^3x^2+x^3y^2+z^2y^n$ & $ t +g(z,x,y)
\subset\PP(9n+8,n+4,3n-2,7)$  &8,12,14
\\
&$n=5;7;9$ &$\PP(n+4,3n-2,7)$,\
$\Diff$=\scriptsize{$(\frac12,0,0,0)$}
& \\
\hline
73& $z^3x^2+x^3y^2+z^2y^6$ & $ t^2 +g(z,x,y)
\subset\PP(31,10,16,7)$  &5
\\
& &$\Diff=\emptyset$
& \\
\hline

\end{supertabular}
}
\end{center}

\footnotesize{
\par (1) The exceptionality condition is $|c|+|d|+|l|+|n|+|e(|a|+|b|)|\ne 0$.
}

\begin{center}


\vspace{1cm}

\begin{center}
\small{13. Singularity -- $t^3+g(t,z,x,y)$}
\end{center}

\tabletail{\hline}
\tablehead{\hline}
\scriptsize{
\begin{supertabular}{|c|c|c|c|}
\hline
1& $z^2x^2+x^iy^jf_{5-i-j}(x,y)$ & $t +zx^2+x^iy^jf_{5-i-j}(x,y)
\subset \PP(5,3,1,1)$  & 6
\\
&$2i+j\le 3$ & $\PP(3,1,1)$,\
$\Diff$=\scriptsize{$(\frac{2}{3},\frac12,0,0)$} &
\\
\hline
2& $z^2x^2+x^5+zy^4$ &
$t +z^2x^2+x^5+zy\subset \PP(10,3,2,7)$  & 21
\\
& & $\PP(3,2,7)$,\
$\Diff$=\scriptsize{$(\frac{2}{3},0,0,\frac34)$} &
\\
\hline
3& $z^2x^2+x^4y+zy^4$ &
$t +z^2x+x^2y+zy^4\subset \PP(17,5,7,3)$  & 18
\\
& & $\PP(5,7,3)$,\
$\Diff$=\scriptsize{$(\frac{2}{3},0,\frac12,0)$} &
\\
\hline
4& $z^2x^2+x^3y^2+zy^4$ &
$t +z^2x^2+x^3y+zy^2\subset \PP(14,4,3,5)$  & 15
\\
& & $\PP(4,3,5)$,\
$\Diff$=\scriptsize{$(\frac{2}{3},0,0,\frac12)$} &
\\
\end{supertabular}
}

\end{center}

\footnotesize{
\par (0) $n=1\&k\le 1 \& (a\ne 0 || (i=0\& j\le 1))$,
$n=2\&k\le 1 \& (a\ne 0 || (i\le 1\& j\le 4))$.

\par (1) The singularity is
$t^3+z^2x+\underline{a}tx^i(x+y^2)^jf_{3-i-j}(x,y^2)+
\underline{b}x^ky(x+y^2)^lf_{4-k-l}(x,y^2)$.
The common exceptionality condition for two cases is
$(i=0 || k\le 1)$. The first case is $(b\ne 0 \ ||\ j\le 2)$ and the coefficient
$\underline{a}$ is absent (see the comments to the tables).
If the coefficient $\underline{b}$ is absent then we have
the second case.
}

\begin{center}


\vspace{1cm}

\begin{center}
\small{14. Singularity -- $tz^if_{2-i}(t,z)+g(z,x,y)(i\le 1)$}
\end{center}

\tabletail{\hline}
\tablehead{\hline}
\scriptsize{
\begin{supertabular}{|c|c|c|c|}
\hline 1 & $x^4+y^n$ & $tz^if_{2-i}(t,z)+x+y \subset \PP(1,1,3,3)$
&5,8,12
\\
& $n=5;7;11$ &$\PP(1,1,3)$,\
$\Diff$=\scriptsize{$(0,0,\frac{3}{4},\frac{n-1}{n})$}
&\\
\hline
2 & $f_3(t,z)+y^2f_2(t,z)+$ &
$tz^if_{2-i}(t,z)+g(z,x^{1/2},y) \subset\PP(2,2,3,1)$  &2
\\
(1)& $(t+z)(x^2y+y^4)+f_2(x^2,y^3) $ &
$\Diff$=\scriptsize{$(0,0,\frac{1}{2},0)$}
&\\
\hline
3& $x^jf_{4-j}(x,y^2)$ & $tz^if_{2-i}(t,z)+g(z,x,y^{1/2}) \subset
\PP(4,4,3,3)$  & 3
\\
&$j\le 2$ &$\Diff$=\scriptsize{$(0,0,0,\frac{1}{2})$} &
\\
\hline
4 & $t^3+ktz^2+lz^3+(at+bz)y^6+$ & $tz^if_{2-i}(t,z)+g(z,x^{1/4},y^{1/3}) \subset
\PP(1,1,3,1)$  &4
\\
&$ctzy^3+dy^9+x^4$&$\PP^2$,\
$\Diff$=\scriptsize{$(0,0,\frac{3}{4},\frac{2}{3})$}
&\\
\hline
5& $x^{2j}f_{2-j}(x^2,y^5)$ & $tz^if_{2-i}(t,z)+x^jf_{2-j}(x,y)
\subset \PP(2,2,3,3)$  & 6
\\
&$j\le 1$ &$\Diff$=\scriptsize{$(0,0,\frac{1}{2},\frac{4}{5})$} &
\\
\hline
6 & $x^4+xy^n$ & $tz^if_{2-i}(t,z)+x^4+xy \subset \PP(4,4,3,9)$
&5,9
\\
& $n=4,5;7,8$ &
$\Diff$=\scriptsize{$(0,0,0,\frac{n-1}{n})$}
&\\
\hline
7 & $x^4+zy^n$ & $tz^if_{2-i}(t,z)+x+zy \subset \PP(1,1,3,2)$
&4,5,8
\\
& $n=3\ i=0;5;7$ &$\PP(1,1,2)$,\
$\Diff$=\scriptsize{$(0,0,\frac{3}{4},\frac{n-1}{n})$}
&\\
\hline
8 & $x^4+zxy^n$ & $tz^if_{2-i}(t,z)+x^4+zxy \subset \PP(4,4,3,5)$
&2,5
\\
& $n=2\ i=0; n=3,4$ &
$\Diff$=\scriptsize{$(0,0,0,\frac{n-1}{n})$}
&\\
\hline
9& $x^jf_{5-j}(x,y)$ & $tz^if_{2-i}(t,z)+x^jf_{5-j}(x,y) \subset
\PP(5,5,3,3)$  & 3
\\
&$j=1$ &$\Diff=\emptyset$ &
\\
\hline
10 & $t^3+ktz^2+lz^3+(at+bz)y^4+$ & $tz^if_{2-i}(t,z)+g(z,x^{1/5},y^{1/2})
\subset \PP(1,1,3,1)$  &5
\\
&$ctzy^2+dy^6+x^5$ &$\PP^2$,\
$\Diff$=\scriptsize{$(0,0,\frac{4}{5},\frac{1}{2})$}
&\\
\hline
11 & $x^5+y^7$ & $tz^if_{2-i}(t,z)+x+y \subset \PP(1,1,3,3)$  &15
\\
&  &$\PP(1,1,3)$,\
$\Diff$=\scriptsize{$(0,0,\frac{4}{5},\frac{6}{7})$}
&\\
\hline
12 & $x^5+xy^5$ & $tz^if_{2-i}(t,z)+x^5+xy \subset \PP(5,5,3,12)$
&6
\\
&  &$\Diff$=\scriptsize{$(0,0,0,\frac{4}{5})$}
&\\
\hline
13 & $x^5+x^2y^4$ & $tz^if_{2-i}(t,z)+x^5+x^2y \subset
\PP(5,5,3,9)$  &9
\\
&  &$\Diff$=\scriptsize{$(0,0,0,\frac{3}{4})$}
&\\
\hline
14 & $x^5+zxy^3$ & $tz^if_{2-i}(t,z)+x^5+zxy \subset \PP(5,5,3,7)$
&7
\\
&  &$\Diff$=\scriptsize{$(0,0,0,\frac{2}{3})$}
&\\
\hline
15 & $t^3+ktz^2+lz^3+(at+bz)y^4+$ & $tz^if_{2-i}(t,z)+g(z,x^{1/4},y)
\subset \PP(2,2,5,1)$
&4
\\
&$ctzy^2+dy^6+x^4y$  & $\Diff$=\scriptsize{$(0,0,\frac{3}{4},0)$}
&\\
\hline
16 & $x^4y+ax^2y^4+by^7$ &
$tz^if_{2-i}(t,z)+g(z,x^{1/2},y) \subset\PP(7,7,9,3)$  &6
\\
& $b\ne 0$ &
$\Diff$=\scriptsize{$(0,0,\frac{1}{2},0)$}
&\\
\hline
17 & $x^4y+y^8$ & $tz^if_{2-i}(t,z)+g(z,x^{1/4},y)
\subset \PP(8,8,21,3)$
&12
\\
&  & $\Diff$=\scriptsize{$(0,0,\frac{3}{4},0)$}
&\\
\hline
18& $x^4y+xy^n$ & $tz^if_{2-i}(t,z)+g(z,x,y) \subset
\PP(4n-1,4n-1,3n-3,9)$  & 6,9
\\
&$n=5;6$ &$\Diff=\emptyset$ &
\\
\hline
19 & $x^4y+zy^5$ & $tz^if_{2-i}(t,z)+g(z,x^{1/4},y)
\subset \PP(5,5,13,2)$
&8
\\
&  & $\Diff$=\scriptsize{$(0,0,\frac{3}{4},0)$}
&\\
\hline
20& $x^4y+zxy^3$ & $tz^if_{2-i}(t,z)+x^4y+zxy^3 \subset
\PP(11,11,7,5)$  & 5
\\
& &$\Diff=\emptyset$ &
\\
\hline
21 & $zx^3+y^n$ & $tz^if_{2-i}(t,z)+zx+y \subset \PP(1,1,2,3)$
&5,7,9
\\
& $n=5\ i=0;7;8$ &$\PP(1,1,2)$,\
$\Diff$=\scriptsize{$(0,0,\frac{2}{3},\frac{n-1}{n})$}
&\\
\hline
22 & $t^3+ktz^2+lz^3+tx^3+aty^4+sy^6$ &
$tz^if_{2-i}(t,z)+g(z,x^{1/3},y^{1/2}) \subset \PP(1,1,2,1)$  & 3
\\
(2)&$+z(bx^3+cy^4)+dx^3y^2+etzy^2$ &
$\Diff$=\scriptsize{$(0,0,\frac{2}{3},\frac{1}{2})$}
&\\
\hline
23 & $zx^3+xy^n$ & $tz^if_{2-i}(t,z)+zx^3+xy \subset \PP(3,3,2,7)$
&4,5,7
\\
& $n=4\ i=0;5;6$ &
$\Diff$=\scriptsize{$(0,0,0,\frac{n-1}{n})$}
&\\
\hline
24 & $zx^3+x^2y^n$ & $tz^if_{2-i}(t,z)+zx^3+x^2y \subset
\PP(3,3,2,5)$  &3,5
\\
& $n=3\ i=0;4$ &
$\Diff$=\scriptsize{$(0,0,0,\frac{n-1}{n})$}
&\\
\hline

\end{supertabular}
}

\end{center}

\footnotesize{
\par (1) See example 3.26.
\par (2) If $d=e=0$ then we require that the polynomial without
monomial $y^6$ is irreducible.
}

\begin{center}



\begin{center}
\small{15. Singularity -- $t^2z+g(t,z,x,y)$}
\end{center}
\tabletail{\hline}
\tablehead{\hline}
\scriptsize{
\begin{supertabular}{|c|c|c|c|}
\hline 1& $z^ix^jf_{4-i-j}(x,y)+y^n$ & $tz+z^ix^jf_{4-i-j}(x,y)+y
\subset \PP(3,1,1,4)$  & 5,8
\\
&$n=5 \ i=0 \& j\le 2; 7$ & $\PP(3,1,1)$,\
$\Diff$=\scriptsize{$(\frac12,0,0,\frac{n-1}{n})$} &
\\
\hline
2& $f_4(z,x)+y^3f_2(z,x)+$ & $tz+g(z,x,y^{1/3})
\subset \PP(3,1,1,2)$  & 4
\\
(1)&$+y^6$ &
$\Diff$=\scriptsize{$(\frac12,0,0,\frac23)$} &
\\
\hline
3& $z^ix^jf_{4-i-j}(x,y)+ty^n$ & $t^2z+z^ix^jf_{4-i-j}(x,y)+ty
\subset \PP(3,2,2,5)$  & 3,5
\\
&$n=3 \ i=0 \& j\le 1; 4$ &
$\Diff$=\scriptsize{$(0,0,0,\frac{n-1}{n})$} &
\\
\hline
4& $z^4+x^3y+y^n$ & $tz+z^4+xy+y^n
\subset \PP(3n,n,4n-4,4)$  & 7,15,16
\\
&$n=5;9;11;15$ &
$\Diff$=\scriptsize{$(\frac12,0,\frac23,0)$} &
24\\
\end{supertabular}
}

\end{center}

\footnotesize{
\par (1) The singularity is $t^2z+z^ix^jf_{4-i-j}(x,y^2)+
az^kx^ly^3f_{2-k-l}(x,y)+by^6$. The exceptionality condition is\\
$(i=0 \ ||\ k=0)\ \& \ \Big( \big(b\ne 0\ \&\ (j,l\le 2)\big)\ ||\
\big(b=0\ \&\ (j,l\le 1)\big)      \Big)$
\par (2) There are three cases. A). The singularity is
$t^2z+z^2x+azx^iy^{j+1}(x+y)^kf_{2-i-j-k}(x,y)+
x^{i_1}y(x+y)^{k_1}f_{4-i_1-k_1}(x,y)$. The exceptionality condition is
$(i=0 || i_1=0)\& (k=0 || k_1\le 1)$.\\
B). The singularity is
$t^2z+z^2x+zx^3+azx^iy^{j+1}(x+y)^kf_{2-i-j-k}(x,y)+
bx^{i_1}y^{j_1+1}(x+y)^{k_1}f_{4-i_1-j_1-k_1}(x,y)$. The exceptionality condition is
$(i=0\ ||\ i_1=0)\ \& \ (j_1\le 2\ \&\ k_1\le 3\ \& \ i_1\le 3)$.\\
C). The singularity is
$t^2z+z^2x+zx^3+ax^5+bzx^iy^{j+1}(x+y)^kf_{2-i-j-k}(x,y)+
cx^{i_1}y^{j_1+1}(x+y)^{k_1}f_{4-i_1-j_1-k_1}(x,y)$, where $a\ne 0$.
The exceptionality condition is
$(i=0\ ||\ i_1=0)\ \& \ i_1\le 3$.\\
\par (3). The following notation $\underline{z^2x+zx^3+x^5}$ means one of the
following
polynomials: $z^2x+zx^3$; $z^2x+zx^3+ax^5$, where $a\ne \frac14$;
$x(z+x^2)^2$. The symbol (*) means that the third case is
impossible.
}

\begin{center}


\vspace{1cm}

\begin{center}
\small{16. Singularity -- $t^2x+g(t,z,x,y)$}
\end{center}

\tabletail{\hline}
\tablehead{\hline}
\scriptsize{
\begin{supertabular}{|c|c|c|c|}
\hline 1& $z^4+x^3y+y^n$ & $t^2x+g(z^{1/2},x,y) \subset
\PP(2n+1,3n,2n-2,6)$  & 3,12
\\
&$n=5;9$ &
$\Diff$=\scriptsize{$(0,\frac12,0,0)$} &
\\
\hline
2& $z^4+x^3y+az^2y^3+by^6$ & $ t^2x +g(z,x,y) \subset\PP(13,9,10,6)$  &
3
\\
& &$\Diff=\emptyset$
& \\
\hline
3& $z^4+x^iyf_{3-i}(x,y^2)+atz^2y$ & $t^2x+g(z^{1/2},x,y) \subset
\PP(5,7,4,2)$  & 4
\\
& &
$\Diff$=\scriptsize{$(0,\frac12,0,0)$} &
\\
\hline
4& $z^if_{4-i}(z,y^2)+x^3y$ & $tx+g(z,x,y) \subset
\PP(17,6,7,3)$  & 6
\\
&$i\le 2$ &
$\Diff$=\scriptsize{$(\frac12,0,0,0)$} &
\\
\hline
5& $z^4+x^3y+zy^n$ & $tx+g(z,x,y) \subset \PP(8n+3,3n,4n-3,9)$  &
5,6,18
\\
&$n=4;5;7$ &
$\Diff$=\scriptsize{$(\frac12,0,0,0)$} &
\\
\hline
6& $z^4+x^3y+zxy^4$ & $tx+g(z,x,y) \subset
\PP(33,11,13,5)$  & 10
\\
& &
$\Diff$=\scriptsize{$(\frac12,0,0,0)$} &
\\
\hline
7 & $z^4+x^3y+ty^n$ & $t^2x+z+x^3y+ty^n \subset
\PP(2n+1,6n+1,2n-1,4)$  & 4,5,9
\\
&$n=3;4;5;6$ & $\PP(2n+1,2n-1,4)$,\
$\Diff$=\scriptsize{$(0,\frac34,0,0)$} &
16\\
\hline
8& $z^4+x^3y+tzy^n$ & $ t^2x +g(z,x,y) \subset\PP(8n+3,6n+1,8n-2,10)$  &
5
\\
&$n=2,3$ &$\Diff=\emptyset$
& \\
\hline
9 & $z^4+x^iy^jf_{5-i-j}(x,y)+aty^3$ & $t^2x+g(z^{1/4},x,y) \subset
\PP(2,5,1,1)$  & 4
\\
&$(a\ne 0 || i=0)\& j\le 2$ & $\PP(2,1,1)$,\
$\Diff$=\scriptsize{$(0,\frac34,0,0)$} &
\\
\hline
10 & $z^4+az^2y^3+by^6+x^5$ & $t^2x+z^2+azy+by^2+x^5 \subset
\PP(4,5,2,5)$  & 10
\\
& &
$\Diff$=\scriptsize{$(0,\frac12,0,\frac23)$} &
\\
\hline
11 & $z^4+x^5+zy^4$ & $t^2x+z^4+x^5+zy \subset
\PP(8,5,4,15)$  & 5
\\
& &
$\Diff$=\scriptsize{$(0,0,0,\frac34)$} &
\\
\hline
12 & $z^4+\underline{a}x^5+tzy^2+\underline{c}zx^2y^2$ &
$t^2x+g(z,x,y^{1/2}) \subset \PP(8,5,4,7)$  & 7
\\
& &
$\Diff$=\scriptsize{$(0,0,0,\frac12)$} &
\\
\hline
13 & $z^4+az^2y^3+by^6+\underline{d}x^4y+\underline{c}zx^2y^2$ &
$tx+g(z,x,y) \subset \PP(19,6,5,4)$  & 8
\\
&$|a|+|b|\ne 0$ &
$\Diff$=\scriptsize{$(\frac12,0,0,0)$} &
\\
\hline
14 & $z^4+x^4y+ty^4$ & $t^2x+z+x^4y+ty^4 \subset
\PP(13,33,7,5)$  & 20
\\
& & $\PP(13,7,5)$,\
$\Diff$=\scriptsize{$(0,\frac34,0,0)$} &
\\
\hline
15 & $z^4+x^3y^2+az^2y^3+by^6$ & $t^2x+g(z^{1/2},x,y) \subset
\PP(7,9,4,3)$  & 6
\\
& &
$\Diff$=\scriptsize{$(0,\frac12,0,0)$} &
\\
\hline
16 & $z^4+x^3y^2+y^7$ & $t^2x+z+x^3y^2+y^7 \subset
\PP(8,21,5,3)$  & 12
\\
& & $\PP(8,5,3)$,\
$\Diff$=\scriptsize{$(0,\frac34,0,0)$} &
\\
\hline
17 & $z^4+x^3y^2+aty^4+bxy^5$ & $t^2x+g(z^{1/2},x,y) \subset
\PP(10,13,6,4)$  & 8
\\
& &
$\Diff$=\scriptsize{$(0,\frac12,0,0)$} &
\\
\hline
18& $z^4+x^3y^2+zy^n$ & $ t^2x +g(z,x,y) \subset
\PP(4n+3,3n,4n-6,9)$  &5,9
\\
&$n=4;5$ &$\Diff=\emptyset$
& \\
\hline
19& $z^4+x^3y^2+tzy^2+bzxy^3$ & $ t^2x +g(z,x,y) \subset
\PP(11,7,6,5)$  &5
\\
& &$\Diff=\emptyset$
& \\
\hline
20& $z^4+zx^2y+y^n$ & $tx+g(z,x,y) \subset
\PP(5n+4,2n,3n-4,8)$  &6,10,12
\\
&$n=5;7;9;11$ &
$\Diff$=\scriptsize{$(\frac12,0,0,0)$} &
16\\
\hline
21& $z^4+zx^2y+az^2y^n+by^{2n}$ & $tx+g(z,x,y) \subset
\PP(5n+2,2n,3n-2,4)$  &4,8
\\
&$n=3\ b\ne 0;5$ &
$\Diff$=\scriptsize{$(\frac12,0,0,0)$} &
\\
\hline
22& $z^if_{4-i}(z,y^2)+zx^2y+ax^2y^3$ & $tx+g(z,x,y) \subset
\PP(11,4,5,2)$  & 4
\\
&$(a=0 \& i\le 2) ||(a\ne 0 \& i\le 3)$ &
$\Diff$=\scriptsize{$(\frac12,0,0,0)$} &
\\
\hline
23& $z^4+zx^2y+ty^n$ & $ t^2x +g(z,x,y) \subset
\PP(5n+4,4n+1,6n-4,11)$  &5,7,11
\\
&$n=3,4;5;6,7$ &$\Diff=\emptyset$
& \\
\hline
24& $z^4+zx^2y+tzy^n$ & $ t^2x +g(z,x,y) \subset
\PP(5n+3,4n+1,6n-2,7)$  &7
\\
&$n=3,4$ &$\Diff=\emptyset$
& \\
\hline
25& $z^4+az^2y^3+by^6+zx^4$ & $tx+ z^4+az^2y+by^2+zx^4\subset
\PP(13,4,3,8)$  & 16
\\
& &
$\Diff$=\scriptsize{$(\frac12,0,0,\frac23)$} &
\\
\hline
26& $z^4+zx^4+tzy^2$ & $t^2x+ z^4+zx^4+tzy\subset
\PP(13,8,6,11)$  & 11
\\
& &
$\Diff$=\scriptsize{$(0,0,0,\frac12)$} &
\\
\hline
27& $z^4+az^2y^3+by^6+zx^3y$ & $tx+g(z,x,y)\subset
\PP(29,9,7,6)$  & 12
\\
& &
$\Diff$=\scriptsize{$(\frac12,0,0,0)$} &
\\
\hline
28& $z^4+zx^2y^2+y^7$ & $tx+g(z,x,y)\subset
\PP(43,14,13,8)$  & 16
\\
& &
$\Diff$=\scriptsize{$(\frac12,0,0,0)$} &
\\
\hline
29& $z^4+zx^2y^2+ty^4$ & $t^2x+ z^4+zx^2y+ty^2\subset
\PP(14,9,8,11)$  & 11
\\
& &
$\Diff$=\scriptsize{$(0,0,0,\frac12)$} &
\\
\hline
30& $z^3y+x^5+y^n$ & $t^2x+g(z^{2/3},x,y) \subset
\PP(2n,5n-5,n,5)$  & 9,12
\\
&$n=6;8$ &
$\Diff$=\scriptsize{$(0,\frac23,0,0)$} &
\\
\hline
31& $z^3y+x^5+ty^4+bx^2y^4$ & $t^2x+g(z^{2/3},x,y) \subset
\PP(8,17,4,3)$  & 6
\\
& &
$\Diff$=\scriptsize{$(0,\frac23,0,0)$} &
\\
\hline
32& $z^3y+z^2(x^2+y^3)+zx^2y^2+$ & $tx+g(z,x,y) \subset
\PP(11,4,3,2)$  & 4 \\
(1)&$zy^5+x^4y+x^2y^4+y^7$ &
$\Diff$=\scriptsize{$(\frac12,0,0,0)$} &
\\
\hline
34& $z^iyf_{3-i}(z,y^2)+x^6$ & $tx+g(z,x,y) \subset
\PP(35,12,7,6)$  & 12
\\
& &
$\Diff$=\scriptsize{$(\frac12,0,0,0)$} &
\\
\hline
35& $z^3y+x^6+x^2y^5$ & $tx+zy+x^6+x^2y^5 \subset
\PP(25,26,5,4)$  & 24
\\
& &
$\Diff$=\scriptsize{$(\frac12,\frac23,0,0)$} &
\\
\hline
36& $z^3y+x^6+zxy^4$ & $tx+g(z,x,y) \subset
\PP(55,19,11,9)$  & 18
\\
& &
$\Diff$=\scriptsize{$(\frac12,0,0,0)$} &
\\
\hline
37& $z^3y+x^6+ty^4$ & $t^2x+zy+x^6+ty^4 \subset
\PP(20,41,8,7)$  & 12
\\
& &
$\Diff$=\scriptsize{$(0,\frac23,0,0)$} &
\\
\hline
38& $z^3y+x^6+tzy^2$ & $ t^2x +g(z,x,y) \subset
\PP(25,17,10,9)$  &9
\\
& &$\Diff=\emptyset$
& \\
\hline
39& $z^iyf_{3-i}(z,y^2)+\underline{a}x^4y^2+\underline{b}zx^4$
& $tx+g(z,x,y) \subset \PP(23,8,5,4)$  & 8
\\
&$i\le 2$ &
$\Diff$=\scriptsize{$(\frac12,0,0,0)$} &
\\
\hline
40& $z^3y+x^{2i}y^2f_{2-i}(x^2,y^3)$ & $tx+g(z^{1/3},x,y) \subset
\PP(13,14,3,2)$  & 12
\\
&$i\le 1$ &
$\Diff$=\scriptsize{$(\frac12,\frac23,0,0)$} &
\\
\hline
41& $z^3y+x^4y^2+y^9$ & $tx+zy+x^4y^2+y^9 \subset
\PP(29,32,7,4)$  & 24
\\
& &
$\Diff$=\scriptsize{$(\frac12,\frac23,0,0)$} &
\\
\hline
42& $z^3y+x^4y^2+zxy^4$ & $tx+g(z,x,y) \subset
\PP(31,11,7,5)$  & 10
\\
& &
$\Diff$=\scriptsize{$(\frac12,0,0,0)$} &
\\
\hline
43& $z^3y+x^4y^2+ty^n$ & $t^2x+g(z^{1/3},x,y) \subset
\PP(3n+2,8n-3,2n-2,5)$  &9,15
\\
&$n=4;5$ &
$\Diff$=\scriptsize{$(0,\frac23,0,0)$} &
\\
\hline
44& $z^3y+x^4y^2+tzy^2$ & $ t^2x +g(z,x,y) \subset
\PP(19,13,8,7)$  &7
\\
& &$\Diff=\emptyset$
& \\
\hline
45& $z^3y+zx^4+y^n$ & $tx+g(z,x,y) \subset
\PP(10n-1,4n-4,2n+1,12)$  &8,24
\\
&$n=6;8$ &
$\Diff$=\scriptsize{$(\frac12,0,0,0)$} &
\\
\hline
46& $z^3y+z^2x^2+y^n$ & $tx+g(z,x,y) \subset
\PP(5n-2,2n-2,n+2,6)$  &6,12
\\
&$n=5;9$ &
$\Diff$=\scriptsize{$(\frac12,0,0,0)$} &
\\
\hline
47& $z^3y+z^2x^2+atzy^2+bzxy^3+$ & $ t^2x +g(z,x,y) \subset
\PP(7,5,4,3)$  &3
\\
&$cx^3y^2+dy^6, d\ne 0 ||( ba\ne 0|| ac\ne 0)$ &$\Diff=\emptyset$
& \\
\hline
48& $z^3y+z^2x^2+azxy^4+bx^3y^3+$ & $tx+g(z,x,y) \subset
\PP(19,7,5,3)$  &6
\\
&$cy^8, |a|+|c|\ne 0$ &
$\Diff$=\scriptsize{$(\frac12,0,0,0)$} &
\\
\hline
49& $z^3y+z^2x^2+azy^{n+1}+bx^2y^n$ & $tx+g(z,x,y) \subset
\PP(5n+2,2n,n+2,4)$  & 4,8
\\
&$n=3\ ab\ne 0;5$ &
$\Diff$=\scriptsize{$(\frac12,0,0,0)$} &
\\
\hline
50& $z^3y+z^2x^2+ty^n$ & $ t^2x +g(z,x,y) \subset
\PP(5n-2,4n-3,2n+2,7)$  &5,7
\\
&$n=4;5$ &$\Diff=\emptyset$
& \\
\hline

\end{supertabular}
}

\end{center}
\footnotesize{
\par (1) There are two cases. The common exceptionality condition is
$|b|+|l|\ne 0$.
A). The singularity is $t^2x+z^3y+z^2x^2+azx^2y^2+bzy^5+cx^4y+dx^2y^4+ly^7$.\\
B). The singularity is $t^2x+z^3y+x^4y+az^2x^2+czx^2y^2+bzy^5+dx^2y^4+ly^7$.
}

\begin{center}


\vspace{1cm}

\begin{center}
\small{17. Singularity -- $t^2y+g(t,z,x,y)$}
\end{center}

\tabletail{\hline}
\tablehead{\hline}
\scriptsize{
\begin{supertabular}{|c|c|c|c|}
\hline 1& $z^4+x^5+az^2y^3+by^6$ & $t^2y+z^4+x+az^2y^3+by^6
\subset \PP(5,3,12,2)$  & 6
\\
& & $\PP(5,3,2)$,\
$\Diff$=\scriptsize{$(0,0,\frac45,0)$} &
\\
\hline
2& $z^4+x^5+y^{2n+1}$ & $t^2y+z+x+y^{2n+1}
\subset \PP(n,2n+1,2n+1,1)$  & 16,20
\\
&$n=3;4$ & $\PP(n,2n+1,1)$,\
$\Diff$=\scriptsize{$(0,\frac34,\frac45,0)$} &
\\
\hline
3& $f_4(z,y^2)+x^5$ & $ty+f_4(z,y^2)+x
\subset \PP(7,2,8,1)$  & 10
\\
& & $\PP(7,2,1)$,\
$\Diff$=\scriptsize{$(\frac12,0,\frac45,0)$} &
\\
\hline
4& $z^4+x^5+xy^n$ & $t^2y+g(z^{1/2},x,y)
\subset \PP(5n-4,5n,2n,8)$  & 12,16
\\
&$n=5;7$ &
$\Diff$=\scriptsize{$(0,\frac12,0,0)$} &
\\
\hline
5& $z^4+x^{2i+1}f_{2-i}(x^2,y^3)$ & $t^2y+g(z^{1/2},x,y)
\subset \PP(13,15,6,4)$  & 8
\\
& &
$\Diff$=\scriptsize{$(0,\frac12,0,0)$} &
\\
\hline
6& $z^4+x^5+az^2xy^2+bx^2y^4+$ & $ty+g(z,x,y)
\subset \PP(17,5,4,3)$  & 6
\\
&$czy^5+dzx^3y, |a|+|b|+|c|\ne 0$ &
$\Diff$=\scriptsize{$(\frac12,0,0,0)$} &
\\
\hline
7& $z^4+x^5+x^2y^5$ & $t^2y+z+x^5+x^2y^5
\subset \PP(11,25,5,3)$  & 12
\\
& &$\PP(11,5,3)$,\
$\Diff$=\scriptsize{$(0,\frac34,0,0)$} &
\\
\hline
8& $z^4+x^5+zy^n$ & $ty+z^4+x+zy^n
\subset \PP(4n-3,n,4n,3)$  & 10,30
\\
&$n=4;7$ & $\PP(4n-3,n,3)$,\
$\Diff$=\scriptsize{$(\frac12,0,\frac45,0)$} &
\\
\hline
9& $z^4+x^5+zxy^n$ & $ty+g(z,x,y)
\subset \PP(20n-11,5n,4n,11)$  & 6,18,22
\\
&$n=3;4;5$ &
$\Diff$=\scriptsize{$(\frac12,0,0,0)$} &
\\
\hline
10& $z^4+x^5+zx^2y^n$ & $ty+g(z,x,y)
\subset \PP(20n-7,5n,4n,7)$  & 10,14
\\
&$n=2;3$ &
$\Diff$=\scriptsize{$(\frac12,0,0,0)$} &
\\
\hline
11& $z^4+zx^4+az^2y^3+by^6$ & $t^2y+z^4+zx+az^2y^3+by^6
\subset \PP(5,3,9,2)$  & 8
\\
& &
$\Diff$=\scriptsize{$(0,0,\frac34,0)$} &
\\
\hline
12& $z^4+zx^4+y^7$ & $t^2y+z^4+zx+y^7
\subset \PP(12,7,21,4)$  & 16
\\
& &
$\Diff$=\scriptsize{$(0,0,\frac34,0)$} &
\\
\hline
13& $z^4+zx^4+xy^n$ & $ty+g(z,x,y)
\subset \PP(16n-13,4n,3n,13)$  & 10,26
\\
&$n=5;6$ &
$\Diff$=\scriptsize{$(\frac12,0,0,0)$} &
\\
\hline
14& $z^4+zx^4+x^3y^3$ & $ty+g(z,x,y)
\subset \PP(41,12,9,7)$  & 14
\\
& &
$\Diff$=\scriptsize{$(\frac12,0,0,0)$} &
\\
\hline
15& $z^3x+x^5+azx^2y^n+by^{3n}$ & $ty+g(z,x,y)
\subset \PP(15n-5,4n,3n,5)$  & 4,10
\\
&$n=2\ b\ne 0;4$ &
$\Diff$=\scriptsize{$(\frac12,0,0,0)$} &
\\
\hline
16& $z^3x+x^5+y^{2n+1}$ & $t^2y+g(z^{1/3},x,y)
\subset \PP(5n,8n+4,2n+1,5)$  & 7,15
\\
&$n=3;5$ &
$\Diff$=\scriptsize{$(0,\frac23,0,0)$} &
\\
\hline
17& $z^3x+x^5+y^n$ & $ty+g(z^{1/3},x,y)
\subset \PP(5n-5,4n,n,5)$  & 15,30
\\
&$n=8;14$ &
$\Diff$=\scriptsize{$(\frac12,\frac23,0,0)$} &
\\
\hline
18& $z^3x+x^if_{5-i}(x,y^2)+az^3y^2$ & $ty+g(z^{1/3},x,y)
\subset \PP(9,8,2,1)$  & 6
\\
&$a\ne 0 || i\le 3$ &
$\Diff$=\scriptsize{$(\frac12,\frac23,0,0)$} &
\\
\hline
19& $z^3x+x^5+x^2y^n$ & $ty+g(z^{1/3},x,y)
\subset \PP(5n-3,4n,n,3)$  & 10,18
\\
&$n=4;8$ &
$\Diff$=\scriptsize{$(\frac12,\frac23,0,0)$} &
\\
\hline
20& $z^3x+x^5+x^2y^5$ & $t^2y+zx+x^5+x^2y^5
\subset \PP(11,20,5,3)$  & 9
\\
& &
$\Diff$=\scriptsize{$(0,\frac23,0,0)$} &
\\
\hline
21& $z^3x+x^5+zy^n$ & $ty+g(z,x,y)
\subset \PP(15n-11,4n,3n,11)$  & 7,14,16
\\
&$n=4;6;8;10$ &
$\Diff$=\scriptsize{$(\frac12,0,0,0)$} &
22\\
\hline
22& $z^3x+x^5+zy^7$ & $ t^2y +g(z,x,y)
\subset\PP(47,28,21,11)$  & 11
\\
& &$\Diff=\emptyset$
& \\
\hline
23& $z^3x+x^5+z^2y^3$ & $ t^2y +g(z,x,y)
\subset\PP(19,12,9,7)$  & 6
\\
& &$\Diff=\emptyset$
& \\
\hline
24& $z^3x+x^5+z^2y^n$ & $ty+g(z,x,y)
\subset \PP(15n-7,4n,3n,7)$  & 12,14
\\
&$n=4;6$ &
$\Diff$=\scriptsize{$(\frac12,0,0,0)$} &
\\
\hline
25& $z^3x+x^6+y^7$ & $t^2y+zx+x^6+y^7
\subset \PP(18,35,7,6)$  & 9
\\
& &
$\Diff$=\scriptsize{$(0,\frac23,0,0)$} &
\\
\hline
26& $z^3x+x^6+ax^3y^4+by^8$ & $ty+g(z^{1/3},x,y)
\subset \PP(21,20,4,3)$  & 18
\\
& &
$\Diff$=\scriptsize{$(\frac12,\frac23,0,0)$} &
\\
\hline
27& $z^3x+x^6+x^2y^5$ & $t^2y+zx+x^6+x^2y^5
\subset \PP(13,25,5,4)$  & 12
\\
& &
$\Diff$=\scriptsize{$(0,\frac23,0,0)$} &
\\
\hline
28& $z^3x+x^6+zy^n$ & $ty+g(z,x,y)
\subset \PP(18n-13,5n,3n,13)$  & 10,26
\\
&$n=5;6$ &
$\Diff$=\scriptsize{$(\frac12,0,0,0)$} &
\\
\hline
29& $z^3x+x^6+zx^2y^3$ & $ty+g(z,x,y)
\subset \PP(47,15,9,7)$  & 14
\\
& &
$\Diff$=\scriptsize{$(\frac12,0,0,0)$} &
\\
\hline
30& $z^3x+x^6+z^2y^3$ & $ t^2y +g(z,x,y)
\subset\PP(23,15,9,8)$  & 8
\\
& &$\Diff=\emptyset$
& \\
\hline
31& $z^3x+zx^4+y^8$ & $ty+g(z,x,y)
\subset \PP(77,24,16,11)$  & 16
\\
& &
$\Diff$=\scriptsize{$(\frac12,0,0,0)$} &
\\
\hline
32& $z^3x+zx^4+x^2y^n$ & $ty+g(z,x,y)
\subset \PP(11n-7,3n,2n,7)$  & 8,14
\\
&$n=4;6$ &
$\Diff$=\scriptsize{$(\frac12,0,0,0)$} &
\\
\hline
33& $z^3x+zx^4+x^3y^4$ & $ty+g(z,x,y)
\subset \PP(39,12,8,5)$  & 10
\\
& &
$\Diff$=\scriptsize{$(\frac12,0,0,0)$} &
\\
\hline
34& $z^3x+tx^3+y^n$ & $t^2y+g(z^{1/3},x,y)
\subset \PP(3n-3,5n-1,n+1,6)$  & 9,18
\\
&$n=6;10$ &
$\Diff$=\scriptsize{$(0,\frac23,0,0)$} &
\\
\hline
35& $z^3x+tx^3+ax^3y^3+by^7$ & $t^2y+g(z^{1/3},x,y)
\subset \PP(9,17,4,3)$  & 6
\\
&$b\ne 0$ &
$\Diff$=\scriptsize{$(0,\frac23,0,0)$} &
\\
\hline
36& $z^3x+tx^3+ax^4y^2+bx^2y^5+$ & $t^2y+g(z^{1/3},x,y)
\subset \PP(7,13,3,2)$  & 6
\\
&$cy^8, |b|+|c|\ne 0$ &
$\Diff$=\scriptsize{$(0,\frac23,0,0)$} &
\\
\hline
37& $z^3x+tx^3+az^2y^3+bzx^3y+$ & $ t^2y +g(z,x,y)
\subset\PP(11,7,5,4)$  & 4
\\
&$cx^2y^4, |a|+|c|\ne 0$ &$\Diff=\emptyset$
& \\
\hline
38& $z^3x+tx^3+x^2y^6$ & $t^2y+g(z^{1/3},x,y)
\subset \PP(17,31,7,4)$  & 12
\\
& &
$\Diff$=\scriptsize{$(0,\frac23,0,0)$} &
\\
\hline
39& $z^3x+tx^3+zy^n$ & $ t^2y +g(z,x,y)
\subset\PP(9n-7,5n-1,3n+2,13)$  & 7,10
\\
&$n=4;6$ &$\Diff=\emptyset$
& \\
\hline
40& $z^3x+tx^3+z^2y^4$ & $ t^2y +g(z,x,y)
\subset\PP(31,19,13,8)$  & 8
\\
& &$\Diff=\emptyset$
& \\
\hline
\end{supertabular}
}

\end{center}

\newpage
\hoffset=0cm

\end{document}